\documentclass[11pt,a4paper]{amsart}
\usepackage[T1]{fontenc}
\usepackage[utf8]{inputenc}
\usepackage{lmodern,amsmath,amssymb,amsthm,mathtools,booktabs,array,tabularx,longtable,mathrsfs,needspace}
\usepackage[a4paper,margin=25mm]{geometry}
\usepackage[expansion=false]{microtype}
\usepackage{xcolor}
\usepackage[colorlinks=true,linkcolor=blue!40!black,urlcolor=blue!40!black,citecolor=blue!40!black]{hyperref}
\hypersetup{pdftitle={The Kao-Lai-Osting conjecture for the second Laplace eigenvalue on Riemannian two-tori},
 pdfauthor={Alexander Pigazzini},pdfsubject={Spectral geometry},
 pdfkeywords={Laplace eigenvalues, two-tori, conformal folding, conformal area}}
\newtheorem{theorem}{Theorem}[section]
\newtheorem{lemma}[theorem]{Lemma}
\newtheorem{proposition}[theorem]{Proposition}
\theoremstyle{remark}
\newtheorem{remark}[theorem]{Remark}
\DeclareMathOperator{\Area}{Area}
\newcommand{\Teq}{T_{\mathrm{eq}}}
\newcommand{\Aeq}{A_{\mathrm{eq}}}
\newcommand{\Cstar}{C_*}
\newcommand{\E}{\mathcal E}
\newcommand{\Abar}{\overline A}
\newcommand{\Av}{\mathbb E}
\newcommand{\dd}{\,\mathrm d}
\newcommand{\pos}[1]{\bigl(#1\bigr)_+}
\newcommand{\centerpar}{\boldsymbol{\zeta}}
\newcommand{\othercenter}{\boldsymbol{\eta}}
\newcommand{\weights}{\mathbf r}
\allowdisplaybreaks[2]
\renewcommand{\arraystretch}{1.2}
\title[The second Laplace eigenvalue on two-tori]{The Kao--Lai--Osting conjecture for the second Laplace eigenvalue on Riemannian two-tori}
\author{Alexander Pigazzini}
\date{}
\subjclass[2020]{Primary 58J50; Secondary 35P15, 53C42}
\keywords{Laplace--Beltrami eigenvalues, two-tori, spectral geometry, conformal folding, conformal area}
\begin{document}
\begin{abstract}
We prove that the supremum of the second positive Laplace eigenvalue,
normalized by area, over all smooth Riemannian metrics on a two-torus is
\[
 \frac{8\pi^2}{\sqrt3}+8\pi.
\]
No smooth metric attains this value. The supremum is approached in the
equilateral conformal class by metrics developing one round spherical
bubble, and this is the only conformal class with the same conformal
supremum. The proof combines spherical folding with a sharp inequality
for the truncated conformal area of fundamental embeddings into \(S^5\).
Normal-curvature and normal-chord estimates yield radial bounds for
height caps. For the associated product distributions, we estimate both
the full conformal integral and the part removed by truncation.
Convexity, moment inequalities, and explicit polynomial identities with positive coefficients
give estimates throughout the conformal parameter domain. The sharp value
is uniformly separated from the conformal spectral supremum outside
every neighbourhood of the equilateral class in moduli space.
\end{abstract}
\maketitle
\begingroup
\small
\tableofcontents
\endgroup
\medskip

\section{Introduction}\label{sec:global-spectral}

Let \((\Sigma,g)\) be a closed connected Riemannian surface.
We enumerate the eigenvalues of the nonnegative Laplace--Beltrami
operator with multiplicity, beginning with zero:
\[
 0=\lambda_0(\Sigma,g)<\lambda_1(\Sigma,g)
 \le\lambda_2(\Sigma,g)\le\cdots.
\]
For an integer \(k\ge1\) and a conformal class \(\mathcal C\), define
\[
 \Lambda_k(\Sigma,\mathcal C)
 =\sup_{g\in\mathcal C}\lambda_k(\Sigma,g)\Area(\Sigma,g),
 \qquad
 \Lambda_k(\Sigma)=\sup_{\mathcal C}\Lambda_k(\Sigma,\mathcal C).
\]
All metrics in these definitions are smooth and positive.
The products $\lambda_k(\Sigma,g)\Area(\Sigma,g)$ are invariant under
homothetic rescaling. Their optimization includes classical problems
studied by Berger~\cite{Berger1973}; the distinction between the conformal
and topological spectra was developed by
Colbois and El Soufi~\cite{ColboisElSoufi}.

Exact values of the first topological eigenvalue are known for a small
number of closed surfaces. Hersch~\cite{Hersch} proved
\(\Lambda_1(S^2)=8\pi\), and Li and Yau~\cite{LiYau} proved
\(\Lambda_1(\mathbb{RP}^2)=12\pi\); in both cases the round metric is
maximizing. Nadirashvili~\cite{Nadirashvili} settled Berger's torus
problem by proving
\(\Lambda_1(T^2)=8\pi^2/\sqrt3\), attained by the flat equilateral
metric. On the Klein bottle \(\mathbb K\), the work of Jakobson,
Nadirashvili and Polterovich~\cite{JakobsonNadirashviliPolterovich}
and El Soufi, Giacomini and Jazar~\cite{ElSoufiGiacominiJazar} identifies
the maximizing metric and the value
\[
 \Lambda_1(\mathbb K)=12\pi\,\mathrm E\!\left(\frac{2\sqrt2}{3}\right),
 \qquad
 \mathrm E(\kappa)=\int_0^{\pi/2}
       \sqrt{1-\kappa^2\sin^2\theta}\,d\theta.
\]
Cianci, Karpukhin and Medvedev~\cite{CianciKarpukhinMedvedev}
address the branched-immersion issues in the complete proofs for the
torus and Klein bottle. For the closed orientable surface \(\Sigma_2\)
of genus two, Nayatani and Shoda~\cite{NayataniShoda} established
\(\Lambda_1(\Sigma_2)=16\pi\). A maximizing singular metric is obtained by pulling back the unit
round sphere metric under the branched double covering of $S^2$ by the Bolza surface; the resulting value is also the supremum over smooth metrics.

For higher eigenvalues, concentration of area changes the form of the
extremal configuration. Nadirashvili~\cite{Nadirashvili02}, with a
subsequent proof by Petrides~\cite{Petrides14}, determined
\(\Lambda_2(S^2)=16\pi\), and Nadirashvili and
Sire~\cite{NadirashviliSire} proved \(\Lambda_3(S^2)=24\pi\).
Karpukhin, Nadirashvili, Penskoi and Polterovich~\cite{KNPPsphere}
then obtained \(\Lambda_k(S^2)=8\pi k\) for every \(k\ge1\).
On the real projective plane, Nadirashvili and
Penskoi~\cite{NadirashviliPenskoi} proved
\(\Lambda_2(\mathbb{RP}^2)=20\pi\), and
Karpukhin~\cite{Karpukhin2021} established
\(\Lambda_k(\mathbb{RP}^2)=4\pi(2k+1)\) for all \(k\ge1\).
For \(k\ge2\), neither supremum is attained by a smooth metric.
The respective limiting configurations consist of \(k\) round spheres,
or one round projective plane and \(k-1\) round spheres.

Based on their numerical investigation of conformal and topological
spectral optimization, Kao, Lai and Osting~\cite{KLO} conjectured
\[
 \Lambda_k(T^2)=\frac{8\pi^2}{\sqrt3}+8\pi(k-1),
 \qquad k\ge1.
\]
The lower bound follows from the spherical gluing construction of
Colbois and El Soufi~\cite{ColboisElSoufi}. It is approached by metrics
degenerating to an equilateral flat torus and \(k-1\) round spherical
components, scaled so that their first positive eigenvalues agree.
The remaining question is the matching upper bound over all metrics
and all conformal classes.

For the second eigenvalue, the conformal-volume estimates of Karpukhin
and Stern~\cite{KS} and Eddaoudi and Girouard~\cite{EG} give
\[
 \Lambda_2(\Sigma,\mathcal C)\le4V_c(\Sigma,\mathcal C),
\]
where \(V_c\) denotes conformal volume. We use the smooth-immersion
convention, defined as follows. Writing
$g_{S^n}$ for the unit round metric and $\operatorname{Conf}(S^n)$ for
its conformal automorphism group, we use
\[
 V_c(n,\Sigma,\mathcal C)=
 \inf_{F}\sup_{M\in\operatorname{Conf}(S^n)}
 \Area\bigl(\Sigma,(M\circ F)^*g_{S^n}\bigr),
\]
where the infimum is taken over smooth conformal immersions
$F:(\Sigma,\mathcal C)\to S^n$, with value $+\infty$ if the family is
empty, and $V_c(\Sigma,\mathcal C)=\inf_n V_c(n,\Sigma,\mathcal C)$.
Li and Yau's definition~\cite{LiYau} allows branched conformal
immersions. Karpukhin and Stern~\cite[Section~2.4 and Propositions~3.3, 4.3]{KS}
also allow branched conformal immersions in the infima used in their
proofs. Allowing branch points can only decrease the infimum, so their
second-eigenvalue bound holds with the smooth-immersion convention above.
We recall this inequality for comparison; the proof of the main theorem
uses the folding construction directly.
Kang~\cite{Kang} obtained bounds depending on the torus modulus,
including
\[
 \lambda_2(T^2,g)\Area(T^2,g)<\frac{16\pi^2}{\sqrt3}
\]
for every smooth metric, and reduced the conjectured upper bound to a
compact region~\eqref{eq:K} of conformal moduli.

Theorem~\ref{thm:global-spectral} proves the Kao--Lai--Osting conjecture for \(k=2\),
together with nonattainment by smooth metrics, uniqueness of the
maximizing conformal class and uniform separation in moduli space. To our knowledge, this gives the first
exact determination of a higher topological Laplace eigenvalue on a
closed surface other than \(S^2\) and \(\mathbb{RP}^2\).

Set
\[
 \Aeq=\frac{4\pi^2}{\sqrt3},\qquad
 L=\Aeq-4\pi,\qquad \Cstar=2\Aeq+8\pi.
\]
Thus $\Cstar=8\pi^2/\sqrt3+8\pi$ is the conjectured value for \(k=2\).
Let $\Teq=\mathbb R^2/\langle(1,0),(\tfrac12,\tfrac{\sqrt3}2)\rangle$ be the
equilateral torus, and let $g_{\mathrm{eq}}$ denote its Euclidean
quotient metric, regarded as a metric on $T^2$. Its conformal class $[g_{\mathrm{eq}}]$ is
the equilateral conformal class; in the notation of
Section~\ref{sec:global-geometry} it has modulus
${a_{\mathrm{mod}}}+i {b_{\mathrm{mod}}}=\tfrac12+\tfrac{\sqrt3}2\,i$ and weights $r_1=r_2=r_3=\tfrac13$, and
$\Aeq$ is the area of its fundamental embedding.
Statements of uniqueness of a conformal class are understood in the
moduli space: classes are identified under pullback by diffeomorphisms.
No uniqueness of a marked conformal class on a fixed parametrized torus
is asserted.

\begin{theorem}\label{thm:global-spectral}
For every smooth metric $g$ on the torus,
\[
 \lambda_2(T^2,g)\Area(T^2,g)<\Cstar=\frac{8\pi^2}{\sqrt3}+8\pi.
\]
Nevertheless $\Lambda_2(T^2)=\Cstar$. Moreover,
\[
 \Lambda_2(T^2,\mathcal C)=\Cstar
 \quad\Longleftrightarrow\quad
 \mathcal C=[g_{\mathrm{eq}}].
\]
The supremum is approached by a sequence of metrics in the equilateral
conformal class that develops one round spherical bubble. For every open neighbourhood $\mathcal V$ of $[g_{\mathrm{eq}}]$ in moduli space, there exists
$\varepsilon_{\mathcal V}>0$ such that
\[
 \Lambda_2(T^2,\mathcal C)\le\Cstar-\varepsilon_{\mathcal V}
 \qquad(\mathcal C\notin\mathcal V).
\]
In particular, for every sequence of smooth metrics $g_j$ satisfying
$\lambda_2(T^2,g_j)\Area(T^2,g_j)\to\Cstar$, the classes $[g_j]$
converge to $[g_{\mathrm{eq}}]$ in moduli space.
\end{theorem}

Theorem~\ref{thm:global-spectral} is proved in Section~\ref{sec:sharpness}.
The conformal-volume method of Li--Yau \cite{LiYau}, the conformal-area
theory of Montiel--Ros \cite{MontielRos}, and Bryant's explicit torus
calculations \cite{Bryant} provide the geometric framework.
Folding supplies test functions orthogonal both to constants and to
a first eigenfunction \cite{Nadirashvili02,GNP,EG,KS,Kang}.
Kang \cite{Kang} reduces the unresolved
moduli to
\begin{equation}\label{eq:K}
 \mathcal K=\{0\le {a_{\mathrm{mod}}}\le1/2,\quad\sqrt{1-{a_{\mathrm{mod}}}^2}\le {b_{\mathrm{mod}}}\le44/25\}.
\end{equation}
Section~\ref{sec:exterior} treats the complementary moduli and proves
the strict comparison at the endpoint \({b_{\mathrm{mod}}}=44/25\).

The proof uses the standard fundamental conformal embedding $F_{\weights}$ into
$S^5$, described in Section~\ref{sec:global-geometry}. Its weights are \(\weights=(r_1,r_2,r_3)\), ordered so that
\(r_1\ge r_2\ge r_3\ge0\) and \(r_1+r_2+r_3=1\); we write \(R=r_1\).
Its induced area element is \(dA_{\weights}\), or \(dA\) when the embedding is
fixed, and its area is \(A(\weights)\).
For ${\centerpar}$ in the open unit ball $B^6$ of $\mathbb R^6$ and $y\in S^5$, define
\[
 \chi_{\centerpar}(y)=\frac{1-|{\centerpar}|^2}{1+|{\centerpar}|^2+2{\centerpar}\cdot y},\qquad
 J(\weights,{\centerpar})=\int_{T^2}\left(\chi_{\centerpar}(F_{\weights})^2-\frac{L}{A(\weights)}\right)_+\,dA_{\weights}.
\]
Here $x_+=\max\{x,0\}$. Orthogonal transformations of the target do
not change the area formula; every conformal transformation of $S^5$
has a conformal factor of this form after such a transformation.

\begin{theorem}[Truncated conformal area]\label{thm:global-truncation}
For every fundamental embedding with modulus in $\mathcal K$ and every
${\centerpar}\in B^6$,
\begin{equation}\label{eq:TR-main}
 J(\weights,{\centerpar})\le4\pi.
\end{equation}
Equality occurs precisely when $r_1=r_2=r_3=1/3$ and ${\centerpar}=0$.
\end{theorem}

Theorem~\ref{thm:global-truncation} is proved in
Section~\ref{sec:analytic-closure}, using the geometric reductions of
Sections~\ref{sec:global-geometry} and~\ref{sec:global-majorants} and the
estimates established in Appendices~\ref{app:small-moderate}
and~\ref{app:large-centres}.

\paragraph{The role of the truncation.}
The quantity $J$ measures the positive excess of a conformal-area
density above the common threshold $c=L/A(\weights)$.
The choice $L=\Aeq-4\pi$ is calibrated by the untransformed equilateral
embedding: when $\weights=(1/3,1/3,1/3)$,
\[
 J(\weights,0)=\Aeq-L=4\pi.
\]
Its usefulness comes from the two-piece structure of folding.

Fix a smooth metric $g$ with conformal class in $\mathcal K$.
Lemma~\ref{lem:global-fold}, applied in Section~\ref{sec:folding-implication}, supplies a map
$G=(G_1,\ldots,G_6):T^2\to S^5$ whose coordinates are orthogonal
in $L^2(g)$ to constants and to a first eigenfunction. On a measurable
set $D\subset T^2$ and its complement $D^c=T^2\setminus D$, the map
has conformal factors $\chi_{\centerpar}$ and $\chi_{\othercenter}$,
for centres $\centerpar,\othercenter\in B^6$.
For its Dirichlet energy, conformal invariance gives
\begin{align*}
 \E(G)&:=\frac12\int_{T^2}|dG|_g^2\,dA_g\\
 &=\int_D\chi_{\centerpar}(F_{\weights})^2\,dA_{\weights}
   +\int_{D^c}\chi_{\othercenter}(F_{\weights})^2\,dA_{\weights}\\
 &\le L+J(\weights,\centerpar)+J(\weights,\othercenter)\\
 &\le L+8\pi=\Aeq+4\pi.
\end{align*}
Indeed, apply $x\le c+(x-c)_+$ on each piece and then enlarge the
two nonnegative excess integrals to the whole torus.
Since $D$ and $D^c$ partition the source, the constant contribution
is counted once: $c\bigl(\Area(D)+\Area(D^c)\bigr)=cA(\weights)=L$,
where these areas are induced by $F_{\weights}$.
Thus one baseline contribution and two truncated excesses control
the energy of the folded map. The single-map alternative is included
by taking $D=T^2$ and identical centres.

The bound is strict unless both truncated integrals equal $4\pi$.
By Theorem~\ref{thm:global-truncation}, that exceptional case requires
the equilateral embedding and two zero centres; then
$\E(G)=\Aeq<\Aeq+4\pi$ as well.
Summing the coordinate Rayleigh inequalities therefore gives
\[
 \lambda_2(T^2,g)\Area(T^2,g)
 \le2\E(G)<2\Aeq+8\pi=\Cstar.
\]
The bubbling construction and the compactness argument proving
uniqueness of the maximizing class are given in
Section~\ref{sec:sharpness}.

The proof uses the spherical folding lemma and the conformal gluing
inequality, with precise references below. The geometric reductions are
proved in the main text; the scalar estimates are proved in
Appendices~\ref{app:small-moderate} and~\ref{app:large-centres}.
The exact algebraic identities used in those estimates are collected in
Appendix~\ref{app:explicit-arithmetic}. Bounds for transcendental
constants follow from explicit finite series and rigorous remainder
estimates, including the logarithmic bounds in
Section~\ref{sec:fixed-log-evaluations} and Machin's formula for $\pi$.

The main analytic task is the truncated-area estimate in
Theorem~\ref{thm:global-truncation}. Normal-curvature and normal-chord
inequalities yield radial bounds for height caps. Convexity and
Fourier moments reduce the remaining ranges to product distributions,
for which the proof controls both the complete conformal integral
and the part removed by truncation.

\section{Organization and notation}\label{sec:analytic-roadmap}

We write \(\centerpar=-sp\), where \(s=|\centerpar|\in[0,1)\) and
\(p\in S^5\). At \(s=0\), the choice of \(p\) is immaterial.
We use \(J_s(\weights,p):=J(\weights,-sp)\), and abbreviate this quantity to \(J\)
when \(\weights,s,p\) are fixed. The principal scalar parameters are
\[
 R=r_1,\qquad R_*=\frac{1936}{2561},\qquad
 \sigma=\frac{s}{1+s^2},\qquad
 \beta=\frac{2s}{1+s^2},\qquad
 \tau=\frac{1-s^2}{1+s^2}.
\]
Thus \(\beta=2\sigma\) and \(\beta^2+\tau^2=1\).
For product estimates we also write \(u=\tau^2\) and
\(\omega=\beta^2=1-u\).
The area envelope and its cutoff are
\[
 \Abar(R)=2\pi^2\sqrt{(1-R)(1+3R)},\qquad
 c_R=\frac{L}{\Abar(R)}.
\]
The actual cutoff is \(c(\weights)=L/A(\weights)\); inequality~\eqref{eq:global-area-envelope}
gives \(A(\weights)\le\Abar(R)\) and hence \(c(\weights)\ge c_R\).
Write $g_{\weights}=F_{\weights}^*g_{S^5}$ for the induced metric and
\[
 d\mu_{\weights}=\frac{dA_{\weights}}{A(\weights)}.
\]
This is the normalized Haar probability measure on the source angular
torus. Expectation denotes integration against this measure for
fundamental-torus heights, or against the product Haar probability
measure for independent circle variables, as specified in each formula.
Unqualified vector norms and scalar products are Euclidean.
Primes denote derivatives with respect to the displayed scalar
argument; partial derivatives are indicated by subscripts or by
\(\partial\). The principal geometric and analytic symbols have the following meanings:
\begin{center}
\begin{tabularx}{\linewidth}{@{}lX@{}}
\toprule
Symbol & Meaning\\\midrule
$a_{\mathrm{mod}},b_{\mathrm{mod}}$ & Real and imaginary parts of the torus modulus.\\
$\weights,R,R_*$ & Weights of the fundamental embedding, largest weight, and its upper bound on $\mathcal K$.\\
$\centerpar,\othercenter$ & Centres of the two conformal transformations.\\
$p,q,m$ & Unit ambient direction, height $q=p\cdot F_{\weights}$, and its maximum $m$. In product estimates $m$ denotes the maximum of the specified product height.\\
$s,\sigma,\beta,\tau,u,\omega$ & The conformal radius and its scalar reparametrizations above.\\
$A(\weights),\Abar(R),L,c_R$ & Induced area, area envelope, fixed area offset, and envelope cutoff.\\
$g_{\weights},\mu_{\weights}$ & Induced metric and normalized induced area measure.\\
$J,J_s,\mathcal P$ & Truncated conformal area, its radius notation, and the product majorant.\\
$Q_{R,\alpha}$ & Product height, with direction amplitude $\alpha\in[0,1]$.\\
$D_{\mathrm{lat}},d_{\mathrm{lat}}$ & Lattice denominator and auxiliary lattice parameter.\\
$D_{\mathrm{height}},D_{\mathrm{prod}}$ & Positive denominators of the height and product Hermite polynomials.\\
$B_{\mathrm{height}},B_{\mathrm{prod}}$ & The distinct moment envelopes of the height and product Hermite majorants.\\
\bottomrule
\end{tabularx}
\end{center}
The symbols in the preceding table have the indicated meanings in the
main geometric reductions. Auxiliary functions, polynomials, coefficients,
interpolation nodes, interval endpoints, and integration variables are
local to the proof or numbered subsection in which they are defined,
unless explicitly identified with an earlier quantity. In the appendices,
each numbered subsection is a separate proof unit, including its
preparatory estimates and any auxiliary lemmas. Local notation in
such a unit is not inherited from the preceding unit except by an explicit
reference. The three angular
coordinates of $F_{\weights}$ are taken modulo $2\pi$. For product
estimates the squared amplitude is $x=\alpha^2$, and the angles $U,V$
are independent with the uniform distribution on $[0,2\pi)$.
Write $\chi_s(q)=(1-s^2)/(1+s^2-2sq)$ for the conformal kernel.
For $s>0$ and $c>0$, the formal truncation height is
\begin{equation}\label{eq:truncation-height}
 q_0(s,c)=\frac{1-\sqrt{u/c}}{\beta}
          =\frac{1+s^2-(1-s^2)/\sqrt c}{2s}.
\end{equation}
It is the solution on the increasing branch of
$\chi_s(q)^2=c$; on $[-1,1]$ the positive part is supported on
$q>q_0(s,c)$, with the evident empty or full-support cases if this
height lies outside $[-1,1]$. At $s=0$ the kernel is constant.
When $c$ is fixed, we write $q_0$ or $q_0(s)$.

Table~\ref{tab:parameter-cover} records the numbered estimates used
in the proof of Theorem~\ref{thm:global-truncation}.
The largest weight $R$ lies in $[1/3,R_*]$, and all joining endpoints are included.
Product estimates imply the corresponding bound for $J$ by
Proposition~\ref{prop:global-product}.

\begingroup\small
\setlength{\tabcolsep}{4pt}
\begin{longtable}{@{}p{.24\linewidth}p{.28\linewidth}p{.42\linewidth}@{}}
\caption{Parameter ranges for the truncated-area bound.}
\label{tab:parameter-cover}\\
\toprule
Radius \(s\)&Weight \(R\)&Numbered estimate\\
\midrule
\endfirsthead
\toprule
Radius \(s\)&Weight \(R\)&Numbered estimate\\
\midrule
\endhead
\bottomrule
\endfoot
\([0,1/4]\)&\([1/3,R_*]\)&Theorem~\ref{thm:small}\\
\([1/4,3/10]\)&\([1/3,R_*]\)&Theorem~\ref{bs:theorem}\\
\midrule
\([3/10,2/5]\)&\([1/3,7/20]\)&Theorem~\ref{thm:quartic-small-weights}\\
\([3/10,2/5]\)&\([7/20,9/20]\)&Theorem~\ref{thm:central-product-strip}\\
\([3/10,2/5]\)&\([9/20,R_*]\)&Theorem~\ref{thm:whole-low-radius-product}\\
\midrule
\([2/5,17/40]\)&\([1/3,7/20]\)&Theorem~\ref{thm:quartic-small-weights}\\
\([17/40,1/2]\)&\([1/3,7/20]\)&Theorem~\ref{mid:strip}\\
\([2/5,1/2]\)&\([7/20,9/20]\)&Theorem~\ref{thm:octic-product-strip}\\
\([2/5,1/2]\)&\([9/20,R_*]\)&Theorem~\ref{thm:whole-low-radius-product}\\
\midrule
\([1/2,3/5]\)&\([1/3,4/9]\)&Theorem~\ref{thm:fixed-product-mid}\\
\([1/2,3/5]\)&\([4/9,R_*]\)&Theorem~\ref{thm:wide-product-three-fifths}\\
\midrule
\([3/5,13/20]\)&\([1/3,1/2]\)&Theorem~\ref{thm:central-small-weights}\\
\([3/5,13/20]\)&\([1/2,R_*]\)&Theorem~\ref{thm:sixth-moment-product}\\
\([13/20,1)\)&\([1/3,1/2]\)&Theorem~\ref{thm:central-small-weights}\\
\midrule
\([13/20,(1+R)/2]\)&\([1/2,R_*]\)&
Lemma~\ref{lem:product-direction}; Propositions~\ref{prop:inner-product},
\ref{prop:middle-product}, \ref{prop:outer-product}\\
\([(1+R)/2,1)\)&\([1/2,R_*]\)&Theorem~\ref{thm:global-highcap-linear}\\
\end{longtable}
\endgroup

The quartic estimates in Theorem~\ref{thm:quartic-small-weights} meet
the radial estimates in Theorem~\ref{mid:strip} at \(s=17/40\).
Section~\ref{sec:analytic-closure} assembles the parameter ranges and
proves Theorem~\ref{thm:global-truncation}.
The complementary moduli and the sharpness argument then complete
Theorem~\ref{thm:global-spectral}.

\section{The folding implication}
\label{sec:folding-implication}

For a closed spherical cap $C\subset S^n$ let $\iota_C$ be the conformal
reflection of $S^n$ across the sphere $\partial C$, and let
$\operatorname{Fold}_C$ be the map equal to the identity on $C$ and to
$\iota_C$ on $S^n\setminus C$. The reflection $\iota_C$ reverses
orientation. A conformal automorphism $M$ of $S^n$ maps closed caps to
closed caps and satisfies $M\circ\iota_{M^{-1}(C)}=\iota_C\circ M$; hence
\[
 \operatorname{Fold}_C\circ M=M\circ\operatorname{Fold}_{M^{-1}(C)},
\]
so folding and composition with a conformal automorphism may be
performed in either order.

\begin{lemma}[Conformal folding]\label{lem:global-fold}
Let $F:T^2\to S^n$ be a smooth embedding, let $g$ be a smooth metric,
and fix a first eigenfunction $u_1$ of $g$. There is a conformal
automorphism $M$ of $S^n$ such that either $G=M\circ F$, or
$G=M\circ\operatorname{Fold}_C\circ F$ for a closed spherical cap $C$, is
a map in $H^1(T^2,S^n)$ whose coordinate functions are orthogonal in
$L^2(g)$ to both $1$ and $u_1$. In the second case $G$ agrees with
$M\circ F$ on $F^{-1}(C)$ and with $(M\circ\iota_C)\circ F$ on its
complement.
\end{lemma}

\begin{proof}[Proof (by reference)]
This is the folding lemma of Eddaoudi and Girouard
\cite[Lemma~2.8]{EG}, in the form stated and proved in the paper of Kang \cite[Lemma~2.2]{Kang}. The method goes back to
Nadirashvili \cite{Nadirashvili02} and to Girouard, Nadirashvili and
Polterovich \cite{GNP}; see also \cite{Petrides14}.
\end{proof}

The same two-piece construction, a conformal dilation composed with a map that is the
identity on a cap and a conformal reflection on its complement, gives
the bound $\Lambda_2(\Sigma,[g])\le4V_c(\Sigma,[g])$ of Karpukhin and Stern
\cite[Proposition~1.20]{KS}, also valid with the smooth-immersion
convention of Section~\ref{sec:global-spectral}. One first conformally balances the
push-forward of $dA_g$. If its correlation with $u_1$ already vanishes,
the first alternative applies. Otherwise one folds across a spherical
cap, balances again, and a degree argument over the family of caps
yields a cap for which the folded map is orthogonal to $u_1$. On each
piece, the resulting map is a composition of globally defined conformal
maps; on the complement of $F^{-1}(C)$ the conformal map applied to $F$
is $M\circ\iota_C$, which reverses orientation when $M$ preserves it.
Like every conformal transformation of $S^n$, it is a conformal dilation
followed by an orthogonal transformation, so its conformal factor is again of
the form $\chi_{\othercenter}$ for some ${\othercenter}\in B^{n+1}$; this is what \eqref{eq:two-center-energy} uses. We use
the orthogonality conclusion of Lemma~\ref{lem:global-fold}, together with the two separate
energy integrals; equality of these integrals is neither assumed nor
needed.

Suppose now that $F$ is conformal and write $dA_F$ for its induced area
element. For ${\centerpar}$ in the open unit ball $B^{n+1}$, let $M_{\centerpar}$ denote a
conformal transformation of $S^n$ whose conformal factor is the
function $\chi_{\centerpar}$ defined by the same formula as in Section~\ref{sec:global-spectral}, now for ${\centerpar}\in B^{n+1}$ and $y\in S^n$. The conformal
invariance of two-dimensional Dirichlet energy gives
\[
 \E(M_{\centerpar}\circ F):=\frac12\int|d(M_{\centerpar}\circ F)|_g^2dA_g
                 =\int \chi_{\centerpar}(F)^2\,dA_F.
\]
For a folded map, put $D=F^{-1}(C)$ and write $D^c=T^2\setminus D$.
The two piecewise smooth compositions agree on $F^{-1}(\partial C)$.
The weak derivative of each coordinate difference vanishes almost
everywhere on its zero set. Thus the energy identity is valid also
when the inverse image of the cap boundary has positive measure. We obtain
\begin{align}
 \E(G)
 &=\int_D \chi_{\centerpar}(F)^2\,dA_F+\int_{D^c} \chi_{\othercenter}(F)^2\,dA_F\notag\\
 &\le\int\max\{\chi_{\centerpar}(F)^2,\chi_{\othercenter}(F)^2\}\,dA_F.\label{eq:two-center-energy}
\end{align}
For the spectral estimate, take $F=F_{\weights}$ for the conformal class of $g$,
with modulus in $\mathcal K$; thus $n=5$, $dA_F=dA_{\weights}$, and $F_{\weights}$
is conformal with respect to $g$. The elementary inequality
\[
 \max\{x,y\}\le c+(x-c)_++(y-c)_+
\]
with $c=L/A(\weights)$ therefore combines with Theorem~\ref{thm:global-truncation} to yield
\begin{equation}\label{eq:two-center-bound}
 \E(G)\le L+8\pi=\Aeq+4\pi.
\end{equation}
If either centre has strict inequality in \eqref{eq:TR-main}, this
bound is strict. The only remaining possibility is the equilateral
embedding with both centres zero, in which case $\E(G)=\Aeq$ and
strictness again holds. The single-map case follows by using two
identical centres.

Writing $G=(G_1,\ldots,G_{n+1})$, the Rayleigh principle gives
\[
 \lambda_2(g)\int G_i^2\,dA_g\le\int|dG_i|_g^2\,dA_g.
\]
Summing and using $\sum G_i^2=1$ proves
\begin{equation}\label{eq:rayleigh-sum}
 \lambda_2(g)\Area(g)\le2\E(G)<2\Aeq+8\pi=\Cstar
\end{equation}
for all moduli in $\mathcal K$. Orthogonality to one chosen first
eigenfunction is sufficient even when $\lambda_1$ is multiple, because
the eigenvalues are ordered with multiplicity.

\section{The fundamental embeddings and a one-parameter reduction}
\label{sec:global-geometry}

Let $0\le {a_{\mathrm{mod}}}\le\frac12$ and ${b_{\mathrm{mod}}}\ge\sqrt{1-{a_{\mathrm{mod}}}^2}$. Put
\[
 \begin{gathered} d_{\mathrm{lat}}=a_{\mathrm{mod}}(1-a_{\mathrm{mod}}),\qquad D_{\mathrm{lat}}=1+b_{\mathrm{mod}}^2-d_{\mathrm{lat}},\\
 r_1=\frac{D_{\mathrm{lat}}-1}{D_{\mathrm{lat}}},\quad r_2=\frac{1-a_{\mathrm{mod}}}{D_{\mathrm{lat}}},\quad
 r_3=\frac{a_{\mathrm{mod}}}{D_{\mathrm{lat}}}.\end{gathered}
\]
The weights satisfy $r_1\ge r_2\ge r_3\ge0$ and $\sum r_i=1$.
The standard fundamental immersion \cite{Bryant}, up to conjugation of
one complex coordinate, is
\[
 F_{\weights}(\theta_1,\theta_2)
  =\bigl(\sqrt{r_1}e^{i\theta_1},\sqrt{r_2}e^{i\theta_2},
          \sqrt{r_3}e^{-i(\theta_1+\theta_2)}\bigr)\in S^5.
\]
To identify the source modulus, take
$\theta_1=2\pi y/{b_{\mathrm{mod}}}$ and $\theta_2=2\pi(x-{a_{\mathrm{mod}}} y/{b_{\mathrm{mod}}})$ on
$\mathbb R^2/\langle(1,0),({a_{\mathrm{mod}}},{b_{\mathrm{mod}}})\rangle$.
The induced metric is $g_{\weights}=(4\pi^2/{D_{\mathrm{lat}}})(dx^2+dy^2)$, so this is a
homothetic conformal embedding of the torus of modulus ${a_{\mathrm{mod}}}+i {b_{\mathrm{mod}}}$.
We call $F_{\weights}$ the fundamental embedding with weights $\weights$, or of modulus
${a_{\mathrm{mod}}}+i {b_{\mathrm{mod}}}$, and the torus embedded in this way a fundamental torus.
To see the metric identity explicitly, abbreviate the modulus by $a+ib$
and $D=1+b^2-a(1-a)$ in this calculation. The three angular differentials
are $(2\pi/b)dy$, $2\pi(dx-a\,dy/b)$ and
$-2\pi(dx+(1-a)dy/b)$. Therefore
\begin{align*}
 \frac{F_{\weights}^*g_{S^5}}{4\pi^2}
 &=r_1\frac{dy^2}{b^2}+r_2(dx-a\,dy/b)^2
       +r_3(dx+(1-a)dy/b)^2\\
 &=\frac{dx^2+dy^2}{D},
\end{align*}
since $r_2+r_3=1/D$, $-ar_2+(1-a)r_3=0$, and
$r_1+a^2r_2+(1-a)^2r_3=b^2/D$. The first two weights are positive;
the first two complex coordinates determine $(\theta_1,\theta_2)$ modulo
$2\pi$, which proves embeddedness, also when $r_3=0$.
Its induced metric is flat, and its area is
\[
 A(\weights)=4\pi^2\sqrt{r_1r_2+r_1r_3+r_2r_3}
     =\frac{4\pi^2{b_{\mathrm{mod}}}}{{D_{\mathrm{lat}}}}.
\]
These statements continue to hold when $r_3=0$: the image is then an
embedded product torus in an equatorial $S^3$.

Write $R=r_1$. For fixed $R$, the elementary inequality
$r_2r_3\le(1-R)^2/4$ gives
\begin{equation}\label{eq:global-area-envelope}
 A(\weights)\le \Abar(R):=2\pi^2\sqrt{(1-R)(1+3R)}\le\Aeq.
\end{equation}
On the residual moduli region ${b_{\mathrm{mod}}}\le44/25$ we have
\begin{equation}\label{eq:global-R-range}
 \frac13\le R\le R_*:=\frac{1936}{2561}.
\end{equation}
Indeed ${D_{\mathrm{lat}}}\le1+{b_{\mathrm{mod}}}^2$ and $R=1-{D_{\mathrm{lat}}}^{-1}$. Equality in the upper
bound is attained at ${a_{\mathrm{mod}}}=0$, ${b_{\mathrm{mod}}}=44/25$.

For a unit ambient vector $p$, set $q=p\cdot F_{\weights}$. Fourier orthogonality
and the constant induced area element imply
\begin{equation}\label{eq:global-moments}
 \int q\,dA=0,\qquad
 \frac1{A(\weights)}\int q^2\,dA
   =\frac12\sum_{i=1}^3r_i|p_i|^2\le\frac R2,
\end{equation}
where $p_i$ are the three complex coordinate components of $p$.

\subsection{Exact normal curvature}

\begin{lemma}\label{lem:global-curvature}
If $\mathrm{II}$ denotes the spherical second fundamental form of
$F_{\weights}$, then
\[
 \kappa_{\weights}:=\max_{|X|=1}|\mathrm{II}(X,X)|
       =\sqrt{{D_{\mathrm{lat}}}-1}=\sqrt{\frac R{1-R}}.
\]
\end{lemma}

\begin{proof}
In the flat physical coordinates of the torus, the three directional
frequencies for an induced-unit tangent vector making angle $\vartheta$
are
\begin{align*}
\ell_1&=\sqrt {D_{\mathrm{lat}}}\,\frac{\sin\vartheta}{{b_{\mathrm{mod}}}},\\
\ell_2&=\sqrt {D_{\mathrm{lat}}}\left(\cos\vartheta-
                         \frac{{a_{\mathrm{mod}}}\sin\vartheta}{{b_{\mathrm{mod}}}}\right),\quad
\ell_3=\sqrt {D_{\mathrm{lat}}}\left(\cos\vartheta+
                         \frac{(1-{a_{\mathrm{mod}}})\sin\vartheta}{{b_{\mathrm{mod}}}}\right).
\end{align*}
Their weighted second moment is one. Direct differentiation therefore gives
$|\mathrm{II}(X,X)|^2=\sum r_i\ell_i^4-1$.
Writing, in this proof only, $c_\vartheta=\cos\vartheta$ and $s_\vartheta=\sin\vartheta$,
expansion yields
\[
 \sum r_i\ell_i^4
  ={D_{\mathrm{lat}}}\left[c_\vartheta^4+\frac{6{d_{\mathrm{lat}}}}{{b_{\mathrm{mod}}}^2}c_\vartheta^2s_\vartheta^2
       +\frac{4{d_{\mathrm{lat}}}(1-2{a_{\mathrm{mod}}})}{{b_{\mathrm{mod}}}^3}c_\vartheta s_\vartheta^3
       +\frac{{b_{\mathrm{mod}}}^2-3{d_{\mathrm{lat}}}^2}{{b_{\mathrm{mod}}}^4}s_\vartheta^4\right].
\]
Subtracting the bracket from $(c_\vartheta^2+s_\vartheta^2)^2$ gives
\begin{align*}
 &s_\vartheta^2(A_0c_\vartheta^2-B_0c_\vartheta s_\vartheta+C_0s_\vartheta^2),\\
 &A_0=2-\frac{6{d_{\mathrm{lat}}}}{{b_{\mathrm{mod}}}^2},\qquad B_0=\frac{4{d_{\mathrm{lat}}}(1-2{a_{\mathrm{mod}}})}{{b_{\mathrm{mod}}}^3},\qquad
 C_0=1-\frac1{{b_{\mathrm{mod}}}^2}+\frac{3{d_{\mathrm{lat}}}^2}{{b_{\mathrm{mod}}}^4}.
\end{align*}
Put ${u_{\mathrm{lat}}}={b_{\mathrm{mod}}}^2-(1-{a_{\mathrm{mod}}}^2)\ge0$. Then $A_0\ge0$, since
${b_{\mathrm{mod}}}^2-3{d_{\mathrm{lat}}}={u_{\mathrm{lat}}}+(1-{a_{\mathrm{mod}}})(1-2{a_{\mathrm{mod}}})$. The numerator of $C_0$ equals
$2{a_{\mathrm{mod}}}^2(1-{a_{\mathrm{mod}}})(1-2{a_{\mathrm{mod}}})$ at ${u_{\mathrm{lat}}}=0$ and increases with ${b_{\mathrm{mod}}}^2\ge3/4$.
Finally,
\[
 \frac{{b_{\mathrm{mod}}}^6}{2}\left(A_0C_0-\frac{B_0^2}{4}\right)
       ={u_{\mathrm{lat}}}({u_{\mathrm{lat}}}+1-2{a_{\mathrm{mod}}})({u_{\mathrm{lat}}}+1-{a_{\mathrm{mod}}})\ge0.
\]
The quadratic form is nonnegative, proving $\sum r_i\ell_i^4\le {D_{\mathrm{lat}}}$.
Equality holds for $\sin\vartheta=0$.
\end{proof}

\subsection{A global normal-chord estimate}

\begin{lemma}\label{lem:global-chord}
For any two points $F_0,F$ of the image,
\[
 |\operatorname{proj}_{N_{F_0}}F|
       \le\kappa_{\weights}(1-F_0\cdot F),
\]
where $N_{F_0}$ is the normal space in the sphere. The spherical normal
injectivity radius is at least $\arctan(1/\kappa_{\weights})$.
\end{lemma}

\begin{proof}
By homogeneity take $F_0=F_{\weights}(0)$. Initially suppose all weights are
positive. Write $\theta_1+\theta_2+\theta_3=0$ and set
\[
 x_i=1-\cos\theta_i,\quad z=\sum r_ix_i=1-F_0\cdot F,
 \quad P=\sum\sin\theta_i,\quad S=\sum r_i^{-1}.
\]
The trigonometric identities
\[
 P^2=(x_1+x_2+x_3)^2-2\sum x_i^2=2x_1x_2x_3
\]
hold because $\theta_1+\theta_2+\theta_3=0$. Identify $\mathbb C^3$ with
$\mathbb R^6$. At $F_0=(\sqrt{r_1},\sqrt{r_2},\sqrt{r_3})$ the tangent
vectors are $i(\sqrt{r_j}\,a_j)_j$ with $\sum a_j=0$, so the spherical
normal space is
\[
 N_{F_0}=\Bigl\{\bigl(\xi_j+i\lambda r_j^{-1/2}\bigr)_j:\
 \sum_j\sqrt{r_j}\,\xi_j=0,\ \lambda\in\mathbb R\Bigr\}.
\]
The real part $(\sqrt{r_j}\cos\theta_j)_j$ of $F$ projects to
$(\sqrt{r_j}(z-x_j))_j$, whose squared length is $\sum r_jx_j^2-z^2$.
The imaginary part $(\sqrt{r_j}\sin\theta_j)_j$ projects onto the unit
vector $S^{-1/2}(r_j^{-1/2})_j$ with coefficient $P/\sqrt S$. Hence
\begin{align*}
|\operatorname{proj}_{N_{F_0}}F|^2
   &=\sum r_ix_i^2-z^2+\frac{P^2}{S}=x^{\mathsf T}M_{\weights}x-z^2,\\
M_{\weights}&=\operatorname{diag}(r_i-2/S)+S^{-1}\mathbf1\mathbf1^{\mathsf T}.
\end{align*}
The matrix $M_{\weights}$ is positive definite. Its scaled matrix $SM_{\weights}$ has
positive diagonal entries; its principal two-by-two minors are
$Sr_ir_j/r_k>0$, and its determinant is $4$.

The set
\[
 \{x:\ \weights\cdot x=1,\ (\sum x_i)^2\ge2\sum x_i^2,\ \sum x_i\ge0\}
\]
is a compact elliptical disk, a section of a Lorentz cone. This forward
cone is contained in the nonnegative octant. The convex function
$x\mapsto x^{\mathsf T}M_{\weights}x$ attains its maximum on this compact convex set
at an extreme point, hence on the boundary $(\sum x_i)^2=2\sum x_i^2$.
This boundary equation is equivalent to $4x_1x_2=(x_3-x_1-x_2)^2$. For a
nonnegative solution put $t_1=\sqrt{x_1}$, $t_2=\pm\sqrt{x_2}$ with the
sign of $x_3-x_1-x_2$, and $t_3=-t_1-t_2$; then $x_i=t_i^2$,
$\sum t_i=0$ and $\sum r_it_i^2=1$, so that $X=i(\sqrt{r_j}\,t_j)_j$ is a
unit tangent vector at $F_0$. On the boundary the last two terms of
$x^{\mathsf T}M_{\weights}x$ cancel, and the normal-chord quotient is
\[
 \sum r_it_i^4-1=|\mathrm{II}(X,X)|^2\le\kappa_{\weights}^2.
\]
Homogeneity in $x$ proves the asserted inequality. The case $r_3=0$
follows by continuity; the limiting torus is still embedded.

For a unit spherical normal $n$ at $F_0$, put
$p=\cos\eta\,F_0+\sin\eta\,n$. The estimate just proved gives
\[
 p\cdot F\le\cos\eta-
       (\cos\eta-\kappa_{\weights}\sin\eta)(1-F_0\cdot F).
\]
If $\eta<\arctan(1/\kappa_{\weights})$, $F_0$ is the unique maximizer of this
height function. The differential of the normal exponential map in
tangential directions contains \(\cos\eta\,I-\sin\eta\,A_n\),
where the shape operator satisfies \(\|A_n\|\le\kappa_{\weights}\).
This operator is invertible when \(\tan\eta<1/\kappa_{\weights}\).
There are therefore no focal points in the indicated range.
Together with uniqueness of the normal foot, this proves the
normal-tube assertion.
\end{proof}

\subsection{Area density of high caps}

\begin{lemma}\label{lem:global-cap}
Let $q=p\cdot F_{\weights}$ and $m=\max q$. If
$\sqrt R<t<m$, the superlevel
$\{q>t\}$ is a disk and its area density satisfies
\begin{equation}\label{eq:global-cap-density}
 -\frac{d}{dt}\Area\{q>t\}
 \le\frac{2\pi(1-R)}{t-Rm-\sqrt{R(1-R)(1-m^2)}}.
\end{equation}
\end{lemma}

\begin{proof}
We use the flat induced metric $g_{\weights}$ of $F_{\weights}$, with its gradient $\nabla$
and Hessian $\operatorname{Hess}$, the normal connection $\nabla^\perp$ of
the normal bundle of $F_{\weights}$ in $S^5$, and we write $\kappa=\kappa_{\weights}$. At
each point decompose the fixed vector
\[
 p=q\,F_{\weights}+\nabla q+n,
\]
where $n$ is normal to $F_{\weights}$ in $S^5$. Differentiating along a tangent
vector $X$ and separating tangential and normal components gives
\begin{equation}\label{eq:cap-hess}
 \operatorname{Hess}q(X,Y)=-q\langle X,Y\rangle+\langle\mathrm{II}(X,Y),n\rangle,
 \qquad \nabla^\perp_Xn=-\mathrm{II}(X,\nabla q).
\end{equation}
For a unit normal $\nu$ the symmetric form
$\langle\mathrm{II}(\cdot,\cdot),\nu\rangle$ satisfies
$|\langle\mathrm{II}(X,X),\nu\rangle|\le\kappa|X|^2$, so its eigenvalues
lie in $[-\kappa,\kappa]$; hence $|\mathrm{II}(X,Y)|\le\kappa|X||Y|$.

\emph{Critical points and topology.} Since
$|n|^2=1-q^2-|\nabla q|^2\le1-q^2$, \eqref{eq:cap-hess} gives
\[
 -\operatorname{Hess}q\ge\bigl(q-\kappa\sqrt{1-q^2}\bigr)g_{\weights},
\]
and $q-\kappa\sqrt{1-q^2}>0$ exactly when $q>\sqrt R$,
since $\kappa/\sqrt{1+\kappa^2}=\sqrt R$.
Thus every critical point of $q$ in $\{q>\sqrt R\}$ is a nondegenerate
local maximum. At such a point $x_0$ one has $p=q\,F_{\weights}(x_0)+n$, that is,
$p=\cos\eta\,F_{\weights}(x_0)+\sin\eta\,\nu$ with a unit normal $\nu$ and
$\cos\eta=q>\sqrt R$, so that $\eta<\arctan(1/\kappa)$; by
Lemma~\ref{lem:global-chord}, $x_0$ is the unique global maximum point.
Hence, for $\sqrt R<t<m$, the only critical point of $q$ in
$\{q\ge t\}$ is the maximum point $x_m$, and every value in $[t,m)$ is
regular. The set $\{q\ge t\}$ is connected. Indeed, on each of its
connected components $q$ attains a maximum, which exceeds $t$ because
$t$ is a regular value; a point where it is attained is therefore a
local maximum of $q$, hence a critical point in $\{q>t\}$, hence the
point $x_m$. Near $x_m$ the superlevel sets are small disks by the Morse
lemma, and the flow of $\nabla q/|\nabla q|^2$ carries each level set
$\{q=t'\}$ diffeomorphically onto $\{q=t''\}$ for $t\le t'<t''<m$.
Therefore $D_t=\{q\ge t\}$ is a closed disk with smooth boundary
$\Gamma_t=\{q=t\}$, and $\{q>t\}$ is its interior.

\emph{The normal component.} Let $x\in\Gamma_t$ and let $\gamma$ be the
integral curve of $\nabla q/|\nabla q|^2$ through $x$. Along $\gamma$ the
function $q$ increases with unit speed, and $\gamma$ tends to $x_m$ as
$q\to m$, because $D_t$ contains no other critical point. By
\eqref{eq:cap-hess} the function $|n|$ is Lipschitz along $\gamma$ and
satisfies, almost everywhere,
\[
 \Bigl|\frac{d|n|}{dq}\Bigr|\le|\nabla^\perp_{\gamma'}n|
 =|\mathrm{II}(\gamma',\nabla q)|\le\kappa|\gamma'||\nabla q|=\kappa,
\]
and $|n|=\sqrt{1-m^2}$ at $x_m$, where $\nabla q=0$. Integration from
$q=t$ to $q=m$ gives $|n|\le\sqrt{1-m^2}+\kappa(m-t)$ at $x$.

\emph{Curvature of the level curves.} Inserting this bound into
\eqref{eq:cap-hess}, we obtain on $\Gamma_t$
\[
 -\operatorname{Hess}q\ge{b_{\mathrm{cap}}}(t)\,g_{\weights},\qquad
 {b_{\mathrm{cap}}}(t)=(1+\kappa^2)t-\kappa^2m-\kappa\sqrt{1-m^2}.
\]
The function $\mu\mapsto\kappa^2\mu+\kappa\sqrt{1-\mu^2}$ is decreasing
for $\mu>\sqrt R$, so its value at $m$ is at most its value at $t$, and
${b_{\mathrm{cap}}}(t)\ge t-\kappa\sqrt{1-t^2}>0$. Let $T$ be a unit tangent vector of
$\Gamma_t$ and $\nu_t=\nabla q/|\nabla q|$ the unit normal of $\Gamma_t$
pointing into $D_t$.
Differentiating $\langle\nabla q,T\rangle=0$ along $\Gamma_t$ shows that
the geodesic curvature of $\Gamma_t$ with respect to $\nu_t$ is
$k_g=-\operatorname{Hess}q(T,T)/|\nabla q|\ge{b_{\mathrm{cap}}}(t)/|\nabla q|$.

\emph{Gauss--Bonnet and coarea.} The metric $g_{\weights}$ is flat, so the
Gauss--Bonnet formula for the disk $D_t$ gives
$\int_{\Gamma_t}k_g\,ds=2\pi$, hence
$\int_{\Gamma_t}|\nabla q|^{-1}ds\le2\pi/{b_{\mathrm{cap}}}(t)$. All values in
$(\sqrt R,m)$ are regular, and the coarea formula gives
\[
 -\frac{d}{dt}\Area\{q>t\}=\int_{\Gamma_t}\frac{ds}{|\nabla q|}
 \le\frac{2\pi}{{b_{\mathrm{cap}}}(t)}.
\]
Finally $1+\kappa^2=1/(1-R)$, $\kappa^2/(1+\kappa^2)=R$ and
$\kappa/(1+\kappa^2)=\sqrt{R(1-R)}$ turn $2\pi/{b_{\mathrm{cap}}}(t)$ into the
right-hand side of \eqref{eq:global-cap-density}.
\end{proof}

\subsection{A radial comparison for height caps}

The normal--chord estimate of Lemma~\ref{lem:global-chord} also controls
critical heights below $\sqrt R$. This gives a radial description of a larger family of caps.

\begin{lemma}[Radial cap comparison]\label{lem:radial-cap-density}
Let $q=p\cdot F_{\weights}$, where $|p|=1$, and suppose that $m=\max q>\sqrt R$.
Put $h=1-R$ and
\begin{equation}\label{radial:critical-height}
 d_R(m)=Rm+\sqrt{R(1-R)}\sqrt{1-m^2},\qquad
 e_R(m)=2d_R(m)-m.
\end{equation}
The maximum point of $q$ is unique, and every other critical point has
height at most $e_R(m)$. For $e_R(m)<t<m$, the cap $D_t=\{q\ge t\}$
is a smooth closed disk. Its lift to the Euclidean universal cover,
centred at its maximum, is star-shaped. With
$z=(t-d_R(m))/(m-d_R(m))\in(-1,1)$, its area density satisfies
\begin{equation}\label{radial:oscillator-density}
 -\frac{\mathrm d}{\mathrm dt}\Area\{q>t\}
 \le\frac{2\pi h}{m-d_R(m)}
          \frac{\arccos z}{\sqrt{1-z^2}}
 <\frac{4\pi h}{t-e_R(m)}.
\end{equation}
Moreover,
\begin{equation}\label{radial:cap-area}
 \Area\{q>t\}\le
 \pi h\left(\arccos\frac{t-d_R(m)}{m-d_R(m)}\right)^2.
\end{equation}
\end{lemma}

\begin{proof}
All gradients and Hessians in this proof are taken with respect to
$g_{\weights}$. Write $\kappa=\sqrt{R/h}$. Lemma~\ref{lem:global-curvature} implies
$|\mathrm{II}(X,Y)|\le\kappa|X||Y|$: for every unit spherical normal
$\nu$, the symmetric form $\langle\mathrm{II}(\cdot,\cdot),\nu\rangle$
has all its eigenvalues in $[-\kappa,\kappa]$. Lemma~\ref{lem:global-chord}
gives
\begin{equation}\label{radial:normal-chord}
 |\operatorname{proj}_{N_{F_{\weights}(x)}}F_{\weights}(y)|
 \le\kappa\bigl(1-F_{\weights}(x)\cdot F_{\weights}(y)\bigr)
 \qquad(x,y\in T^2).
\end{equation}
Set $n_0=\sqrt{1-m^2}$, $a=m-\kappa n_0>0$, $d=d_R(m)$ and
$e=e_R(m)$. Thus $m-d=ha>0$ and $e<m$.

\emph{Critical heights.}
Fix a maximum point $F_0$ and decompose
$p=mF_0+\mathbf n_0$, where $\mathbf n_0$ is spherical normal and
$|\mathbf n_0|=n_0$. For another image point $F_1$, write
$\zeta=F_0\cdot F_1<1$. The strict inequality follows from embeddedness.
By \eqref{radial:normal-chord},
\begin{equation}\label{radial:unique-maximum}
 q(F_1)\le m-a(1-\zeta)<m.
\end{equation}
This proves uniqueness of the maximum. If $F_1$ is critical and its
height is $q_1$, write $p=q_1F_1+\mathbf n_1$, where $\mathbf n_1$
is spherical normal and $|\mathbf n_1|=\sqrt{1-q_1^2}$. Applying
\eqref{radial:normal-chord} at $F_1$ and using
\eqref{radial:unique-maximum} gives
\[
 a(1-\zeta)\le m-q_1
 \le\bigl(\kappa\sqrt{1-q_1^2}-q_1\bigr)(1-\zeta).
\]
Hence $q_1+a\le\kappa\sqrt{1-q_1^2}$. If $q_1+a\le0$, then
$q_1\le-a<e$. Otherwise, squaring and taking the upper root yields
\begin{equation}\label{radial:critical-ceiling}
 q_1\le h\bigl[-a+\kappa\sqrt{1+\kappa^2-a^2}\bigr]
 =h\bigl[(\kappa^2-1)m+2\kappa n_0\bigr]=e.
\end{equation}
Here $1+\kappa^2-a^2=(\kappa m+n_0)^2$; the upper-root expression
also shows $e>-a$.

\emph{Topology and Hessian comparison.}
There are no critical points in $\{e<q<m\}$. Every component of $D_t$
contains a point with height greater than $t$, because $t$ is a regular
value. Its maximum is therefore an interior critical point and must
be the unique global maximum. Thus $D_t$ is connected. At its maximum,
$-\operatorname{Hess}q\ge(m-\kappa n_0)g_{\weights}=ag_{\weights}>0$, so the maximum is
nondegenerate. The Morse lemma and regular gradient flow identify
$D_t$ with a closed disk with smooth boundary.

To obtain the Hessian bound throughout this disk, decompose
$p=qF_{\weights}+\nabla q+\mathbf n$. By \eqref{eq:cap-hess},
\[
 \operatorname{Hess}q(X,Y)
 =-q\langle X,Y\rangle+\langle\mathrm{II}(X,Y),\mathbf n\rangle,
 \qquad
 \nabla_X^\perp\mathbf n=-\mathrm{II}(X,\nabla q).
\]
The flow of $\nabla q/|\nabla q|^2$ from any nonmaximum point of
$D_t$ tends to its unique maximum. Almost everywhere along this flow,
$|\mathrm d|\mathbf n|/\mathrm dq|\le\kappa$. Integrating gives
$|\mathbf n|\le n_0+\kappa(m-q)$, and hence
\begin{equation}\label{radial:hessian-comparison}
 -\operatorname{Hess}q\ge
 \bigl((1+\kappa^2)q-\kappa^2m-\kappa n_0\bigr)g_{\weights}
 =\frac{q-d}{h}\,g_{\weights}.
\end{equation}
The right-hand side is allowed to be negative when $q<d$.

\emph{Radial description.}
Lift $D_t$ isometrically to the Euclidean universal cover of the flat
torus, placing its maximum at the origin. This lift is a compact
embedded disk. For a unit vector $u$, let $r_t(u)$ be the first exit
radius of the ray $ru$ from its interior. Put $w(r)=q(ru)-d$ and
$M=m-d>0$. While the ray remains in the disk,
\eqref{radial:hessian-comparison} gives
\begin{equation}\label{radial:radial-oscillator}
 w''+h^{-1}w\le0,\qquad w(0)=M,\qquad w'(0)=0.
\end{equation}
The derivative is initially strictly negative. On any initial interval
where $w'<0$,
\[
 \frac{\mathrm d}{\mathrm dr}\bigl((w')^2+h^{-1}w^2\bigr)
 =2w'(w''+h^{-1}w)\ge0.
\]
Since $-M<t-d\le w<M$ away from the origin in this disk, it follows that
\begin{equation}\label{radial:radial-speed}
 -w'\ge h^{-1/2}\sqrt{M^2-w^2}>0.
\end{equation}
Consequently $w'$ cannot have a first zero before exit. Integrating
\eqref{radial:radial-speed} yields
\begin{equation}\label{radial:exit-radius}
 r_t(u)\le\sqrt h\arccos\frac{t-d}{m-d}.
\end{equation}
The exit is transverse by \eqref{radial:radial-speed}, and its radius
depends smoothly on $u$. The radial first-exit set is therefore a
compact star-shaped disk with boundary contained in $\{q=t\}$. By
transversality this boundary locally coincides with the regular level
curve bounding the lift of $D_t$. It is consequently an entire connected
component of that boundary. The lift of $D_t$ has only one boundary
component; the first-exit disk is therefore all of $D_t$.

\emph{Area and density.}
Polar integration and \eqref{radial:exit-radius} give
\eqref{radial:cap-area}. Write $u(\theta)=(\cos\theta,\sin\theta)$
and $r_t(\theta)=r_t(u(\theta))$ in polar coordinates.
Polar differentiation, followed by
\eqref{radial:radial-speed}--\eqref{radial:exit-radius}, gives
\[
 -\frac{\mathrm d}{\mathrm dt}\Area\{q>t\}
 =\int_0^{2\pi}\frac{r_t(\theta)}{-q_r(r_t(\theta),\theta)}\dd\theta
 \le\frac{2\pi h}{m-d}\frac{\arccos z}{\sqrt{1-z^2}}.
\]
Here $z=(t-d)/(m-d)$.
For $0<\theta<\pi$ one has $\theta<2\tan(\theta/2)$. With
$\theta=\arccos z$ this implies
\[
 \frac{\arccos z}{\sqrt{1-z^2}}<\frac2{1+z},
 \qquad (m-d)(1+z)=t-e,
\]
which proves the second inequality in \eqref{radial:oscillator-density}.
\end{proof}

\begin{remark}
For $t>\sqrt R$, the first bound in \eqref{radial:oscillator-density}
is strictly smaller than \eqref{eq:global-cap-density}, since
$z\arccos z<\sqrt{1-z^2}$ for $0<z<1$. The radial estimate also
applies throughout the larger interval $e_R(m)<t<m$. For a direction
through an image point, $m=1$ and $e_R(1)=2R-1$.
\end{remark}

\section{Scalar majorants and elementary constants}\label{sec:global-majorants}
Put
\begin{align*}
L&=\Aeq-4\pi,\quad c(\weights)=L/A(\weights),\quad {c_R}=L/\Abar(R),\\
\chi_s(q)&=\frac{1-s^2}{1+s^2-2sq}.
\end{align*}
For ${\centerpar}=-sp$, the desired single-centre estimate is
\[
 J(\weights,{\centerpar}):=\int\bigl(\chi_s(p\cdot F_{\weights})^2-c(\weights)\bigr)_+\,dA\le4\pi.
\]
Set $f_{R,s}(q)=(\chi_s(q)^2-{c_R})_+$. This is convex in $q$.
The inequality $c(\weights)\ge{c_R}$ permits replacing the
integrand by $f_{R,s}$.

Let \(U,V\) be independent uniform angles in \([0,2\pi)\). For
\(1/3\le R\le R_*\) and \(0\le\alpha\le1\), define
\[
 Q_{R,\alpha}
 =\sqrt R\,\alpha\cos U+\sqrt{1-R}\sqrt{1-\alpha^2}\cos V,
 \qquad
 \mathcal P(R,s,\alpha)
 =\Abar(R)\,\mathbb E\bigl(\chi_s(Q_{R,\alpha})^2-c_R\bigr)_+.
\]
When \(R,\alpha\) are fixed, we abbreviate \(Q_{R,\alpha}\) to \(Q\).

\begin{proposition}[Product majorant]\label{prop:global-product}
For every centre ${\centerpar}=-sp$,
\[
 J(\weights,{\centerpar})\le\sup_{0\le\alpha\le1}\mathcal P(R,s,\alpha),
\]
where
\[
 \mathcal P(R,s,\alpha)=\frac{\Abar(R)}{4\pi^2}
 \int_0^{2\pi}\!\int_0^{2\pi}
 f_{R,s}\bigl(\sqrt R\alpha\cos u+
       \sqrt{1-R}\sqrt{1-\alpha^2}\cos w\bigr)\,du\,dw.
\]
\end{proposition}

\begin{proof}
Separate the first coordinate contribution, with $\alpha=|p_1|$.
For fixed $u$, the other two terms combine into one cosine in $w$,
whose amplitude is at most
\[
 \sqrt{r_2}|p_2|+\sqrt{r_3}|p_3|
       \le\sqrt{1-R}\sqrt{1-\alpha^2}.
\]
For every convex function $f$, the circular average of $f(x+h\cos w)$
is nondecreasing in $h\ge0$ on its admissible range. Indeed, for
$0\le h_1\le h_2$ with $h_2>0$, set $t=h_1/h_2$. Convexity yields
\[
 f(x+h_1\cos w)\le (1-t)f(x)+t f(x+h_2\cos w).
\]
Jensen's inequality gives $f(x)\le\Av_w f(x+h_2\cos w)$, since
$\Av_w\cos w=0$. Integrating proves the assertion, including $h_2=0$
by equality. This argument also applies to the nonsmooth convex
positive-part kernel. The induced area
element is constant. Finally replace its total area by the larger
$\Abar(R)$.
\end{proof}

\begin{proposition}[Cap plus quadratic majorant]\label{prop:global-quadcap}
Choose $t>\sqrt R$, $\xi\in\mathbb R$, and $C\ge0$ such that
\[
 f_{R,s}(q)\le {Q_{\mathrm{cap}}}(q):=C(q-\xi)^2\qquad(-1\le q\le t).
\]
Then $J(\weights,{\centerpar})$ is bounded by
\begin{align*}
 \mathcal U(R,s,m;t,\xi,C)
  ={}&\Abar(R)C(R/2+\xi^2)\\
   &+2\pi(1-R)\int_t^m
       \frac{[f_{R,s}(q)-C(q-\xi)^2]_+}
       {q-Rm-\sqrt{R(1-R)(1-m^2)}}\,dq
\end{align*}
when $m=\max(p\cdot F_{\weights})>t$. If $m\le t$, the first term alone is a
valid upper bound.
\end{proposition}

\begin{proof}
Integrate ${Q_{\mathrm{cap}}}$ using \eqref{eq:global-moments} and the area envelope
\eqref{eq:global-area-envelope}. Its positive excess is supported above $t$,
where Lemma~\ref{lem:global-cap} applies.
\end{proof}

\subsection{Constants, a conformal square, and the quadratic bulk bound}
\label{sec:analytic-highcaps}
Throughout \(1/3\le R\le R_*\), put
\[
 {c_R}=\frac{L}{\Abar(R)}
 =\frac{2(1/\sqrt3-1/\pi)}{\sqrt{(1-R)(1+3R)}}.
\]
We also use the untruncated kernel notation \(f_s(q)=\chi_s(q)^2\).
The rational bounds used below give
\begin{equation}\label{radial:basic-constants}
 \frac49<{c_R}<\frac35,\qquad
 \frac{\Abar(R)}{4\pi}<\frac{11}{6},\qquad
 \frac7{50}<(1-R){c_R}<\frac23.
\end{equation}
Here \(c_R\) increases and \((1-R){c_R}\) decreases.
For an explicit check, write \(w_*=1/\sqrt3-1/\pi\). The bounds
\[
 \frac{333}{106}<\pi<\frac{355}{113}<\frac{22}{7},
 \qquad \frac{265}{153}<\sqrt3<\frac{1351}{780}
\]
give \(25903/100000<w_*<25905/100000\).
Also \((1-R_*)(1+3R_*)>49/64\), so
\[
 {c_R}<\frac{16}{7}\frac{25905}{100000}<\frac35.
\]
For the lower bound on \((1-R){c_R}\), use \(R_*<19/25\) and
\(\sqrt{3/41}>541/2000\), obtaining
\[
 (1-R){c_R}>2\frac{25903}{100000}\frac{541}{2000}
 =\frac{14013523}{100000000}>\frac7{50}.
\]
The other bounds follow from \({c_{1/3}}=1-\sqrt3/\pi>4/9\),
\(\pi/\sqrt3<11/6\), and monotonicity. In particular, for \(R\le1/2\),
\[
 (1-R){c_R}\ge\frac{2w_*}{\sqrt5}
 >\frac89\frac{25903}{100000}>\frac9{40}.
\]

All stated bounds for \(\pi\) can themselves be obtained from the elementary
Machin identity \(\pi=16\arctan(1/5)-4\arctan(1/239)\).
For $0<x<1$ put $S_n(x)=\sum_{j=0}^n(-1)^j x^{2j+1}/(2j+1)$.
The alternating-series remainder gives
$S_{2k+1}(x)<\arctan x<S_{2k}(x)$. Thus the following are purely
rational comparisons:
\[
 \frac{333}{106}<16S_5(1/5)-\frac4{239}<\pi
 <16S_6(1/5)-4\left(\frac1{239}-\frac1{3\cdot239^3}\right)
 <\frac{355}{113}.
\]
The first and last inequalities follow by multiplication by the positive
common denominators of the displayed finite sums. The square-root
bounds are obtained by squaring positive rational endpoints. 

For \(s>0\), write \(R=\cos^2\vartheta_R\), \(m=\cos\vartheta_m>\sqrt R\), and
\(e=e_R(m)=\cos(2\vartheta_R-\vartheta_m)\), where $\vartheta_R,\vartheta_m\in[0,\pi/2]$. With \(P=1-s^2\), \(B=1+s^2\), and
\(X=\cos(\vartheta_R-\vartheta_m)\), direct expansion gives
\begin{equation}\label{radial:conformal-square}
 (B-2sm)(B-2se)
 =4s^2\left(X-\frac{B\sqrt R}{2s}\right)^2+(1-R)P^2.
\end{equation}
This is the conformal-square identity used in the radial cap estimates.

\subsection{A quadratic estimate for smaller height maxima}

\begin{lemma}\label{radial:bulk-bound}
Suppose $1/2\le s<1$, $0<M\le1$, and $m\le M$. Then
\begin{equation}\label{radial:bulk-simple}
 \frac{J_s(\weights,p)}{4\pi}
 \le\frac{\Abar(R)}{4\pi}
 \pos{\chi_s(M)^2-\frac49}\frac{R}{2M^2+R}.
\end{equation}
\end{lemma}

\begin{proof}
Use the lower truncation $4/9$ from \eqref{radial:basic-constants} and
put $b=R/(2M)>0$. If $\chi_s(M)^2\le4/9$, the integrand vanishes.
Otherwise we claim the quadratic bound
\[
 \pos{\chi_s(q)^2-4/9}
 \le\frac{\chi_s(M)^2-4/9}{(M+b)^2}(q+b)^2
 \qquad(-1\le q\le M).
\]
By \eqref{eq:truncation-height}, the positive support starts at
$q_0(s,4/9)=5s/4-1/(4s)>0$ for $s\ge1/2$, so $q+b>0$ there.
For $z=B-2sq$, the derivative of
$(\chi_s(q)^2-4/9)/(q+b)^2$ has the sign of
\[
 \frac49 z^3-2P^2z+P^2(B+2sb).
\]
Its minimum for $z\ge0$ is
$P^2[B+2sb-2\sqrt{2/3}\,P]$, which is positive because
$B/P\ge5/3>2\sqrt{2/3}$ for $s\ge1/2$.
Thus the ratio is increasing on its positive support, proving the
claimed bound. Integrating against $d\mu_{\weights}$ and using
\eqref{eq:global-moments} gives
$\Av_{\mu_{\weights}} q=0$ and $\Av_{\mu_{\weights}} q^2\le R/2$. Finally,
\[
 \frac{R/2+b^2}{(M+b)^2}=\frac{R}{2M^2+R},
\]
and $A(\weights)\le\Abar(R)$ completes the proof.
\end{proof}

\begin{lemma}[A quadratic cap majorant]\label{lem:quadratic-cap-majorant}
Let \(0<s<1\), \(P=1-s^2\), \(B=1+s^2\), and
\[
 f_s(q)=\frac{P^2}{(B-2sq)^2},\qquad
 q_0=q_0(s,4/9)=\frac{5s^2-1}{4s}.
\]
Here $q_0$ is the truncation height in \eqref{eq:truncation-height}.
Suppose \(b>0\), \(q_0+b>0\), and
\[
 B+2sb\ge\sqrt{\frac83}\,P.
\]
Then \(q\mapsto(f_s(q)-4/9)/(q+b)^2\) is nondecreasing on
\((q_0,1]\). If \(q_0<t<1\), the polynomial
\[
 {Q_{\mathrm{quad}}}(q)=\frac{f_s(t)-4/9}{(t+b)^2}(q+b)^2
\]
majorizes \((f_s(q)-4/9)_+\) for \(-1\le q\le t\), and
\(f_s(q)-4/9-{Q_{\mathrm{quad}}}(q)\ge0\) for \(t\le q\le1\).
\end{lemma}
\begin{proof}
On the positive support \(q>q_0\), both \(q+b\) and \(B-2sq\)
are positive. With \(\xi=B-2sq>0\), the sign of the derivative
of the quotient is the sign of
\[
 P^2(B+2sb-2\xi)+\frac49\xi^3.
\]
The elementary cubic maximum
\[
 \sup_{\xi\ge0}
 \left(2\xi-\frac{4\xi^3}{9P^2}\right)
 =\sqrt{\frac83}\,P
\]
proves the monotonicity. On \(q\le q_0\) the positive part
vanishes and \({Q_{\mathrm{quad}}}\ge0\); on the remaining intervals the assertions
follow from the quotient comparison at \(t\).
\end{proof}

\subsection{Fourier moments}\label{sec:fourier-moments}
The following bounds are elementary Fourier identities followed by AM--GM; they hold for $1/3\le R\le2/3$:
\begin{equation}\label{eq:mom}
 \Av_{\mu_{\weights}} q=0,\qquad \Av_{\mu_{\weights}} q^2\le R/2,\qquad
 |\Av_{\mu_{\weights}} q^3|\le R/4,\qquad \Av_{\mu_{\weights}} q^4\le5R^2/8.
\end{equation}
Indeed, after absorbing phases write $q=\sum_{i=1}^3 a_i\cos\theta_i$ with $\theta_3=-\theta_1-\theta_2+\phi$, $a_i\ge0$, $S=\sum a_i^2\le R$. Then
\[
\Av_{\mu_{\weights}} q^2=S/2,\quad \Av_{\mu_{\weights}} q^3=\tfrac32a_1a_2a_3\cos\phi,\quad
\Av_{\mu_{\weights}} q^4=\tfrac34S^2-\tfrac38\sum_i a_i^4,
\]
where $a_i=\sqrt{r_i}|p_i|$ and $S=\sum_i a_i^2$.
The identities follow by taking the constant Fourier coefficient:
in degree four only terms with even multiplicities have nonzero average.
By AM--GM, $a_1a_2a_3\le(S/3)^{3/2}$ and
$\sum_i a_i^4\ge S^2/3$. Consequently
$\tfrac32(S/3)^{3/2}\le R/4$ for $R\le2/3$ and
$\tfrac34 S^2-\tfrac38(S^2/3)=5S^2/8\le5R^2/8$.
This proves all the bounds in \eqref{eq:mom}.

\begin{theorem}[Small centres]\label{thm:small}
For every fundamental embedding with \(1/3\le R\le R_*\) and
\(0\le s\le1/4\), \(J_s(\weights,p)\le4\pi\). Equality holds precisely for
the equal-weight torus and \(s=0\).
\end{theorem}
Section~\ref{sec:small-centres} proves Theorem~\ref{thm:small}.
Section~\ref{lw:quartic} proves the low-weight quartic estimates of
Theorem~\ref{thm:quartic-small-weights}; both sections give the scalar sign arguments explicitly.

\section{Proof of the sharp truncated-area inequality}\label{sec:analytic-closure}

The proof is organized by the conformal radius and the largest weight.
The geometric input is the radial cap comparison in
Lemma~\ref{lem:radial-cap-density}; the reduction to independent circle
variables is Proposition~\ref{prop:global-product}. For small and
moderate radii, Fourier moments and polynomial majorants control the
truncated kernel, supplemented by radial estimates for small weights.
For large radii, the radial comparison is used directly where the
height cap is sufficiently high; elsewhere the product estimate controls
the full conformal integral and the part removed by truncation. These
arguments are given in Appendices~\ref{app:small-moderate} and
\ref{app:large-centres}, respectively, with their full algebraic
comparisons. The exact interval-positivity identities are in
Appendix~\ref{app:explicit-arithmetic}.

\begin{proof}[Proof of Theorem~\ref{thm:global-truncation}]
The parameter ranges are listed in Table~\ref{tab:parameter-cover}.
Theorem~\ref{thm:small} treats $0\le s\le1/4$, and
Theorem~\ref{bs:theorem} treats $1/4\le s\le3/10$, in both cases
for every $1/3\le R\le R_*$.

For $3/10\le s\le2/5$, Theorem~\ref{thm:quartic-small-weights}
covers $1/3\le R\le7/20$,
Theorem~\ref{thm:central-product-strip} covers $7/20\le R\le9/20$,
and Theorem~\ref{thm:whole-low-radius-product} covers $9/20\le R\le R_*$.
For $2/5\le s\le1/2$, the small-weight range $1/3\le R\le7/20$
is covered by Theorem~\ref{thm:quartic-small-weights} up to
$s=17/40$ and by Theorem~\ref{mid:strip} thereafter.
Theorem~\ref{thm:octic-product-strip} covers $7/20\le R\le9/20$,
and Theorem~\ref{thm:whole-low-radius-product} covers $9/20\le R\le R_*$.

For $1/2\le s\le3/5$, use Theorem~\ref{thm:fixed-product-mid}
on $1/3\le R\le4/9$ and Theorem~\ref{thm:wide-product-three-fifths}
on $4/9\le R\le R_*$.
For $3/5\le s\le13/20$, Theorem~\ref{thm:central-small-weights}
covers $R\le1/2$, and Theorem~\ref{thm:sixth-moment-product}
covers $R\ge1/2$.

At larger radii, Theorem~\ref{thm:central-small-weights} still covers
$R\le1/2$.
For $R\ge1/2$ and $13/20\le s\le(1+R)/2$,
Lemma~\ref{lem:product-direction} reduces the product supremum to
$x=\alpha^2\ge R$. Put
\[
 \varkappa=4R(1-R),\qquad
 y=\omega(x-(1-R))-(2R-1).
\]
The three ranges
\[
 y\le\frac3{10}\sqrt{{\varkappa}u},\qquad
 \frac3{10}\sqrt{{\varkappa}u}\le y\le\frac12\sqrt{{\varkappa}u},\qquad
 y\ge\frac12\sqrt{{\varkappa}u}
\]
cover every real $y$.
Propositions~\ref{prop:inner-product}, \ref{prop:middle-product},
and \ref{prop:outer-product} give the strict product bound on
these respective ranges.
Finally, Theorem~\ref{thm:global-highcap-linear} covers
$s\ge(1+R)/2$.

All endpoints in these joins are included.
For fixed $R$ and $s<1$, the function
$\alpha\mapsto\mathcal P(R,s,\alpha)$ is continuous on $[0,1]$:
the conformal kernel has a positive denominator uniformly on the
compact set of amplitudes and angles. Thus the supremum of the
pointwise strict product bounds is also strictly less than $4\pi$.
Proposition~\ref{prop:global-product} transfers each product bound to
$J$, while the cap estimates and the moment estimates for the three
Fourier modes of $F_{\weights}$ already bound $J$ directly. All bounds for $s\ge1/4$ are strict.
The equality statement in Theorem~\ref{thm:small} therefore proves
that equality holds precisely for $s=0$ and
$\weights=(1/3,1/3,1/3)$.
\end{proof}

\section{The complementary conformal moduli}
\label{sec:exterior}

Let $0\le {a_{\mathrm{mod}}}\le1/2$ and ${b_{\mathrm{mod}}}\ge44/25$. On the flat torus of modulus
${a_{\mathrm{mod}}}+i {b_{\mathrm{mod}}}$, use the embedding into $S^3$
\[
 \Psi_\rho(x,y)=\left(\sqrt\rho\,e^{2\pi i y/{b_{\mathrm{mod}}}},
       \sqrt{1-\rho}\,e^{2\pi i(x-{a_{\mathrm{mod}}} y/{b_{\mathrm{mod}}})}\right),
 \qquad \rho=1-\frac1{{b_{\mathrm{mod}}}^2+8/3}>\frac23.
\]
Its Dirichlet-energy and induced-area densities are constant. Its total
Dirichlet energy and total induced area are
\[
 E_0=\frac{2\pi^2}{{b_{\mathrm{mod}}}}\{\rho+({a_{\mathrm{mod}}}^2+{b_{\mathrm{mod}}}^2)(1-\rho)\},\qquad
 A_0=4\pi^2\sqrt{\rho(1-\rho)}.
\]

The energy and area densities remain proportional after target conformal
transformations. The image is a product torus with weights
\((\rho,1-\rho)\), and its conformal area equals
\[
 \frac{8\pi^2}{3\sqrt3\sqrt\rho}
 =A_0\,\frac{2}{3\sqrt3\,\rho\sqrt{1-\rho}}.
\]
This formula is proved in Lemma~\ref{lem:exterior-product-area},
following Kang~\cite[Lemmas~2.3--2.4]{Kang} with the elimination and
boundary calculations made explicit.
Therefore
\[
 \sup_M\E(M\circ\Psi_\rho)
 =E_0\frac{2}{3\sqrt3\,\rho\sqrt{1-\rho}}.
\]
Lemma~\ref{lem:global-fold} and the coordinatewise Rayleigh inequalities give
$\lambda_2\Area\le4\sup_M\E(M\circ\Psi_\rho)$, since each of the two
piece energies is at most the full single-map supremum. Hence, with
$B={b_{\mathrm{mod}}}^2$,
\[
 \lambda_2\Area
 \le\frac{16\pi^2}{3\sqrt3}\,H(B),\qquad
 H(B)=\frac{(2B+23/12)\sqrt{(B+8/3)/B}}{B+5/3}.
\]
The logarithmic derivative is
\[
 \frac{H'(B)}{H(B)}
 =\frac{17/12}{(2B+23/12)(B+5/3)}
       -\frac{8/3}{2B(B+8/3)}.
\]
After multiplication by the positive common denominator, its numerator
is $-(5/2)B^2-(58/9)B-230/27$. Thus $H$ decreases for $B>0$.
It remains to bound the value at $B_0=1936/625$. At this point
\[
 \sqrt{\frac{B_0+8/3}{B_0}}=\sqrt{\frac{1351}{726}}
                  <\frac{6821}{5000}.
\]
The rational inequalities $333/106<\pi<22/7$ and
$265/153<\sqrt3<97/56$ yield
\begin{align*}
 \frac{16\pi^2}{3\sqrt3}H(B_0)
 &<\frac{16(22/7)^2}{3(265/153)}
      \frac{2B_0+23/12}{B_0+5/3}\frac{6821}{5000}
 <\frac{707}{10},\\
 \Cstar&>
    \frac{8(333/106)^2}{97/56}+8\frac{333}{106}
 >\frac{707}{10}.
\end{align*}
The displayed endpoint comparisons are exact rational inequalities
after applying the stated bounds for $\pi$ and $\sqrt3$. The exterior region thus has
a uniform strict upper bound below $\Cstar$.

\subsection{Conformal area of a product torus}
\label{sec:exterior-conformal-area-proof}

We record the product-torus computation of Bryant and Kang
\cite[\S3 and Appendix~A]{Bryant}\cite[Lemmas~2.3--2.4]{Kang}, with
the algebraic elimination and the boundary limit made explicit.

\begin{lemma}\label{lem:exterior-product-area}
For $2/3<\rho<1$, the product torus
\[
 X_\rho(u_1,u_2)=
 (\sqrt\rho\,e^{iu_1},\sqrt{1-\rho}\,e^{iu_2})\subset S^3
\]
has conformal area
\[
 \sup_{M\in\operatorname{Conf}(S^3)}\Area(M\circ X_\rho)
 =\frac{8\pi^2}{3\sqrt3\sqrt\rho}.
\]
\end{lemma}

\begin{proof}
Put $h=1-\rho$ and $A_0=4\pi^2\sqrt{\rho h}$.
The squared conformal factor can equivalently be parametrized as
$(1-|\boldsymbol{\xi}|^2)/(1-\boldsymbol{\xi}\cdot X_\rho)^2$, with $|\boldsymbol{\xi}|<1$;
the relation to the ball parameter ${\centerpar}$ in $\chi_{\centerpar}$ is
$\boldsymbol{\xi}=-2{\centerpar}/(1+|{\centerpar}|^2)$.
Using independent rotations in the two coordinate planes, the area
ratio takes the form
\[
 \frac{\Area(M\circ X_\rho)}{A_0}
 =N(\lambda,\mu)\frac1{4\pi^2}
 \int_0^{2\pi}\!\int_0^{2\pi}
 \frac{du_1\,du_2}{(1-\lambda\cos u_1-\mu\cos u_2)^2},
\]
where
\[
 N=1-\frac{\lambda^2}{\rho}-\frac{\mu^2}{h},\qquad
 (\lambda,\mu)\in\Omega_\rho
 :=\{\lambda,\mu\ge0:\lambda^2/\rho+\mu^2/h<1\}.
\]
Cauchy--Schwarz gives $\lambda+\mu<1$ in this domain. Define
\[
 \Delta_-=1-(\lambda-\mu)^2,\qquad
 \Delta_+=1-(\lambda+\mu)^2,\qquad
 z=\frac{4\lambda\mu}{\Delta_-}\in[0,1).
\]
For completeness, the circle integration formula of
\cite[(3.6), Appendix~A]{Bryant} follows from the same conditioning
argument as Lemma~\ref{lem:full-average-series}. After integrating the
common angle, using the elementary circular integral, and halving the
difference angle, the double average becomes
\[
 \frac{2}{\pi\Delta_-^{3/2}}\int_0^{\pi/2}
       (1-z\sin^2\theta)^{-3/2}\,d\theta.
\]
For $y=\sin^2\theta$ one has the identity
\[
 (1-z)(1-zy)^{-3/2}-(1-zy)^{1/2}
 =-z\frac{d}{d\theta}\frac{\sin\theta\cos\theta}{\sqrt{1-zy}}.
\]
The boundary term vanishes at both ends. Integrating gives
$\int_0^{\pi/2}(1-z\sin^2\theta)^{-3/2}\,d\theta=\mathcal E_{\rm ell}(z)/(1-z)$.
As $1-z=\Delta_+/\Delta_-$, we obtain
\[
 \frac1{4\pi^2}\int_0^{2\pi}\!\int_0^{2\pi}
 \frac{du_1\,du_2}{(1-\lambda\cos u_1-\mu\cos u_2)^2}
 =\frac{2\mathcal E_{\rm ell}(z)}{\pi\sqrt{\Delta_-}\Delta_+},
 \quad
 \mathcal E_{\rm ell}(z)=\int_0^{\pi/2}
       \sqrt{1-z\sin^2\theta}\,d\theta.
\]
Thus $\mathcal E_{\rm ell}(z)=\mathrm E(\sqrt z)$ in the notation of
Section~\ref{sec:global-spectral}. Concavity of the square root gives
$\mathcal E_{\rm ell}(z)\le(\pi/2)(1-z/4)$. Consequently the area
ratio is at most
\[
 I(\lambda,\mu)
 :=\frac{N(\lambda,\mu)(1-\lambda^2+\lambda\mu-\mu^2)}
         {\Delta_-^{3/2}\Delta_+}.
\]

We first exclude critical points with $\lambda,\mu>0$.
Set $B=\rho/h>2$, $t=\mu/\lambda>0$, $x=\lambda^2$, and
$D_0=1-\lambda^2+\lambda\mu-\mu^2$.
Consider the following quadratic polynomials in $x$:
\begin{align*}
 P_1(t,x)={}&\frac1{2t}\bigl[
 (1-Bt)(1+t)(1-(1-t)^2x)^2\\
 &\hspace{32mm}-(1+Bt)(1-t)
 (1-(1+t^2)x)(1-(1+t)^2x)\bigr],\\
 P_2(t,x)={}&2(1+B)(B^{-1}+t)
 (1-(1+t)^2x)(1-(1-t+t^2)x)\\
 &-3(1+t)(1-(1-t)^2x)
 (1-(1+B)(B^{-1}+t^2)x).
\end{align*}
Direct differentiation gives the identities
\[
\begin{aligned}
 P_1(t,x)&=\frac{\Delta_-\Delta_+D_0}{3\mu}
 \left(\partial_\mu-\frac{B\mu}{\lambda}\partial_\lambda\right)\log I,\\
 P_2(t,x)&=-\frac{N\Delta_+D_0}{\lambda}
 (\partial_\lambda+\partial_\mu)\log I.
\end{aligned}
\]
They are valid throughout the interior with $\lambda,\mu>0$, where
$N,\Delta_-,\Delta_+>0$ and $D_0=\Delta_++3\lambda\mu>0$.
Thus both polynomials vanish at every such critical point.
Here is the algebraic elimination explicitly. Write
$P_1=\alpha_2x^2+\alpha_1x+\alpha_0$ and
$P_2=\beta_2x^2+\beta_1x+\beta_0$. To avoid large expanded expressions,
set
\begin{align*}
 U_1&=(1-Bt)(1+t),&U_2&=(1+Bt)(1-t),\\
 \ell&=(1-t)^2,&m_1&=1+t^2,&m_2&=(1+t)^2,\\
 m_3&=1-t+t^2,&m_4&=(1+B)(B^{-1}+t^2),\\
 V_1&=2(1+B)(B^{-1}+t),&V_2&=3(1+t).
\end{align*}
Expanding the two products in $P_1,P_2$ gives
\begin{align*}
 (\alpha_0,\alpha_1,\alpha_2)
 &=\frac1{2t}\bigl(U_1-U_2,-2U_1\ell+U_2(m_1+m_2),
              U_1\ell^2-U_2m_1m_2\bigr),\\
 (\beta_0,\beta_1,\beta_2)
 &=\bigl(V_1-V_2,-V_1(m_2+m_3)+V_2(\ell+m_4),
             V_1m_2m_3-V_2\ell m_4\bigr).
\end{align*}
The first numerator is divisible by $2t$ coefficientwise. Put
$D_{ij}=\alpha_i\beta_j-\alpha_j\beta_i$. The determinant of the
Sylvester matrix is $D_{20}^2-D_{21}D_{10}$. Alternatively, a common
root $x$ implies $D_{21}x+D_{20}=D_{20}x+D_{10}=0$, which also
implies $D_{20}^2-D_{21}D_{10}=0$, including degree drops. Substitution
of the preceding coefficients gives
\begin{align*}
 D_{20}&=\frac{(t+1)(Bt+1)}{B}
 (3B^2t^3+B^2t-4Bt^2+4Bt-t^2-3),\\
 D_{21}&=-\frac{(t-1)(t+1)(Bt+1)}{B}
 \bigl(3B^2t^4-7B^2t^3+B^2t^2-B^2t\\
 &\hspace{35mm}-8Bt^3+8Bt^2-8Bt-t^3+t^2-7t+3\bigr),\\
 D_{10}&=-\frac{Bt+1}{B}(3B^2t^2+B^2t-t-3).
\end{align*}
Squaring the first expression, subtracting the product of the other
two, and collecting powers of $B$ yields the exact factorization
\[
 \operatorname{Res}_x(P_1,P_2)
 =\frac{4t(t+1)(Bt-1)^2(Bt+1)^2}{B^2}\,Q_B(t),
\]
where
\[
 Q_B(t)=B^2t(9t^2-2t+1)+8Bt^2+(8B-2)t+t^2+9>0.
\]
Indeed $9t^2-2t+1=9(t-1/9)^2+8/9>0$ and all the other
displayed summands are positive or nonnegative. Thus a common root
requires $t=1/B$. At that value,
\[
 P_2(1/B,x)=\frac{B+1}{B^5}\bigl((B+1)^2x-B^2\bigr)^2.
\]
Its only root gives $(\lambda,\mu)=(\rho,h)$ and
$\lambda+\mu=1$, outside $\Omega_\rho$. Hence $I$ has no such
interior critical point.

The only singular point on the elliptical part of the boundary is
$(\rho,h)$. The identity
\[
 \Delta_+-N=
 \left(\sqrt{h/\rho}\,\lambda-\sqrt{\rho/h}\,\mu\right)^2\ge0
\]
implies $0<N/\Delta_+\le1$ in $\Omega_\rho$, and therefore
\[
 \limsup_{(\lambda,\mu)\to(\rho,h)} I(\lambda,\mu)
 \le\frac3{8\sqrt{\rho h}}
 <\frac2{3\sqrt3\,\rho\sqrt h}.
\]
The strict comparison is equivalent to $243\rho<256$.
At the other points of the elliptical boundary the limit is zero.
On the axis $\mu=0$, setting $z_0=\lambda^2$ gives
\[
 I(\lambda,0)=\frac{1-z_0/\rho}{(1-z_0)^{3/2}},
 \qquad 0\le z_0<\rho.
\]
Its maximum occurs at $z_0=3\rho-2$ and equals
$2/(3\sqrt3\,\rho\sqrt h)$. On the other axis the same expression
with $\rho$ replaced by $h<1/3$ is decreasing and is at most $1$.
Compactness away from the singular boundary point, together with the
absence of interior critical points, now proves
\[
 \sup_{\Omega_\rho} I=\frac2{3\sqrt3\,\rho\sqrt h}.
\]
At $\mu=0$, $\lambda^2=3\rho-2$, the elliptic-integral bound is an
equality, so this is also the supremum of the actual area ratio.
Multiplication by $A_0$ proves the formula.
\end{proof}

\section{Proof of the main theorem}
\label{sec:sharpness}

\begin{proof}[Proof of Theorem~\ref{thm:global-spectral}]
\emph{The upper bound.}
By uniformization and pullback by a diffeomorphism, the conformal class
of any smooth metric $g$ on $T^2$ has a modulus ${a_{\mathrm{mod}}}+i {b_{\mathrm{mod}}}$ with
$0\le {a_{\mathrm{mod}}}\le1/2$ and ${b_{\mathrm{mod}}}\ge\sqrt{1-{a_{\mathrm{mod}}}^2}$.
Theorem~\ref{thm:global-truncation}, proved in
Section~\ref{sec:analytic-closure}, and the folding argument in
Section~\ref{sec:folding-implication} give
\[
 \lambda_2(T^2,g)\Area(T^2,g)<\Cstar
\]
when ${b_{\mathrm{mod}}}\le44/25$, that is, when the modulus lies in $\mathcal K$.
Section~\ref{sec:exterior} gives the same strict inequality for
${b_{\mathrm{mod}}}\ge44/25$. These two ranges cover every modulus, including their
common endpoint. Taking the supremum over all smooth metrics therefore
yields $\Lambda_2(T^2)\le\Cstar$.

\emph{Sharpness and nonattainment.}
Write
\[
 \Lambda_1(T^2,\mathcal C)
 =\sup_{g\in\mathcal C}\lambda_1(T^2,g)\Area(T^2,g).
\]
For the reverse inequality, the conformal bubbling construction gives
\begin{equation}\label{eq:conformal-gluing}
 \Lambda_2(T^2,\mathcal C)
       \ge\Lambda_1(T^2,\mathcal C)+8\pi;
\end{equation}
see \cite[Theorem~B]{ColboisElSoufi} and \cite{KNPP}. The equilateral embedding has area $\Aeq$ and its induced metric is
$g_{(1/3,1/3,1/3)}=(8\pi^2/3)g_{\mathrm{eq}}$. The spectrum of this
induced metric is
\[
 \lambda_{m,n}=2(m^2-mn+n^2),\qquad (m,n)\in\mathbb Z^2.
\]
The quadratic form is a positive integer for every nonzero lattice
pair and equals \(1\) at \((1,0)\), so the first positive eigenvalue
is \(2\).
Consequently \eqref{eq:conformal-gluing} gives
$\Lambda_2(T^2,[g_{\mathrm{eq}}])\ge2\Aeq+8\pi=\Cstar$.
More precisely, apply the conformal gluing construction in the proof
of \cite[Theorem~B]{ColboisElSoufi}, together with the spectral
convergence in \cite[Lemma~2.1]{ColboisElSoufi}, to the equilateral
flat torus with the above induced metric and the unit round sphere.
By taking a diagonal sequence of the gluing and smoothing parameters,
one obtains smooth metrics \(g_j\in[g_{\mathrm{eq}}]\) whose areas
converge to the sum of the areas of the two components and whose
eigenvalues converge to the ordered spectrum of their disjoint union,
with multiplicities.

Both limiting components have first positive eigenvalue \(2\).
Since their disjoint union has two zero eigenvalues, it follows that
\[
\Area(T^2,g_j)\longrightarrow\Aeq+4\pi,
\qquad
\lambda_1(T^2,g_j)\longrightarrow0,
\qquad
\lambda_2(T^2,g_j)\longrightarrow2.
\]
Consequently,
\[
\lambda_2(T^2,g_j)\Area(T^2,g_j)
\longrightarrow 2\Aeq+8\pi=\Cstar.
\]
This proves sharpness in the equilateral conformal class, with one
round spherical bubble. The strict upper bound established above
excludes attainment by any smooth metric.

\emph{Uniform separation of the other conformal classes.}
Let $\mathcal V$ be an open neighbourhood of the equilateral class in
moduli space. It suffices to obtain a uniform gap on the compact subset
$\mathcal K_{\mathcal V}$ of $\mathcal K$ consisting of moduli whose
classes are outside $\mathcal V$, since Section~\ref{sec:exterior}
already supplies a uniform gap for $b_{\mathrm{mod}}\ge44/25$. If
$\mathcal K_{\mathcal V}$ is empty there is nothing to prove in $\mathcal K$.
Use the fixed angular torus of Section~\ref{sec:global-geometry} to
parametrize all the embeddings $F_{\weights}$; the weights and the
area $A(\weights)$ vary continuously with the modulus, and the area
measure is $A(\weights)$ times the fixed Haar probability measure.
In particular $F_{\weights}$ depends uniformly continuously on the
modulus and the angular variables on this compact family.

Write
\[
 I_{\mathrm{conf}}(\weights,\centerpar)=\int_{T^2}\chi_{\centerpar}(F_{\weights})^2\,dA_{\weights}.
\]
For any sequence of moduli in $\mathcal K_{\mathcal V}$ converging to a
limit modulus with weights $\weights_\infty$, and any centres
$\centerpar_j\to\centerpar_\infty\in\overline B^6$, one has
\begin{equation}\label{eq:single-boundary}
 \limsup_{j\to\infty} I_{\mathrm{conf}}(\weights_j,\centerpar_j)\le4\pi
 \qquad\text{if } |\centerpar_\infty|=1.
\end{equation}
Indeed, decompose the integral as
\[
 I_{\mathrm{conf}}(\weights_j,\centerpar_j)=J(\weights_j,\centerpar_j)+
 \int_{T^2}\min\{\chi_{\centerpar_j}(F_{\weights_j})^2,
 L/A(\weights_j)\}\,dA_{\weights_j}.
\]
The second term tends to zero by dominated convergence on the fixed
angular torus: its integrand with respect to Haar probability is bounded
by $L$, and tends pointwise to zero away from the inverse image of
$-\centerpar_\infty$ under $F_{\weights_\infty}$. This inverse image has
measure zero by embeddedness, including the case $r_3=0$.
Theorem~\ref{thm:global-truncation} bounds the first term by $4\pi$.
When $|\centerpar_\infty|<1$, the integral $I_{\mathrm{conf}}$ is continuous in the weights
and the centre, since its denominator is bounded away from zero locally.
Moreover for every interior centre with weights arising in
$\mathcal K_{\mathcal V}$, the equality statement of
Theorem~\ref{thm:global-truncation} gives
\[
 I_{\mathrm{conf}}(\weights,\centerpar)\le J(\weights,\centerpar)+L<\Aeq.
\]
Now consider the two-centre quantity
\[
 V_2(\weights)=\sup_{\centerpar,\othercenter\in B^6}
 \int_{T^2}\max\{\chi_{\centerpar}(F_{\weights})^2,
                    \chi_{\othercenter}(F_{\weights})^2\}\,dA_{\weights}.
\]
The supremum of $V_2$ over $\mathcal K_{\mathcal V}$ is strictly below
$\Aeq+4\pi$. To see this, suppose a sequence of moduli and pairs of
centres has integrals tending to $\Aeq+4\pi$ (the upper bound follows
from \eqref{eq:TR-main}). Pass to a subsequence so that the moduli converge
in $\mathcal K_{\mathcal V}$ and the centres converge in
$\overline B^6\times\overline B^6$. If both limit centres are interior,
continuity and the strict inequality in the folding argument, applied to
the integral of the maximum using \eqref{eq:TR-main}, give a contradiction.
If one limit centre lies on the boundary and one is interior, the bound
$\max\{x,y\}\le x+y$, \eqref{eq:single-boundary}, and the strict bound
$I_{\mathrm{conf}}<\Aeq$ for the interior limit give a limsup strictly smaller than
$\Aeq+4\pi$. If both lie on the boundary the limsup is at most
$8\pi<\Aeq+4\pi$. These cases exhaust the possibilities.

The folding argument gives $\Lambda_2(T^2,\mathcal C)\le2V_2(\weights)$
for every class represented in $\mathcal K$. Combining the uniform gap
just proved with the exterior gap proves the assertion for
$\varepsilon_{\mathcal V}$. Since the moduli space is Hausdorff, every
class distinct from the equilateral class is outside some such
neighbourhood. This yields uniqueness; applying the uniform gap to every
neighbourhood also proves convergence in moduli space for every maximizing
sequence.
\end{proof}

\appendix
\section{Small and moderate conformal centres}\label{app:small-moderate}

This appendix proves the bounds used for $0\le s\le3/5$ in
Table~\ref{tab:parameter-cover}. It also establishes the fixed-amplitude
estimate for an adjoining larger-radius range. Each subsection records
its own parameter hypotheses; the common geometric notation is that of
Section~\ref{sec:analytic-roadmap}.

\subsection{Small conformal centres}\label{sec:small-centres}

The following preparatory estimates for Theorem~\ref{thm:small} use
Hermite interpolation, the Fourier moment bounds \eqref{eq:mom}, the Hermite remainder formula, and the area envelope
\eqref{eq:global-area-envelope}, namely
\[
\Abar(R)=2\pi^2\sqrt{(1-R)(1+3R)},\qquad
\Aeq=4\pi^2/\sqrt3,\qquad L=\Aeq-4\pi.
\]
The normalized height variable satisfies
\[
\mathbb E q=0,\quad \mathbb E q^2\le R/2,\quad
|\mathbb E q^3|\le R/4,\quad \mathbb E q^4\le5R^2/8
\qquad(1/3\le R\le2/3).
\]
We use ${\sigma}=s/(1+s^2)$, so $s\le1/4$ is equivalent to ${\sigma}\le4/17$.

\subsubsection*{Notation and identities}

For \(0\le {\sigma}<1/2\), define the kernel and Hermite polynomial below. In this subsection put \(T=4/17\) and restrict to \(0\le {\sigma}\le T\). Write
\[
f_{\sigma}(h)=\frac{1-4{\sigma}^2}{(1-2{\sigma}h)^2},\quad
{D_{\mathrm{height}}}({\sigma})=(2-{\sigma})^3(1-2{\sigma})(5+2{\sigma})^3>0.
\]
Let \(H_{\sigma}\) be the degree-four Hermite interpolant, with double nodes
\(x=-1/5\), \(y=1/4\), and simple node \(1\). Since the fifth derivative of the kernel is nonnegative (and positive for \({\sigma}>0\)), the Hermite remainder proves \(H_{\sigma}\ge f_{\sigma}\) on \([-1,1]\). Its Newton representation is
\[
H_{\sigma}(h)=f_{\sigma}(x)+f_{\sigma}'(x)(h-x)+(h-x)^2\{u_2+u_3(h-y)+u_4(h-y)^2\},
\]
where direct interpolation gives
\[
\begin{aligned}
u_2&=\frac{1200{\sigma}^2(5-{\sigma})(1-2{\sigma})(1+2{\sigma})}{(2-{\sigma})^2(5+2{\sigma})^3},\\
u_3&=\frac{1600{\sigma}^3(20-{\sigma})(1-2{\sigma})(1+2{\sigma})}{(2-{\sigma})^3(5+2{\sigma})^3},\\
u_4&=\frac{1600{\sigma}^4(1+2{\sigma})(2{\sigma}^2-83{\sigma}+50)}{{D_{\mathrm{height}}}({\sigma})}.
\end{aligned}
\]
For clarity, the ordinary monomial coefficients \(H_{\sigma}(h)=\sum h_j({\sigma})h^j\) satisfy
\[
\begin{aligned}
{D_{\mathrm{height}}} h_0&=-(1+2{\sigma})(336{\sigma}^6-488{\sigma}^5+2642{\sigma}^4-1199{\sigma}^3-4630{\sigma}^2+4300{\sigma}-1000),\\
{D_{\mathrm{height}}} h_1&=4{\sigma}(1+2{\sigma})(80{\sigma}^5-2522{\sigma}^4+1199{\sigma}^3+4630{\sigma}^2-4300{\sigma}+1000),\\
{D_{\mathrm{height}}} h_2&=12{\sigma}^2(1+2{\sigma})(294{\sigma}^4-441{\sigma}^3+4630{\sigma}^2-4300{\sigma}+1000),\\
{D_{\mathrm{height}}} h_3&=160{\sigma}^3(1+2{\sigma})(-42{\sigma}^3+923{\sigma}^2-860{\sigma}+200),\\
{D_{\mathrm{height}}} h_4&=1600{\sigma}^4(1+2{\sigma})(2{\sigma}^2-83{\sigma}+50).
\end{aligned}
\]
Each identity follows either by expansion of the Newton expression or by checking the five interpolation conditions.

\subsubsection*{The coefficient signs need only two monotonicities}

Let
\[
p_2({\sigma})=294{\sigma}^4-441{\sigma}^3+4630{\sigma}^2-4300{\sigma}+1000,
\quad p_3({\sigma})=-42{\sigma}^3+923{\sigma}^2-860{\sigma}+200.
\]
Since \(T<1/4\),
\[
p_2'({\sigma})\le\frac{1176}{64}+\frac{9260}{4}-4300<0,
\qquad p_3'({\sigma})\le\frac{1846}{4}-860<0.
\]
At the right endpoint,
\[
p_2(T)=\frac{20021976}{83521}>0,
\qquad p_3(T)=\frac{236808}{4913}>0.
\]
Also \(2{\sigma}^2-83{\sigma}+50\ge50-83/4=117/4>0\). Thus
\(h_2,h_3,h_4>0\) for \(0<{\sigma}\le T\); at zero they vanish.

The coefficient-sign argument is extended to the endpoint $3/8$ in
Section~\ref{lw:quartic}, for the low-weight range.

\subsubsection*{\texorpdfstring{The lower bound \(H_{\sigma}\ge49/100\)}{The lower bound H_sigma >=49/100}}

\begin{lemma}\label{sa:floor}
For \(0\le {\sigma}\le T\) and \(-1\le h\le1\),
\[
H_{\sigma}(h)\ge\frac{49}{100}.
\]

\end{lemma}
\begin{proof}
If \({\sigma}\le1/6\), then
\[
H_{\sigma}(h)\ge f_{\sigma}(h)\ge f_{\sigma}(-1)=\frac{1-2{\sigma}}{1+2{\sigma}}\ge\frac12.
\]
For every \({\sigma}\le T\), if \(h\ge h_*=-31/56\), monotonicity in height gives
\[
f_{\sigma}(h)\ge f_{\sigma}(h_*)\ge f_T(h_*)=49/100.
\]
Here \(\partial_{\sigma} f_{\sigma}(h_*)=4(h_*-2{\sigma})/(1-2{\sigma}h_*)^3<0\).

It remains to consider \(1/6\le {\sigma}\le T\) and \(h\le h_*\). All \(u_j\) are positive. Moreover
\[
\frac{u_3}{u_4}=
\frac{(20-{\sigma})(1-2{\sigma})^2}{{\sigma}(2{\sigma}^2-83{\sigma}+50)}\le\frac32.
\]
Indeed this inequality is equivalent to
\[
G({\sigma})=7{\sigma}^3-\frac{417}{2}{\sigma}^2+156{\sigma}-20\ge0.
\]
On the interval, \(G'({\sigma})=21{\sigma}^2-417{\sigma}+156\ge156-417T>0\), and
\(G(1/6)=13/54>0\).

For \(z=h-y\le-45/56\), the quadratic
\(q_{\sigma}(z)=u_2+u_3z+u_4z^2\) is decreasing as \(z\) increases through this range, because
\[
q_{\sigma}'(z)=u_3+2u_4z\le u_4\left(\frac32-\frac{90}{56}\right)<0.
\]
Consequently
\[
q_{\sigma}(z)\ge q_{\sigma}(-45/56)=A({\sigma}),
\qquad
A({\sigma})=\frac{75{\sigma}^2(1+2{\sigma})P({\sigma})}{49{D_{\mathrm{height}}}({\sigma})},
\]
where
\[
P({\sigma})=7846{\sigma}^4-151673{\sigma}^3+155886{\sigma}^2-53648{\sigma}+7840.
\]
This is positive: since \({\sigma}<1/4\),
\[
P({\sigma})\ge\frac{471871}{4}{\sigma}^2-53648{\sigma}+7840>0.
\]
The last quadratic has positive leading coefficient and discriminant
\(53648^2-471871\cdot7840=-821360736<0\).

Completing the square now gives the \emph{uniform} lower bound
\[
H_{\sigma}(h)\ge f_{\sigma}(x)-\frac{f_{\sigma}'(x)^2}{4A({\sigma})}
=\frac{49}{100}+\frac{N({\sigma})}{300(5+2{\sigma})^3P({\sigma})},
\]
where the exact numerator is
\[
\begin{aligned}
N({\sigma})={}&-675986896{\sigma}^7+9404645728{\sigma}^6+10968909624{\sigma}^5\\
&-16385495072{\sigma}^4+2779263095{\sigma}^3+1452889350{\sigma}^2\\
&-542290000{\sigma}+51940000.
\end{aligned}
\]
To determine its sign, write \(z={\sigma}-1/6\), so
\(0\le z\le7/102<7/100\). Direct differentiation gives
\[
\begin{aligned}
N'(1/6+z)={}&-4731908272z^6+51695966096z^5
+\frac{299688444940}{3}z^4\\
&-\frac{371052314536}{27}z^3
-\frac{343839835028}{27}z^2
+\frac{117684807104}{81}z
-\frac{58575012928}{729}.
\end{aligned}
\]
The positive fourth- and fifth-degree terms can be absorbed into the negative quadratic term. More explicitly,
\[
\begin{aligned}
-\frac{343839835028}{27}
&+\frac{299688444940}{3}\left(\frac7{100}\right)^2
+51695966096\left(\frac7{100}\right)^3\\
&=-\frac{5158512785377046}{421875}<-12000000000.
\end{aligned}
\]
Dropping the remaining negative terms and slightly enlarging the constant and linear coefficients yields
\[
N'(1/6+z)
\le 10^6(-80+1453z-12000z^2)<0.
\]
The final quadratic is negative on the whole real line, since
\(1453^2=2111209<3840000=4\cdot80\cdot12000\).
Therefore
\[
N({\sigma})\ge N(T)=\frac{98568848626560}{410338673}>0.
\]
This proves the lower-floor lemma.
\end{proof}

\subsubsection*{A moment-envelope comparison for every real weight}

Let
\[
B_{\mathrm{height}}(R,{\sigma})=h_0+\left(\frac{h_2}{2}+\frac{h_3}{4}\right)R+
\frac58h_4R^2.
\]
For every real $R$, and \(0 \le  {\sigma} \le  T\),
\[
B_{\mathrm{height}}(R,{\sigma})\le1+\frac98(R-1/3)^2,
\]
with strict inequality for \({\sigma}>0\).

To prove this put \(y=R-1/3\). The numerator of the difference is
\[
{D_{\mathrm{height}}}({\sigma})\left[1+\frac98y^2-B_{\mathrm{height}}(1/3+y,{\sigma})\right]
=a_2({\sigma})y^2+a_1({\sigma})y+a_0({\sigma}),
\]
where
\[
\begin{aligned}
a_0={}&\frac{2{\sigma}^2}{9}
(844{\sigma}^5-23756{\sigma}^4-19549{\sigma}^3+38870{\sigma}^2-32700{\sigma}+9000),\\
a_1={}&-\frac{2{\sigma}^2(1+2{\sigma})}{3}
(2126{\sigma}^4-31589{\sigma}^3+40070{\sigma}^2-26700{\sigma}+9000),\\
a_2={}&\frac18(9000-20700{\sigma}+270{\sigma}^2+11331{\sigma}^3\\
&\hspace{25mm}-401116{\sigma}^4-138160{\sigma}^5+1312144{\sigma}^6-31856{\sigma}^7).
\end{aligned}
\]
A direct lower bound is
\[
a_2({\sigma})\ge
\frac{9000-20700T-401116T^4-138160T^5-31856T^7}{8}
=\frac{143569259111}{410338673}>0.
\]
An exact algebraic simplification gives the discriminant identity
\[
4a_0a_2-a_1^2={\sigma}^2{Q_{\mathrm{disc}}}({\sigma}),
\]
with
\[
\begin{aligned}
{Q_{\mathrm{disc}}}({\sigma})={}&9000000-53400000{\sigma}+78350000{\sigma}^2-29000000{\sigma}^3\\
&-346912500{\sigma}^4+931568460{\sigma}^5+839388049{\sigma}^6\\
&-4485501248{\sigma}^7+5404516544{\sigma}^8+24528864{\sigma}^9\\
&-5247338288{\sigma}^{10}+437885056{\sigma}^{11}-11022720{\sigma}^{12}.
\end{aligned}
\]
It remains to check \({Q_{\mathrm{disc}}}>0\). Its derivative is
\[
\begin{aligned}
{Q_{\mathrm{disc}}}'({\sigma})={}&-53400000+156700000{\sigma}-87000000{\sigma}^2
-1387650000{\sigma}^3\\
&+4657842300{\sigma}^4+5036328294{\sigma}^5-31398508736{\sigma}^6\\
&+43236132352{\sigma}^7+220759776{\sigma}^8-52473382880{\sigma}^9\\
&+4816735616{\sigma}^{10}-132272640{\sigma}^{11}.
\end{aligned}
\]
The quadratic \(156700000{\sigma}-87000000{\sigma}^2\) increases on \([0,T]\). Replace its value and every other positive monomial by their values at \(T\), and drop the remaining negative monomials. This yields the exact upper bound
\[
\begin{aligned}
{Q_{\mathrm{disc}}}'({\sigma})\le{}&-53400000+156700000T-87000000T^2\\
&+4657842300T^4+5036328294T^5+43236132352T^7\\
&+220759776T^8+4816735616T^{10}\\
={}&-\frac{3439632755465815104}{2015993900449}<0.
\end{aligned}
\]
Thus
\[
{Q_{\mathrm{disc}}}({\sigma})\ge {Q_{\mathrm{disc}}}(T)=\frac{8889551284064040000}{582622237229761}>0.
\]
For \({\sigma}>0\), the quadratic in \(y\) has positive leading coefficient and strictly negative discriminant, so it is everywhere positive. At \({\sigma}=0\) the desired assertion is immediate (\(H=1\), \(B_{\mathrm{height}}=1\)).

\subsubsection*{Completion of the small-centre argument}

\begin{proof}[Proof of Theorem~\ref{thm:small}]
For $1/3\le R\le2/3$, we use the preceding coefficient signs,
Hermite lower bound and moment-envelope comparison.
For a fundamental torus let \(A=A(\weights)\le\Abar(R)\),
\({c_R}=L/\Abar(R)\), and \(q=p\cdot F_{\weights}\).
The Fourier moment bounds \eqref{eq:mom}, together with the coefficient signs established in this proof, give
\(\mathbb E H_{\sigma}(q)\le B_{\mathrm{height}}(R,{\sigma})\) for \(1/3\le R\le2/3\).
Lemma~\ref{sa:floor} gives the pointwise lower bound
\[
H_{\sigma}(h)\ge\frac{49}{100}
\qquad\text{for every } h\in[-1,1].
\]
In particular, this bound applies to \(h=q\).
If \({c_R} \le  49/100\), then
\[
(f_{\sigma}(q)-L/A)_+\le H_{\sigma}(q)-{c_R}.
\]
The right side is nonnegative. Consequently
\[
J\le\Abar(R)B_{\mathrm{height}}(R,{\sigma})-L
\le\Abar(R)\left[1+\frac98(R-1/3)^2\right]-L
\le \Aeq-L=4\pi.
\]
The last step uses the inequality \((1-u)^{-1/2} \ge  1+u/2\).

If \({c_R} \ge  49/100\), the same reasoning gives
\[
J\le\Abar(R)(B_{\mathrm{height}}(R,{\sigma})-49/100).
\]
Put \(z=(R-1/3)^2\in[0,1/9]\). Then
\[
\frac{\Abar(R)}{2\pi^2}
\left[1+\frac98(R-1/3)^2-\frac{49}{100}\right]
=\frac2{\sqrt3}\sqrt{1-\frac94z}
\left(\frac{51}{100}+\frac98z\right).
\]
Its derivative has the sign of
\[
\frac98-\frac98\frac{51}{100}-\frac{243}{64}z
=\frac{441}{800}-\frac{243}{64}z>0
\quad (0\le z\le1/9).
\]
Therefore this expression is at most its value at \(R=2/3\), namely
\(127/200\). It follows that
\[
J\le\frac{127\pi^2}{100}<4\pi,
\]
where the last strict inequality follows from \(\pi<22/7\), since
\(127\cdot22=2794<2800\).

It remains to treat the range \(2/3\le R\le R_*\), for which we use
 the quadratic Hermite argument below.

\subsubsection{\texorpdfstring{Quadratic proof for $2/3\le R\le R_*$}{Quadratic proof for 2/3 <= R <= R_*}}\label{sec:small-centres-large-weights}
We use a quadratic Hermite majorant that avoids a larger collection of coefficient inequalities.
Assume first $0<s\le1/4$. Let $V=R/2\in[1/3,R_*/2]$. Denote by $H_{s,V}$ the quadratic Hermite polynomial tangent to $f_s$ at $-V$ and agreeing with $f_s$ at $1$. The positivity of $f_s'''$ makes $H_{s,V}\ge f_s$ on $[-1,1]$. Write $P=1-s^2$, $B=1+s^2$, $d=B+2sV$ and $e=B-2s=(1-s)^2$.
The interpolant has the form
\[
 H_{s,V}(h)=f_s(-V)+f_s'(-V)(h+V)+A_{s,V}(h+V)^2,
\]
where
\begin{align*}
 A_{s,V}&=\frac{f_s(1)-f_s(-V)-(1+V)f_s'(-V)}{(1+V)^2}\\
 &=\frac{4s^2 P^2(d+2e)}{e^2d^3}>0.
\end{align*}
For the second identity use $d-e=2s(1+V)$ and
$d^3-3de^2+2e^3=(d-e)^2(d+2e)$. Completing the square gives
\[
 \min_{h\in\mathbb R} H_{s,V}(h)
 =\frac{P^2}{d^2}-\frac{P^2e^2}{d^3(d+2e)}.
\]
Subtracting $1/2$ and using $d+2e=3s^2+2Vs-4s+3$ yields
\begin{equation}\label{eq:N}
 \min_{h\in\mathbb R}H_{s,V}(h)-\frac12
 =\frac{N(s,V)}{2(1+s^2+2Vs)^3(3s^2+2Vs-4s+3)},
\end{equation}
with
\[
\begin{aligned}
N(s,V)&=1+a_7(s+s^7)+a_6(s^2+s^6)+a_5(s^3+s^5)+a_4s^4+s^8,\\
a_7&=4-4V,\quad a_6=-20+8V-40V^2,\\
a_5&=12-76V+48V^2-48V^3,\\
a_4&=-10+80V-112V^2+32V^3-16V^4.
\end{aligned}
\]
For $1/6\le V\le19/50$ and $0\le s\le1/4$, the bounds
$a_6\le-98/5$, $a_5\le1$, $a_4\le23$, $a_7\le4$ follow from the
coefficient formulas (the first two decrease in $V$, while
$a_4\le-10+80(19/50)+32(19/50)^3<23$). Subscripts denote partial derivatives.
In particular,
\begin{align*}
 N_{ss}&=42a_7s^5+a_6(2+30s^4)+a_5(6s+20s^3)+12a_4s^2+56s^6\\
 &\le\frac{168}{4^5}-\frac{196}{5}+\frac64+\frac{20}{4^3}
       +\frac{276}{4^2}+\frac{56}{4^6}<-19.
\end{align*}
Furthermore
\[
N_V\le-4s+80s^4\le-\frac{11}{4}s,
\]
because $8-80V<0$, $-76+96V-144V^2<0$, and $a_4'\le80$.
Thus $N$ is concave in $s$ and nonincreasing in $V$, and
\begin{equation}\label{eq:Npos}
N(s,V)\ge\min\left\{1,N(1/4,19/50)\right\}
=\frac{173885049}{25600000000}>0.
\end{equation}
Therefore $H_{s,V}\ge1/2$.

For $R\ge2/3$, the threshold obeys
\[
{c_R}\ge {c_{2/3}}=\frac2{\sqrt3}-\frac2\pi
>\frac{15}{13}-\frac{20}{31}>\frac12,
\]
using $\sqrt3<26/15$ and $\pi>31/10$. Hence
$J\le A\bigl(\int H_{s,V}\dd\mu_{\weights}-1/2\bigr)$, where $A=A(\weights)$.
The coefficient $A_{s,V}$ is positive. Since $\Av_{\mu_{\weights}} q=0$ and
$\Av_{\mu_{\weights}} q^2\le V$, the expectation of the quadratic is bounded by
$f_s(-V)+Vf_s'(-V)+A_{s,V}(V+V^2)$. This is also its expectation under
the two-point probability measure assigning masses $1/(1+V)$ and
$V/(1+V)$ to $-V$ and $1$, respectively, because that measure has mean
zero and second moment $V$. At these two nodes $H_{s,V}=f_s$; hence
\[
 \int H_{s,V}\dd\mu_{\weights}\le F(s,V)
 :=\frac{f_s(-V)+V f_s(1)}{1+V}.
\]
An exact simplification gives
\begin{align*}
F(s,V)&=1+\frac{4s^2{Q_{\mathrm{mom}}}(s,V)}{(1-s)^2(1+s^2+2Vs)^2},\\
{Q_{\mathrm{mom}}}(s,V)&=(3V-1)(1+s^2)+(4V^2-2V+2)s.
\end{align*}
For $V\ge1/3$ all coefficients of ${Q_{\mathrm{mom}}}$ as a polynomial in $s$ are nonnegative. Moreover
\[
\frac\dd{\dd s}\bigl[(1-s)(1+s^2+2Vs)\bigr]
=2V-1+(2-4V)s-3s^2
\le-\frac6{25}+\frac16<0
\]
when $V\le19/50$ and $s\le1/4$. Thus $F(s,V)$ increases in $s\in[0,1/4]$, and
\[
F(s,R/2)\le M(R):=F(1/4,R/2)
=\frac{450/(17+4R)^2+25R/9}{2+R}.
\]
A direct differentiation yields
\[
M'(R)=\frac{800(4R+35)}{9(4R+17)^3}
\le M'(2/3)=\frac{90400}{205379}<\frac12.
\]
On the other hand
\begin{multline*}
\frac\dd{\dd R}\left(\frac12+\frac{2}{\pi\sqrt{(1-R)(1+3R)}}\right)\\
=\frac{6R-2}{\pi\bigl((1-R)(1+3R)\bigr)^{3/2}}\ge\frac2\pi>\frac12
\qquad (R\ge2/3).
\end{multline*}
The desired comparison therefore reduces to $R=2/3$, where
\[
 M(2/3)=\frac{35425}{31329}
 <\frac{25}{22}<\frac12+\frac2\pi.
\]
Since $M(R)\ge F(s,R/2)\ge1$ and $A(\weights)\le\Abar(R)$, it follows that
\[
J\le\Abar(R)\bigl(M(R)-1/2\bigr)<4\pi,
\]
proving the final range in Theorem~\ref{thm:small}. The endpoint $s=0$ follows by continuity, or directly from $J=(A-L)_+<4\pi$ for $R\ge2/3$.

\paragraph{Equality.}
For \({\sigma}>0\), the moment comparison in the first case for
$1/3\le R\le2/3$ is strict. For \({\sigma}=0\), the direct identity
\(J=(A-L)_+\le \Aeq-L\) gives equality exactly for the
 equal-weight torus. The second case and the quadratic argument
for $R\ge2/3$ are strict. This proves the equality statement of
Theorem~\ref{thm:small}.
\end{proof}

\subsection{Quartic estimates for small weights}\label{lw:quartic}

The quartic Hermite majorant gives the following estimates for small weights.

\begin{theorem}[Quartic bounds for small weights]\label{thm:quartic-small-weights}
For every fundamental embedding with $1/3\le R\le7/20$ and every
unit direction $p$, one has
\begin{align*}
 J_s(\weights,p)&\le4\pi-\frac{3\Aeq}{3400}
       &&\left(\frac14\le s\le\frac25\right),\\
 J_s(\weights,p)&\le4\pi-\frac{\Aeq}{14000}
       &&\left(\frac25\le s\le\frac{17}{40}\right).
\end{align*}
In particular, $J_s(\weights,p)<4\pi$ throughout $1/4\le s\le17/40$.
\end{theorem}

\begin{proof}
We use the quartic $H_{\sigma}$, its coefficients $h_j,u_j$, the denominator
${D_{\mathrm{height}}}(\sigma)$, and the moment bound $B_{\mathrm{height}}(R,\sigma)$ from
Section~\ref{sec:small-centres}. The four polynomial comparisons
needed below follow from one-dimensional monotonicity and convexity.

Let
\[
a=4/17,\qquad b=10/29,\qquad d=9/25.
\]
The required last endpoint \(\sigma(17/40)=680/1889\) is strictly smaller than \(d\), so proving the second floor and envelope on \([b,d]\) suffices.

\subsubsection*{A quadratic lower bound for both floors}

For \(k=17\) on \(a\le {\sigma}\le b\), or \(k=14\) on \(b\le {\sigma}\le d\), define
\[
G_{k,{\sigma}}(h)=H_{\sigma}(h)+\frac{(1-h)^2}{k}.
\]
We prove
\[
G_{k,{\sigma}}(h)\ge9/20\qquad(-1\le h\le1).
\]

\paragraph{\texorpdfstring{The right interval \(h\ge -1/4\)}{The right interval h >= -1/4}}

For each fixed \(h\), the derivative \(\partial_{\sigma} f_{\sigma}(h)=4(h-2{\sigma})/(1-2{\sigma}h)^3\) changes sign at most once, from positive to negative. Hence for \(0\le {\sigma}\le T\), where $T=10/29$ for $k=17$ and $T=9/25$ for $k=14$,
\[
f_{\sigma}(h)\ge\min\{1,f_T(h)\}.
\]
The function \(f_T(h)+(1-h)^2/k\) is strictly convex in \(h\); its derivative at \(h=-1/4\) is positive, and its value there exceeds \(9/20\), as shown by these exact rational calculations:

\begingroup
\renewcommand{\arraystretch}{2.2}
\[
\begin{array}{c|cc}
(k,T)&\text{gap at }-1/4&\text{derivative at }-1/4\\
(17,10/29) & \dfrac{541}{23120} & \dfrac{2965}{9826}\\
(14,9/25) & \dfrac{29181}{3898720} & \dfrac{1400369}{5750612}\\
\end{array}
\]
\endgroup

It follows that \(G\ge 9/20\) on \(h\ge -1/4\), since \(H\ge f\) and \(1+(1-h)^2/k>9/20\).

\paragraph{\texorpdfstring{The left interval \(h\le -1/4\)}{The left interval h <= -1/4}}

The ratio of Newton coefficients is
\[
\frac{u_3}{u_4}=\frac{(20-{\sigma})(1-2{\sigma})^2}{{\sigma}(2{\sigma}^2-83{\sigma}+50)}\le1
\qquad(a\le {\sigma}\le d).
\]
Indeed the inequality is equivalent to
\[
6{\sigma}^3-167{\sigma}^2+131{\sigma}-20\ge0.
\]
This polynomial is strictly increasing on \([a,d]\), because its derivative is at least
\(131-334d>0\); its value at \(a\) is \(8136/4913>0\).

For \(z=h-1/4\le -1/2\), the Newton quadratic \(q_{\sigma}(z)=u_2+u_3z+u_4z^2\) is therefore decreasing in \(z\). Consequently
\[
q_{\sigma}(z)\ge q_{\sigma}(-1/2)=A({\sigma})
=\frac{400{\sigma}^2(1+2{\sigma})p_{\mathrm{floor}}({\sigma})}{{D_{\mathrm{height}}}({\sigma})},
\]
where
\[
p_{\mathrm{floor}}({\sigma})=22{\sigma}^4-347{\sigma}^3+419{\sigma}^2-181{\sigma}+30>0.
\]
To check positivity use \({\sigma}\le d<3/8\):
\[
p_{\mathrm{floor}}({\sigma})\ge\frac{2311}{8}{\sigma}^2-181{\sigma}+30>0,
\]
whose discriminant is \(181^2-4(2311/8)30=-1904<0\).

Writing \(x=-1/5\), we obtain
\[
G_{k,{\sigma}}(h)\ge
f_{\sigma}(x)+f_{\sigma}'(x)(h-x)+A({\sigma})(h-x)^2+\frac{(1-h)^2}{k}.
\]
The minimum of the right-hand quadratic is exactly
\[
\frac9{20}+
\frac{N_k({\sigma})}{20(5+2{\sigma})^3\{{D_{\mathrm{height}}}({\sigma})+400k{\sigma}^2(1+2{\sigma})p_{\mathrm{floor}}({\sigma})\}}.
\]
The denominator is strictly positive, and the explicit numerators are
\[
\begin{aligned}
N_{17}({\sigma})={}&-1895120512{\sigma}^{10}+19841322688{\sigma}^9+24272581600{\sigma}^8\\
&-28192527952{\sigma}^7-6007153048{\sigma}^6+8521122980{\sigma}^5\\
&-579055750{\sigma}^4-755899375{\sigma}^3+116406250{\sigma}^2\\
&+8487500{\sigma}+1375000,
\end{aligned}
\]
\[
\begin{aligned}
N_{14}({\sigma})={}&-1560118912{\sigma}^{10}+16334173888{\sigma}^9+19930828000{\sigma}^8\\
&-23348675152{\sigma}^7-4974951448{\sigma}^6+7142346980{\sigma}^5\\
&-471985750{\sigma}^4-648649375{\sigma}^3+96906250{\sigma}^2\\
&+8487500{\sigma}+1375000.
\end{aligned}
\]
The following elementary derivative bounds prove both are decreasing on their respective intervals.

\paragraph{First numerator}

Put \(z={\sigma}-a\), so \(0\le z\le b-a<11/100\). Exact expansion gives
\[
\begin{aligned}
N_{17}'(a+z)={}&-18951205120z^9+\frac{2353478986944}{17}z^8
+\frac{142345430390528}{289}z^7\\
&+\frac{1859862557261520}{4913}z^6
+\frac{2844454946285040}{83521}z^5\\
&-\frac{54749773922651340}{1419857}z^4
-\frac{140329823188575704}{24137569}z^3\\
&+\frac{99977586127302627}{410338673}z^2
-\frac{831828279097744020}{6975757441}z
-\frac{1417899634161866100}{118587876497}.
\end{aligned}
\]
Increasing the coefficients of positive terms and dropping the negative cubic and ninth-degree terms yields
\[
\begin{aligned}
N_{17}'(a+z)\le{}&-11\cdot10^6+10^6z(-119+244z)\\
&+10^9z^4(-38+35z+379z^2+493z^3+139z^4)<0.
\end{aligned}
\]
Both bracketed polynomials are negative when \(z\le 11/100\). For the second,
\[
-38+35\frac{11}{100}+379\left(\frac{11}{100}\right)^2
+493\left(\frac{11}{100}\right)^3
+139\left(\frac{11}{100}\right)^4
=-\frac{2888756601}{100000000}<0.
\]
Thus
\[
N_{17}({\sigma})\ge N_{17}(b)
=\frac{100370715218568000000}{420707233300201}>0.
\]

\paragraph{Second numerator}

Put \(z={\sigma}-b\), so \(0\le z\le d-b<1/60\). Exact expansion gives
\[
\begin{aligned}
N_{14}'(b+z)={}&-15601189120z^9+\frac{2859112363968}{29}z^8
+\frac{418987880733440}{841}z^7\\
&+\frac{16026981278175504}{24389}z^6
+\frac{240399986013156912}{707281}z^5\\
&+\frac{1178638286957789700}{20511149}z^4
-\frac{3162050628857147000}{594823321}z^3\\
&-\frac{41456211829205180625}{17249876309}z^2
-\frac{162957531744980100000}{500246412961}z\\
&-\frac{407874526536351600000}{14507145975869}.
\end{aligned}
\]
It follows, on dropping the negative linear, quadratic, and ninth-degree terms, that
\[
N_{14}'(b+z)
\le-28\cdot10^6+
10^9z^3(-5+58z+340z^2+658z^3+499z^4+99z^5)<0.
\]
The bracket increases in \(z\ge 0\) and is still negative at \(z=1/60\) (indeed it is less than \(-3\)). Therefore
\[
N_{14}({\sigma})\ge N_{14}(d)
=\frac{7995624685216792213}{95367431640625}>0.
\]
This completes the proof of both floor inequalities.

\subsubsection*{Convexity gives both moment-envelope bounds}

For completeness, the coefficient signs up to \({\sigma}=3/8\) follow directly
from the monomial coefficient formulas for $H_{\sigma}$ in
Section~\ref{sec:small-centres}: the two numerator polynomials
\[
p_2=294{\sigma}^4-441{\sigma}^3+4630{\sigma}^2-4300{\sigma}+1000,
\quad p_3=-42{\sigma}^3+923{\sigma}^2-860{\sigma}+200
\]
are decreasing there, since
\[
p_2'\le1176(3/8)^3+9260(3/8)-4300=-48991/64<0,
\quad p_3'\le1846(3/8)-860=-671/4<0.
\]
Their endpoint values are \(43319/2048\) and \(1301/256\), respectively, both positive. The fourth-degree coefficient's numerator obeys \(2{\sigma}^2-83{\sigma}+50\ge 151/8>0\). Thus \(h_2,h_3,h_4>0\) on \(0<{\sigma}\le 3/8\), and \(B_{\mathrm{height}}(R,{\sigma})\) increases in \(R\ge 0\). It suffices to bound
\[
\mathfrak b(\sigma):=B_{\mathrm{height}}(7/20,\sigma)
=-\frac{2\sigma+1}{5{D_{\mathrm{height}}}(\sigma)}\,p_6(\sigma),
\]
where
\[
\begin{aligned}
p_6(\sigma)={}&308\sigma^6-11582\sigma^5-5830\sigma^4+25155\sigma^3\\
&-33650\sigma^2+21500\sigma-5000.
\end{aligned}
\]

\paragraph{\texorpdfstring{Convexity on \([a,d]\)}{Convexity on [a,d]}}

Direct differentiation gives
\[
\mathfrak b''({\sigma})=\frac{M({\sigma})}{5(2-{\sigma})^5(1-2{\sigma})^3(5+2{\sigma})^5},
\]
where
\[
\begin{aligned}
M({\sigma})=4\bigl(&187776{\sigma}^{10}-39504{\sigma}^9+9974640{\sigma}^8+3752096{\sigma}^7\\
&+9206756{\sigma}^6+21404901{\sigma}^5-72918050{\sigma}^4\\
&+60601000{\sigma}^3-23760000{\sigma}^2+5187500{\sigma}-475000\bigr).
\end{aligned}
\]
The denominator is positive. Write \(z={\sigma}-a\), so \(0\le z\le d-a<1/8\) and
\(M(a+z)=\sum_{j=0}^{10}c_jz^j\). The exact coefficients are:

\begingroup
\renewcommand{\arraystretch}{2.2}
\[
\begin{array}{c|r}
j&c_j\\
0 & \dfrac{107466284661492960}{2015993900449}\\
1 & \dfrac{307166160092613552}{118587876497}\\
2 & \dfrac{-52704781096106880}{6975757441}\\
3 & \dfrac{11533931614456992}{410338673}\\
4 & \dfrac{-3498608404791288}{24137569}\\
5 & \dfrac{261596014409556}{1419857}\\
6 & \dfrac{10332127521168}{83521}\\
7 & \dfrac{446939259520}{4913}\\
8 & \dfrac{11974772928}{289}\\
9 & \dfrac{27357888}{17}\\
10 & 751104\\
\end{array}
\]
\endgroup

All coefficients of degree at least five are positive, and the lower coefficients satisfy
\[
\begin{gathered}
c_0>53000,\quad c_1>2590000,\quad c_2>-7556000,\\
c_3>28108000,\quad c_4>-144945000.
\end{gathered}
\]
Consequently
\[
M(a+z)>53000+z\left(2590000-\frac{7556000}{8}\right)
+z^3\left(28108000-\frac{144945000}{8}\right)>0.
\]
This proves strict convexity.

\paragraph{Endpoint bounds}

A convex function on a closed interval is bounded above by the larger endpoint value. Here
\[
\mathfrak b(a)=\frac{167916653}{180980325}<\frac{93}{100},
\qquad
\mathfrak b(b)=\frac{37478189}{41399424}<\frac{93}{100}.
\]
For the second interval,
\[
\frac{229}{250}-\mathfrak b(b)=\frac{55460423}{5174928000}>0,
\quad
\frac{229}{250}-\mathfrak b(d)=\frac{414489641091}{352693723632250}>0.
\]
Therefore
\[
B_{\mathrm{height}}(R,{\sigma})\le93/100\quad(a\le {\sigma}\le b),
\qquad B_{\mathrm{height}}(R,{\sigma})\le229/250\quad(b\le {\sigma}\le d),
\]
for every \(1/3\le R\le 7/20\).

\subsubsection*{Integration and the spectral-area gap}

For $R\le7/20$ one has ${c_R}<9/20$. The nonnegative polynomial
\[
Q_{\sigma}(h)=H_{\sigma}(h)-{c_R}+(1-h)^2/k
\]
majorizes $(f_{\sigma}(h)-{c_R})_+$. Using $\mathbb E q=0$ and $\mathbb E q^2\le R/2$ from
\eqref{eq:global-moments}, integration gives
\[
J\le\Abar(R)\{B_{\mathrm{height}}(R,{\sigma})+k^{-1}(1+R/2)\}-L.
\]
Here the replacement of $A(\weights)$ by $\Abar(R)$ is legitimate because the expectation of $Q_{\sigma}$ is nonnegative. Using $\Abar(R)\le \Aeq$ from \eqref{eq:global-area-envelope}
and $L=\Aeq-4\pi$, the endpoint comparisons just proved yield
\[
J\le4\pi-\frac{3\Aeq}{3400}
\quad(1/4\le s\le2/5),
\]
\[
J\le4\pi-\frac{\Aeq}{14000}
\quad(2/5\le s\le17/40).
\]
These are the two asserted bounds.
\end{proof}

\subsection{Radial estimates for small weights}\label{mid:analytic}

\begin{theorem}\label{mid:strip}
For every fundamental torus with $1/3\le R\le7/20$, every unit ambient
direction $p$, and every $17/40\le s\le3/5$, the truncated area satisfies
$J_s(\weights,p)<4\pi$.
\end{theorem}

To prove Theorem~\ref{mid:strip}, we first establish preparatory estimates
using Lemma~\ref{lem:radial-cap-density}. They reduce the possible
height maxima to three scalar integrals with fixed endpoints.

\subsubsection*{Definitions and the radial upper bound}

Set
\[
{t_{\mathrm c}}=3/5,\quad{b_{\mathrm q}}=7/24,\quad c_0=4/9,\quad
K_0=\frac{363}{200}\frac{749}{2880}.
\]
For \(s\) in the above range put
\[
P=1-s^2,\quad B=1+s^2,\quad
f_s(q)=\frac{P^2}{(B-2sq)^2},\quad
C(s)=\frac{f_s({t_{\mathrm c}})-c_0}{({t_{\mathrm c}}+{b_{\mathrm q}})^2},
\]
\[
Q^{\mathrm{quad}}_s(q)=C(s)(q+{b_{\mathrm q}})^2,
\qquad\psi_s(q)=f_s(q)-c_0-Q^{\mathrm{quad}}_s(q)\quad(q\ge{t_{\mathrm c}}).
\]
The actual truncation satisfies $L/A(\weights)\ge {c_R}\ge1-\sqrt3/\pi>4/9$, the last elementary comparison follows from $\sqrt3<26/15$ and $\pi>157/50$. Thus it suffices to bound the positive part with cutoff $c_0$. Lemma~\ref{lem:quadratic-cap-majorant} shows that \(Q^{\mathrm{quad}}_s\) majorizes \((f_s-c_0)_+\) below \({t_{\mathrm c}}\), and that \(\psi_s\ge 0\) above \({t_{\mathrm c}}\). The same argument is valid on this smaller \(s\) range. It uses only
\[
1+s^2+(7/12)s\ge(49/30)(1-s^2),
\]
whose left minus right is increasing and equals \(4331/48000>0\) at \(s=17/40\). Also \(q_0(s,c_0)+{b_{\mathrm q}}>0\) and \(q_0(s,c_0)<{t_{\mathrm c}}\), where \(q_0(s,c_0)=(5s^2-1)/(4s)\) by \eqref{eq:truncation-height}.

The second moment gives
\[
\frac1{4\pi}\int Q^{\mathrm{quad}}_s(q)\,dA_{\weights}\le K_0C(s).
\]
For \(m=\max q\le {t_{\mathrm c}}\) this already proves the desired conclusion, since
\[
K_0C(s)\le K_0C(17/40)
=\frac{9201107531}{14783064360}<\frac{623}{1000}<1.
\]
We henceforth suppose \(m>{t_{\mathrm c}}\). In particular \(m>\sqrt{R}\), so Lemma~\ref{lem:radial-cap-density} applies.

For
\[
e_R(m)=(2R-1)m+2\sqrt{R(1-R)}\sqrt{1-m^2},
\]
the second inequality in \eqref{radial:oscillator-density}, divided by $4\pi$,
gives the density bound \((1-R)/(q-e_R(m))\). This kernel increases with \(R\) when \(m\ge q\):
\[
\frac{\partial}{\partial R}\frac{1-R}{q-e_R(m)}
=\frac{m-q+\sqrt{(1-R)/R}\sqrt{1-m^2}}{(q-e_R(m))^2}\ge0.
\]
Since \(2\sqrt{91}/20 <24/25\), define
\[
\widehat e(m)=-\frac3{10}m+\frac{24}{25}\sqrt{1-m^2}
\ge e_{7/20}(m).
\]
It decreases on \([{t_{\mathrm c}},1]\) and
\(\widehat e({t_{\mathrm c}})=147/250<{t_{\mathrm c}}\), so all denominators below are positive. The radial density bound \eqref{radial:oscillator-density} and the split into
\(Q^{\mathrm{quad}}_s\) plus cap excess give
\begin{equation}
\frac{J}{4\pi}\le K_0C(s)+\frac{13}{20}
\int_{t_{\mathrm c}}^m\frac{\psi_s(q)}{q-\widehat e(m)}\,dq.
\label{eq:radial-excess-bound}
\end{equation}

\subsubsection*{\texorpdfstring{Convexity of the excess and monotonicity in \(s\)}{Convexity of the excess and monotonicity in s}}

We use convexity in the height and monotonicity in the conformal
radius to reduce the radial estimate to fixed endpoints.

\paragraph{Convexity in height}

We have \(\psi_s({t_{\mathrm c}})=0\), \(\psi_s\ge 0\) and \(\psi_s\) is strictly convex on \([{t_{\mathrm c}},1]\). Indeed
\[
f_s''(q)=\frac{24s^2P^2}{(B-2sq)^4}
\ge\frac{24s^2P^2}{(B-2s{t_{\mathrm c}})^4}.
\]
Writing \(T=B-2s{t_{\mathrm c}}\), the inequalities
\[
T\le2s,\qquad12({t_{\mathrm c}}+{b_{\mathrm q}})^2>4
\]
imply
\(f_s''({t_{\mathrm c}})>2f_s({t_{\mathrm c}})/({t_{\mathrm c}}+{b_{\mathrm q}})^2>2C(s)\).
The first inequality holds because \(1+s^2-(16/5)s\) decreases on the interval and is negative at \(17/40\). Thus \(\psi_s''>0\).

\paragraph{Monotonicity of the excess in the conformal radius}

The function \(C(s)\) decreases, and \(\psi_s(q)\) increases with \(s\) for each \(q\ge {t_{\mathrm c}}\). Here is an elementary proof of the second, less immediate assertion.

Put \({\sigma}=s/(1+s^2)\), so
\[
1/3<{\sigma}\le15/34,\qquad
f_{\sigma}(q)=\frac{1-4{\sigma}^2}{(1-2{\sigma}q)^2},\qquad
\partial_{\sigma} f_{\sigma}(q)=\frac{4(q-2{\sigma})}{(1-2{\sigma}q)^3}.
\]
Regard $C$ and $\psi_s$ as functions of $\sigma$ through the inverse
$s=s(\sigma)$ of the increasing map $s\mapsto\sigma$. Since
\({t_{\mathrm c}}<2{\sigma}\), $\frac{d}{d\sigma}C(s(\sigma))<0$.
For \(q\ge 2{\sigma}\), the assertion
$\frac{d}{d\sigma}\psi_{s(\sigma)}(q)\ge0$ is immediate. For \({t_{\mathrm c}}\le q<2{\sigma}\), consider
\[
g_{\sigma}(q)=\frac{2{\sigma}-q}{(1-2{\sigma}q)^3(q+{b_{\mathrm q}})^2}.
\]
The sign of \(g_{\sigma}'\) is that of
\[
N({\sigma},q)=-8{\sigma}q^2+(20{\sigma}^2-7{\sigma}/6+1)q+7{\sigma}^2/2-4{\sigma}-7/24.
\]
Its maximum as a quadratic in the whole real \(q\)-axis is
\[
\max_q N({\sigma},q)=\frac{P_{\mathrm{rad}}({\sigma})}{1152{\sigma}},
\quad
P_{\mathrm{rad}}({\sigma})=14400{\sigma}^4+2352{\sigma}^3-3119{\sigma}^2-420{\sigma}+36.
\]
This numerator is negative on \([1/3,15/34]\): it is convex there because
\[
P_{\mathrm{rad}}''({\sigma})=172800{\sigma}^2+14112{\sigma}-6238\ge17666>0,
\]
and both endpoint values are negative,
\[
P_{\mathrm{rad}}(1/3)=-557/3,
\qquad
P_{\mathrm{rad}}(15/34)=-2966751/334084.
\]
Consequently \(g_{\sigma}(q)\le g_{\sigma}({t_{\mathrm c}})\). Written out, this says precisely
\[
\frac{2{\sigma}-q}{(1-2{\sigma}q)^3}
\le\frac{2{\sigma}-{t_{\mathrm c}}}{(1-2{\sigma}{t_{\mathrm c}})^3}
\frac{(q+{b_{\mathrm q}})^2}{({t_{\mathrm c}}+{b_{\mathrm q}})^2},
\]
which is equivalent to $\frac{d}{d\sigma}\psi_{s(\sigma)}(q)\ge0$. Since \({\sigma}\) increases with \(s<1\), the desired monotonicity follows.

\subsubsection*{Three height ranges reduce the problem to fixed integrals}

For any \(b>0\), Cauchy--Schwarz gives
\[
\widehat e(m)+bm
\le\sqrt{(b-3/10)^2+(24/25)^2}.
\]
We use these three rational choices:

\begingroup
\renewcommand{\arraystretch}{2.2}
\[
\begin{array}{c|ccc}
\text{height interval}&b_i&a_i&a_i^2-(b_i-3/10)^2-(24/25)^2\\
{[3/5,4/5]} & 1 & \dfrac{119}{100} & \dfrac{9}{2000}\\
{[4/5,19/20]} & 2 & \dfrac{1953}{1000} & \dfrac{2609}{1000000}\\
{[19/20,1]} & \dfrac{10}{3} & \dfrac{191}{60} & \dfrac{109}{10000}\\
\end{array}
\]
\endgroup

Thus \(\widehat e(m)\le a_i-b_im\) on every corresponding interval. Define
\[
I_i(s,m)=\int_{t_{\mathrm c}}^m\frac{\psi_s(q)}{q-a_i+b_i m}\,dq.
\]
All its denominators are positive on the specified height range. We claim \(I_i(s,m)\) is nondecreasing in \(m\) on that range.

To see this put
\[
{\ell_m}=m-{t_{\mathrm c}},\quad {D_m}={t_{\mathrm c}}+b_i m-a_i>0,\quad {E_m}={D_m}+{\ell_m},
\quad x={\ell_m}/{D_m}.
\]
Convexity of \(\psi_s\) and its zero at \({t_{\mathrm c}}\) give
\(\psi_s(q)\le\psi_s(m)(q-{t_{\mathrm c}})/{\ell_m}\). Differentiation and elementary integration imply
\[
\begin{aligned}
\partial_m I_i
&=\frac{\psi_s(m)}{E_m}-b_i\int_{t_{\mathrm c}}^m
\frac{\psi_s(q)}{(q-a_i+b_i m)^2}\,dq\\
&\ge\frac{\psi_s(m)}{E_m}
\left[1+b_i-b_i\frac{1+x}{x}\log(1+x)\right]\\
&\ge\frac{\psi_s(m)}{E_m}\left(1-\frac{b_i x}{2}\right).
\end{aligned}
\]
Here
\[
\log(1+x)\le\frac{x(2+x)}{2(1+x)}\quad(x\ge0),
\]
as follows by differentiation. The relevant upper bounds for \(x\) are
\[
x\le20/21<2\quad(i=1),\qquad
x\le200/247<1\quad(i=2),\qquad
x\le3/5=2/(10/3)\quad(i=3).
\]
For \(i=1\), \(x\) increases with \(m\); for \(i=2,3\), it decreases. These facts follow directly from
\(x'=((1+b_i){t_{\mathrm c}}-a_i)/{D_m}^2\). Thus the derivative of every \(I_i\) is nonnegative.

Inequality \eqref{eq:radial-excess-bound} is therefore at most one of the three fixed-endpoint functions
\begin{equation}
V_i(s)=K_0C(s)+\frac{13}{20}\int_{t_{\mathrm c}}^{M_i}
\frac{\psi_s(q)}{q-E_i}\,dq,
\label{eq:fixed-endpoint-functions}
\end{equation}
where
\[
(M_1,E_1)=(4/5,39/100),\quad
(M_2,E_2)=(19/20,53/1000),\quad
(M_3,E_3)=(1,-3/20).
\]
It remains to bound these three one-variable functions below \(1\).

\subsubsection*{The first fixed integral}

By monotonicity proved above, \(C(s)\le C(17/40)\) and
\(\psi_s(q)\le \psi_{3/5}(q)\). The convex chord bound on \([{t_{\mathrm c}},M_1]\) gives
\[
\int_{t_{\mathrm c}}^{M_1}\frac{\psi_{3/5}(q)}{q-E_1}\,dq
\le\psi_{3/5}(M_1)\frac{M_1-{t_{\mathrm c}}}{M_1+{t_{\mathrm c}}-2E_1}
=\frac{10}{31}\psi_{3/5}(4/5).
\]
Indeed the exact integral of the chord is
\(\psi_{3/5}(M_1)[1-({t_{\mathrm c}}-E_1)(M_1-{t_{\mathrm c}})^{-1}\log((M_1-E_1)/({t_{\mathrm c}}-E_1))]\),
and \(\log((1+u)/(1-u))\ge 2u\) yields the displayed upper bound.
All values are rational:
\[
\psi_{3/5}(4/5)=\frac{1101533}{858675},
\]
\[
V_1(s)\le
\frac{218548166544199}{245177122410600}<\frac9{10}<1.
\]

\subsubsection*{The second fixed integral}

Put \(M=19/20,E=53/1000\). We use the three radius intervals
\[
[17/40,1/2],\quad[1/2,11/20],\quad[11/20,3/5].
\]
On a generic interval \([s_-,s_+]\), monotonicity gives
\[
V_2(s)\le K_0C(s_-)+\frac{13}{20}\int_{t_{\mathrm c}}^M
\frac{\psi_{s_+}(q)}{q-E}\,dq.
\]
The three remaining integrals have simple rational upper bounds. For a fixed \(s=s_+\), define
\[
D=B-2sE,\ T=B-2s{t_{\mathrm c}},\ U=B-2sM,
\quad K=P^2/D^2,
\]
\[
A_0=K-c_0-C(s)(E+{b_{\mathrm q}})^2<0,
\quad u=\frac{M-{t_{\mathrm c}}}{M+{t_{\mathrm c}}-2E},
\quad v=\frac{s(M-{t_{\mathrm c}})}{B-s(M+{t_{\mathrm c}})}.
\]
The sign of \(A_0\) follows already from \(f_s(E)<4/9\) when \(s\ge 1/2\); at \(s=1/2\) its value before the subtraction of \(Q^{\mathrm{quad}}_s(E)\) is \(62500/159201<4/9\), and it decreases with \(s\).

An explicit antiderivative yields
\[
\begin{aligned}
\mathcal I(s):=\int_{t_{\mathrm c}}^M\frac{\psi_s(q)}{q-E}\,dq
={}&A_0\log\frac{M-E}{{t_{\mathrm c}}-E}+K\log\frac{T}{U}
+\frac{P^2}{D}\left(\frac1U-\frac1T\right)\\
&-C(s)\left[\frac{M^2-{t_{\mathrm c}}^2}{2}+(E+2{b_{\mathrm q}})(M-{t_{\mathrm c}})\right].
\end{aligned}
\]
Use
\[
\log\frac{M-E}{{t_{\mathrm c}}-E}\ge2u,
\qquad
\log\frac{T}{U}\le2\left(v+\frac{v^3}{3(1-v^2)}\right).
\]
Both are immediate from the positive power series for the logarithm. Denote the resulting rational upper bound for \(\mathcal I(s)\) by \(U_2(s)\). Direct substitution gives:

\begingroup
\renewcommand{\arraystretch}{2.2}
\[
\begin{array}{c|c|c}
s_+&U_2(s_+)&\text{strict upper bound}\\
\dfrac{1}{2} & \dfrac{10364714935750}{20424605770881} & \dfrac{51}{100}\\
\dfrac{11}{20} & \dfrac{215297340652147981325}{343634121299686565022} & \dfrac{63}{100}\\
\dfrac{3}{5} & \dfrac{908564378470700}{1257242735138379} & \dfrac{73}{100}\\
\end{array}
\]
\endgroup

For the first terms,
\[
\begin{aligned}
K_0C(17/40)&=9201107531/14783064360<623/1000,\\
K_0C(1/2)&=1142603/2169960<527/1000,\\
K_0C(11/20)&=369607931/848069160<436/1000.
\end{aligned}
\]
Thus on the three radius intervals respectively,
\[
V_2<\frac{623}{1000}+\frac{13}{20}\frac{51}{100}=\frac{1909}{2000}<1,
\]
\[
V_2<\frac{527}{1000}+\frac{13}{20}\frac{63}{100}=\frac{1873}{2000}<1,
\]
\[
V_2<\frac{436}{1000}+\frac{13}{20}\frac{73}{100}=\frac{1821}{2000}<1.
\]

\subsubsection{The terminal scalar cap integral}
\label{sec:terminal-cap-integral}

Set
\[
 a=t_{\mathrm c}=\frac35,\qquad e=\frac3{20},\qquad {b_{\mathrm q}}=\frac7{24},\qquad
 d=a+{b_{\mathrm q}}=\frac{107}{120},\qquad
 K=\frac{363\cdot749}{200\cdot2880}=\frac{90629}{192000},
\]
and write
\[
 f_s(q)=\frac{(1-s^2)^2}{(1+s^2-2sq)^2},\qquad
 C(s)=\frac{f_s(a)-4/9}{d^2}.
\]
The third fixed-endpoint function in \eqref{eq:fixed-endpoint-functions} is
\begin{equation}\label{eq:terminal-cap-integral}
 V_3(s)=K C(s)+\frac{13}{20}\int_a^1
 \frac{f_s(q)-4/9-C(s)(q+{b_{\mathrm q}})^2}{q+e}\,dq.
\end{equation}

\begin{lemma}\label{lem:terminal-cap-integral}
The function $V_3$ is strictly increasing on $[1/3,3/5]$. In particular,
for $17/40\le s\le3/5$,
\[
 V_3(s)\le V_3(3/5)
 <\frac{461457834793}{464306867640}<1.
\]
\end{lemma}

\begin{proof}
Polynomial division gives
\[
 I:=\int_a^1\frac{(q+{b_{\mathrm q}})^2}{q+e}\,dq
 =\frac{37}{75}+\frac{289}{14400}\log\frac{23}{15}.
\]
Put ${\mathfrak c}=[K-(13/20)I]/d^2$. The elementary inequality
\[
 \log x>\frac{2(x-1)}{x+1}\qquad(x>1)
\]
follows from the positive power series for $2\operatorname{arctanh}$.
It gives $\log(23/15)>8/19$, and hence
\begin{equation}\label{eq:terminal-cap-A-bound}
 {\mathfrak c}<\frac{319273}{1740248}<\frac3{16},\qquad
 \frac3{16}-\frac{319273}{1740248}
 =\frac{14047}{3480496}>0.
\end{equation}

Set $\beta=2s/(1+s^2)$ and $u=1-\beta^2$. On $1/3\le s\le3/5$,
one has $a\le\beta\le15/17$, and
\[
 \frac{\partial f_s(q)}{\partial s}
 =\frac{4(1-s^2)}{(1+s^2)^2}
       \frac{q-\beta}{(1-\beta q)^3}.
\]
Thus $V_3'(s)$ is the positive factor $4(1-s^2)/(1+s^2)^2$ times
\[
 \Xi(\beta)=\frac{13}{20}\int_a^1
 \frac{q-\beta}{(q+e)(1-\beta q)^3}\,dq
 +{\mathfrak c}\frac{a-\beta}{(1-a\beta)^3}.
\]
Let
\[
 y=\frac{\beta-a}{1-a\beta}\in[0,3/5].
\]
Then $u/(1-a\beta)^2=(1-y^2)/(1-a^2)$. The substitution
$x=(q-\beta)/(1-\beta q)$ yields
\[
 \int_\beta^1\frac{q-\beta}{(1-\beta q)^3}\,dq=\frac1{2u},\qquad
 \int_a^\beta\frac{q-\beta}{(1-\beta q)^3}\,dq=-\frac{y^2}{2u}.
\]
On the positive part of the integral use $1/(q+e)\ge1/(1+e)$;
on its negative part use $1/(q+e)\le1/(a+e)$. The latter comparison
has the indicated direction because the remaining factor is negative.
Using \eqref{eq:terminal-cap-A-bound}, we obtain
\[
 2u\Xi(\beta)\ge
 \frac{13}{20}\left(\frac1{1+e}-\frac{y^2}{a+e}\right)
 -\frac{2{\mathfrak c}}{1-a^2}y(1-y^2)
 \ge B(y),
\]
where
\[
 B(y)=\frac{13}{23}-\frac{13}{15}y^2
             -\frac{75}{128}y(1-y^2).
\]
For $0\le y\le3/5$,
\[
 B'(y)=\frac{225}{128}y^2-\frac{26}{15}y-\frac{75}{128}
 \le-\frac{1303}{1920}y-\frac{75}{128}<0.
\]
Consequently
\[
 B(y)\ge B(3/5)=\frac{649}{23000}>0,
\]
which proves strict monotonicity.

At $s=3/5$ one has
\[
 f_{3/5}(q)=\frac{64}{(17-15q)^2},\qquad
 f_{3/5}(a)=1,\qquad C(3/5)=C_{\mathrm{end}}:=\frac{8000}{11449}.
\]
Partial fractions give
\[
 \int_a^1\frac{f_{3/5}(q)}{q+3/20}\,dq
 =\frac{96}{77}+\frac{1024}{5929}\log\frac{92}{15}.
\]
Writing $r_*=1024/5929$ and $\ell_*=\log(23/15)$, and using
$\log(92/15)=\log4+\ell_*$, we find
\[
 \begin{aligned}
 V_3(3/5)=KC_{\mathrm{end}}+\frac{13}{20}\biggl[
 &\frac{96}{77}-\frac{37C_{\mathrm{end}}}{75}+r_*\log4\\
 &+\left(r_*-\frac49-\frac{289C_{\mathrm{end}}}{14400}\right)\ell_*\biggr].
 \end{aligned}
\]
The coefficient of $\ell_*$ is
\[
 r_*-\frac49-\frac{289C_{\mathrm{end}}}{14400}
 =-\frac{19397545}{67881121}<0.
\]
Finally,
\[
 \begin{aligned}
 \log2
 &=2\left[\frac13+\sum_{j\ge1}\frac{3^{-(2j+1)}}{2j+1}\right]\\
 &<2\left[\frac13+\frac13\frac{1/27}{1-1/9}\right]
 =\frac{25}{36},
 \end{aligned}
\]
so $\log4<25/18$. Combining this with $\ell_*>8/19$ gives
\[
 \begin{aligned}
 V_3(3/5)
 &<KC_{\mathrm{end}}+\frac{13}{20}\left[
 \frac{96}{77}-\frac{37C_{\mathrm{end}}}{75}
 +\frac{1024}{5929}\frac{25}{18}
 -\frac{19397545}{67881121}\frac8{19}\right]\\
 &=\frac{461457834793}{464306867640}<1.
 \end{aligned}
\]
The final positive gap is $2849032847/464306867640$.
\end{proof}

\begin{proof}[Proof of Theorem~\ref{mid:strip}]
For $m\le{t_{\mathrm c}}$, the quadratic moment bound gives $J/(4\pi)<623/1000$.
For $m>{t_{\mathrm c}}$, \eqref{eq:radial-excess-bound} and the reduction to
\eqref{eq:fixed-endpoint-functions} bound $J/(4\pi)$ by one of
$V_1,V_2,V_3$. The first is less than $9/10$, the second is less than
$1909/2000$, and Lemma~\ref{lem:terminal-cap-integral} gives
\[
V_3(s)<\frac{461457834793}{464306867640}<1.
\]
This proves Theorem~\ref{mid:strip}.
\end{proof}

Together with Theorems~\ref{thm:small} and~\ref{thm:quartic-small-weights},
and Theorem~\ref{thm:central-small-weights} for \(s\ge3/5\),
this covers every \(0\le s<1\) when \(1/3\le R\le7/20\).

\subsection{\texorpdfstring{The interval \(1/4\le s\le3/10\)}{The interval 1/4 <= s <=3/10}}\label{sec:bottom-strip}

A scalar lower bound for the Hermite polynomial and five weight
comparisons yield a uniform estimate at small radii.

\begin{theorem}[Small-radius bound]\label{bs:theorem}
For every fundamental embedding with $1/3\le R\le R_*$, every unit
direction $p$, and every $0\le s\le3/10$, one has
\[
 J_s(\weights,p)\le4\pi.
\]
Equality holds if and only if $s=0$ and $r_1=r_2=r_3=1/3$.
Moreover,
\[
 \mathcal P(R,s,\alpha)<4\pi
 \qquad\left(\frac9{20}\le R\le R_*,\quad
             \frac14\le s\le\frac3{10},\quad0\le\alpha\le1\right).
\]
\end{theorem}

\paragraph{Preparatory Hermite estimates.}
For the range $1/4\le s\le3/10$ we first establish the polynomial
identities and the lower bound used in the proof of
Theorem~\ref{bs:theorem}.

Put \[
a=4/17,\quad b=30/109,\quad z={\sigma}-a,
\qquad 0\le z\le b-a<1/25,
\] where \({\sigma}=s/(1+s^2)\). Put \[
 f_{\sigma}(h)=\frac{1-4{\sigma}^2}{(1-2{\sigma}h)^2},\qquad x=-\frac15,\qquad y=\frac14.
\] Let \(H_{\sigma}(h)=\sum_{j=0}^4h_j({\sigma})h^j\) be the quartic matching
\(f_{\sigma},f_{\sigma}'\) at \(x,y\), and matching \(f_{\sigma}\) at 1. The positive fifth
derivative of \(f_{\sigma}\) and the Hermite remainder give \(H_{\sigma}\ge f_{\sigma}\) on
\([-1,1]\). For a complete rational specification set \[
 {D_{\mathrm{height}}}=(2-{\sigma})^3(1-2{\sigma})(5+2{\sigma})^3.
\] Then \[
\begin{aligned}
{D_{\mathrm{height}}} h_0={}&-(1+2{\sigma})(336{\sigma}^6-488{\sigma}^5+2642{\sigma}^4-1199{\sigma}^3\\
 &\hspace{19mm}-4630{\sigma}^2+4300{\sigma}-1000),\\
{D_{\mathrm{height}}} h_1={}&4{\sigma}(1+2{\sigma})(80{\sigma}^5-2522{\sigma}^4+1199{\sigma}^3\\
 &\hspace{19mm}+4630{\sigma}^2-4300{\sigma}+1000),\\
{D_{\mathrm{height}}} h_2={}&12{\sigma}^2(1+2{\sigma})(294{\sigma}^4-441{\sigma}^3+4630{\sigma}^2-4300{\sigma}+1000),\\
{D_{\mathrm{height}}} h_3={}&160{\sigma}^3(1+2{\sigma})(-42{\sigma}^3+923{\sigma}^2-860{\sigma}+200),\\
{D_{\mathrm{height}}} h_4={}&1600{\sigma}^4(1+2{\sigma})(2{\sigma}^2-83{\sigma}+50).
\end{aligned}
\] Equivalently, \[
 H_{\sigma}(h)=f_{\sigma}(x)+f_{\sigma}'(x)(h-x)
       +(h-x)^2[u_2+u_3(h-y)+u_4(h-y)^2],
\] with \[
 u_2=\frac{1200{\sigma}^2(5-{\sigma})(1-2{\sigma})(1+2{\sigma})}{(2-{\sigma})^2(5+2{\sigma})^3},\qquad
 u_3=\frac{1600{\sigma}^3(20-{\sigma})(1-2{\sigma})(1+2{\sigma})}{(2-{\sigma})^3(5+2{\sigma})^3},
 \qquad u_4=h_4.
\] The signs \(h_2,h_3,h_4>0\) on \(0<{\sigma}\le3/8\) can be checked directly.
For the parenthetical factors \(p_2,p_3,p_4\) in their numerators, \[
 p_2'\le-48991/64<0,\qquad p_2(3/8)=43319/2048>0,
\] \[
 p_3'\le-671/4<0,\qquad p_3(3/8)=1301/256>0,
 \qquad p_4\ge151/8>0.
\]

\subsubsection{A scalar lower bound}\label{bs:scalar-floor}

\begin{lemma}[A uniform Hermite floor]\label{bs:lemma}
For $a\le {\sigma}\le b$ and every real $h$, one has
\[
 H_{\sigma}(h)\ge23/50.
\]
\end{lemma}
\begin{proof}
\emph{Convexity and the right half-line.}
\label{bs:proof-convexity-and-the-right-half-line}

The preceding formulas give \(h_2,h_3,h_4>0\). Direct algebra gives \[
8h_2h_4-3h_3^2=
-\frac{230400{\sigma}^6(1+2{\sigma})^2K({\sigma})}{{D_{\mathrm{height}}}({\sigma})^2},
\] where \[
K({\sigma})=196{\sigma}^6-8988{\sigma}^5+267681{\sigma}^4-258160{\sigma}^3-24000{\sigma}^2+84000{\sigma}-20000.
\] Throughout \([a,b]\), termwise bounds give \[
\begin{aligned}
K({\sigma})\le{}&196b^6+267681b^4+84000b-20000\\
&-8988a^5-258160a^3-24000a^2\\
={}&-\frac{101905810790123033552}{2381242332078369737}<0.
\end{aligned}
\] Thus the discriminant of the quadratic \(H_{\sigma}''(h)\) is negative, and
its leading coefficient is positive. Consequently \(H_{\sigma}\) is strictly
convex on the real line.

Moreover \[
H_{\sigma}'(-1/2)=-\frac{4{\sigma}(1+2{\sigma})P_5({\sigma})}{{D_{\mathrm{height}}}({\sigma})},
\] where \[
P_5({\sigma})=2462{\sigma}^5-43091{\sigma}^4+48491{\sigma}^3-23530{\sigma}^2+7300{\sigma}-1000.
\] The exact shifted derivative is \[
P_5'(a+z)=12310z^4-\frac{2733228}{17}z^3+
\frac{8061201}{289}z^2-\frac{32369868}{4913}z+
\frac{173166660}{83521}.
\] Hence for \(0\le z\le 1/25\), \[
P_5'(a+z)>2000-6600z-161000z^3>0.
\] But \[
P_5(b)=-86826796000/15386239549<0.
\] It follows that \(H_{\sigma}'(-1/2)>0\). Strict convexity therefore implies
\(H_{\sigma}(h)\ge H_{\sigma}(-1/2)\) for \(h\ge -1/2\).

\par\smallskip\noindent\emph{The left half-line.}\label{bs:proof-the-left-half-line}

Write \(x=-1/5,y=1/4\). The Newton quadratic is
\(q_{\sigma}(v)=u_2+u_3v+u_4v^2\). For \(1/6\le {\sigma}\le 3/10\), \[
u_3/u_4\le3/2.
\] Indeed this is equivalent to \(7{\sigma}^3-(417/2){\sigma}^2+156{\sigma}-20\ge0\), whose
derivative is positive on that interval and whose value at \(1/6\) is
\(13/54\). For \(h\le -1/2\), one has \(v=h-y\le -3/4\), so \[
q_{\sigma}(v)\ge q_{\sigma}(-3/4)=A({\sigma})
=\frac{300{\sigma}^2(1+2{\sigma})p_{\mathrm{bottom}}({\sigma})}{{D_{\mathrm{height}}}({\sigma})},
\] where \[
p_{\mathrm{bottom}}({\sigma})=38{\sigma}^4-713{\sigma}^3+750{\sigma}^2-268{\sigma}+40>0.
\] For positivity it suffices to note \[
p_{\mathrm{bottom}}({\sigma})\ge(5361/10){\sigma}^2-268{\sigma}+40>0
\quad({\sigma}\le3/10),
\] because the quadratic has discriminant \(-13952\). Completing the
square in the Newton expression gives \[
H_{\sigma}(h)\ge f_{\sigma}(x)-\frac{f_{\sigma}'(x)^2}{4A({\sigma})}
=\frac{23}{50}+\frac{N({\sigma})}{150(5+2{\sigma})^3p_{\mathrm{bottom}}({\sigma})},
\] where \[
\begin{aligned}
N({\sigma})={}&-1660976{\sigma}^7+22276256{\sigma}^6+24404520{\sigma}^5\\
&-38875264{\sigma}^4+7412065{\sigma}^3+3504450{\sigma}^2-1452500{\sigma}+155000.
\end{aligned}
\] Its shifted derivative is \[
\begin{aligned}
N'(a+z)={}&-11626832z^6+\frac{1993134144}{17}z^5
+\frac{77917653960}{289}z^4\\
&+\frac{148921967232}{4913}z^3
-\frac{2515589043741}{83521}z^2
+\frac{2941403868}{1419857}z
-\frac{3125949300756}{24137569}.
\end{aligned}
\] Thus \[
N'(a+z)<-129000+2100z+31000000z^3+270000000z^4+118000000z^5<0
\] for \(0\le z\le 1/25\). Therefore \[
N({\sigma})\ge N(b)=\frac{33063217703060000}{182803912081669}>0.
\] The floor holds for \(h\le -1/2\), and in particular at \(h=-1/2\);
convexity and the positive derivative at that point prove it on the
other half-line.
\end{proof}

\begin{proof}[Proof of Theorem~\ref{bs:theorem}]
Theorem~\ref{thm:small} proves the assertion and its equality statement
for $0\le s\le1/4$. It remains to consider $1/4\le s\le3/10$.
We use the Hermite polynomial, coefficient signs and lower bound
established above, in particular Lemma~\ref{bs:lemma}.
For $R\le9/20$ we estimate the height distribution of $F_{\weights}$
(the sum of its three Fourier modes) directly.
For $R\ge9/20$ we prove the stated product bound, which implies the
bound for $J$ by Proposition~\ref{prop:global-product}.

\subsubsection{\texorpdfstring{Three-mode heights for
\(R\le 9/20\)}{Three-mode heights for R <= 9/20}}\label{bs:three-mode-heights}

Recall the moment envelope from Section~\ref{sec:small-centres},
\[
 B_{\mathrm{height}}(R,\sigma)
 =h_0+(h_2/2+h_3/4)R+(5h_4/8)R^2.
\] The analytic coefficient signs imply
$B_{\mathrm{height}}(R,\sigma)\le B_{\mathrm{height}}(9/20,\sigma)$
for $R\le9/20$. We show
\[
 B_{\mathrm{height}}(9/20,\sigma)<1999/2000\qquad(a\le\sigma\le b).
\]
It is enough to prove convexity and check its two endpoints. Direct
differentiation gives \[
\partial_{\sigma}^2 B_{\mathrm{height}}(9/20,{\sigma})=\frac{M_{\mathrm{bot}}({\sigma})}{5(2-{\sigma})^5(1-2{\sigma})^3(5+2{\sigma})^5},
\] with \[
\begin{aligned}
M_{\mathrm{bot}}({\sigma})={}&-303872{\sigma}^{10}+6039328{\sigma}^9-33107040{\sigma}^8+153016768{\sigma}^7\\
&-310617512{\sigma}^6+490649358{\sigma}^5-483961900{\sigma}^4\\
&+272558000{\sigma}^3-92580000{\sigma}^2+19250000{\sigma}-1300000.
\end{aligned}
\] For \(0\le z\le 1/25\), expansion about \(a=4/17\) gives \[
\begin{aligned}
 M_{\mathrm{bot}}(a+z)>{}&470000-9200000z^2-102000000z^4-104000000z^6\\
 &-22000000z^8-304000z^{10}>450000.
\end{aligned}
\] To make this short domination directly checkable, the nonpositive
shifted coefficients are \[
\begin{array}{c|r}
j &[z^j]M_{\mathrm{bot}}(a+z)\\\hline
2&-63798078090348000/6975757441\\
4&-2449469858978196/24137569\\
6&-8644311123624/83521\\
8&-6090653664/289\\
10&-303872
\end{array}
\] and all the other positive-degree coefficients are positive. The
constant term is
\[
 \frac{963573418694034720}{2015993900449}>470000.
\]
The complete shifted coefficient list is recorded in Section~\ref{sec:explicit-small-arithmetic}.

At the endpoints, \[
B_{\mathrm{height}}(9/20,a)=533877463/542940975<1999/2000,
\] \[
B_{\mathrm{height}}(9/20,b)=9007541157407/9012510107447<1999/2000.
\] Convexity proves the bound. Write \(w_*=1/\sqrt3-1/\pi\), so \({c_R}=2w_*/S(R)\), where
\(S(R)=\sqrt{(1-R)(1+3R)}\).
The bound \(w_*<25905/100000\) from Section~\ref{sec:analytic-highcaps} and
\[
 S(9/20)^2=\frac{517}{400}>
 \left(\frac{113}{100}\right)^2
\]
give
\[
 {c_R}\le {c_{9/20}}<\frac{5181}{11300}<\frac{23}{50}
 \qquad(R\le9/20).
\] Lemma~\ref{bs:lemma} therefore
yields \[
(f_{\sigma}(q)-L/A(\weights))_+\le H_{\sigma}(q)-{c_R}.
\] Using the three-mode moment bounds \eqref{eq:mom} and the nonnegative
right side gives \[
J\le\Abar(R)B_{\mathrm{height}}(R,\sigma)-L
<\frac{1999}{2000}\Aeq-L<4\pi.
\]

\subsubsection{\texorpdfstring{Product moments for
\(R\ge 9/20\)}{Product moments for R >= 9/20}}\label{bs:product-moments-for-rge-920}

Let \(\rho=\max(R,1-R)\). For the symmetric two-cosine product height
\(Q\), \[
\mathbb E Q^2\le\rho/2,
\qquad
\mathbb E Q^4\le M_4(\rho),
\] where \[
M_4(\rho)=
\begin{cases}
\displaystyle\frac{9\rho^2(1-\rho)^2}{8(6\rho(1-\rho)-1)},&1/2\le\rho\le2/3,\\
\displaystyle\frac{3\rho^2}{8},&\rho\ge2/3.
\end{cases}
\] This is the explicit maximum of a quadratic in the squared direction
amplitude. It increases with \(\rho\). For the first branch, put
\(w=\rho(1-\rho)\); the function \(w^2/(6w-1)\) decreases with
\(w<1/3\), while \(w\) decreases with \(\rho>1/2\). The second branch is
immediate.

Define \[
E_\rho({\sigma})=h_0+\frac\rho2h_2+M_4(\rho)h_4.
\] All odd product moments vanish, so \(\mathbb E H_{\sigma}(Q)\le E_\rho({\sigma})\).
The expression \(E_\rho\) increases in \(\rho\) by positivity of
\(h_2,h_4\). It also increases in \({\sigma}\) on \([a,b]\) whenever
\(\rho\ge 23/40\).

\paragraph{Elementary proof of the latter
monotonicity}\label{bs:elementary-proof-of-the-latter-monotonicity}

It suffices to prove \(E_{23/40}'>0\), \(h_2'>0\), and \(h_4'>0\), since
\(M_4\) increases. Their denominators are positive multiples of \[
\Delta({\sigma})=(2-{\sigma})^4(1-2{\sigma})^2(5+2{\sigma})^4.
\] More exactly, \[
E_{23/40}'=N_E/(14920\Delta),\qquad h_2'=N_2/\Delta,
\qquad h_4'=N_4/\Delta,
\] where \[
\begin{aligned}
N_E={}&{\sigma}(2475394256{\sigma}^7-10082324644{\sigma}^6+64814903040{\sigma}^5\\
&-71988791157{\sigma}^4+24545362980{\sigma}^3-2153818200{\sigma}^2
+722128000{\sigma}-164120000),\\
N_2={}&24{\sigma}(2352{\sigma}^7-72908{\sigma}^6+92320{\sigma}^5-368759{\sigma}^4\\
&+326960{\sigma}^3-31400{\sigma}^2-44000{\sigma}+10000),\\
N_4={}&3200{\sigma}^3(336{\sigma}^5-644{\sigma}^4+8320{\sigma}^3-5517{\sigma}^2-1100{\sigma}+1000).
\end{aligned}
\] Direct expansion at \({\sigma}=a+z\) gives the simple lower bounds \[
N_E(a+z)>7000000-3100000000z^3-10400000000z^4-5500000000z^7>6700000,
\] \[
\begin{aligned}
N_2(a+z)>{}&6100-102000z-59000z^2-1600000z^4\\
&-7800000z^5-600000z^6-1700000z^7>1900,
\end{aligned}
\] \[
N_4(a+z)>22000-13000z^2-3200000z^3-2900000z^4-37000z^7>21700,
\] valid throughout \(0\le z\le 1/25\). All omitted monomials have
positive coefficients. The complete shifted coefficients are displayed in Section~\ref{sec:explicit-small-arithmetic}; the stated estimates follow by termwise comparison.

\paragraph{Five rational weight
comparisons}\label{bs:five-rational-weight-comparisons}

The following endpoint bounds follow by substituting \({\sigma}=b\) in the
explicit \(h_j\) formulas:

\[
\begin{array}{ccl}
\rho&E_\rho(b)&\text{strict upper bound}\\\hline
23/40&574879035704129/572198514055784&201/200\\
13/20&1370163650087557/1315826475687262&521/500\\
7/10&9728489225507/9012510107447&27/25\\
71/100&9801176971352/9012510107447&136/125\\
1936/2561&393669821987973023/349766504759910623&563/500
\end{array}
\]

Write \(S(R)=\sqrt{(1-R)(1+3R)}\), so \(\Abar(R)=2\pi^2 S(R)\) and
\(S\) decreases for \(R\ge 1/3\). The estimates below use the scalar
floor \(23/50\) from Lemma~\ref{bs:lemma}. In the first row we lower the cutoff further to
\(9/20\); both floors are valid because \(H\ge 23/50\).
In each row, $\rho_+$ is the maximum of $\max(R,1-R)$ on the listed
interval, $S_+$ is an upper bound for $S(R)$ there, and $c_*$ is the
lower cutoff used in that row. The last column bounds
$S(R)(E_{\max(R,1-R)}(\sigma)-c_*)$ uniformly.

\[
\begin{array}{ccccc}
R\text{ interval}&\rho_+&S_+&c_*&S(E_\rho-c_*)\text{ upper bound}\\\hline
{[9/20,23/40]}&23/40&57/50&9/20&6327/10000\\
{[23/40,13/20]}&13/20&27/25&23/50&7857/12500\\
{[13/20,7/10]}&7/10&51/50&23/50&1581/2500\\
{[7/10,71/100]}&71/100&193/200&23/50&30301/50000\\
{[71/100,R_*]}&R_*&953/1000&23/50&317349/500000
\end{array}
\]

Every entry in the final column is strictly smaller than \(7/11\). The
area bounds follow by squaring at the left endpoint, then monotonicity.
The cutoff is large enough in each row: using \(\sqrt3<26/15\) and
\(\pi>157/50\), \[
{c_R}>\frac{100}{57}\left(\frac{15}{26}-\frac{50}{157}\right)>\frac9{20}
\quad(R\ge9/20),
\] where the final gap is \(7967/2326740\). For \(R\ge 23/40\), the
sharper area bound gives \[
{c_R}>\frac{50}{27}\left(\frac{15}{26}-\frac{50}{157}\right)>\frac{23}{50},
\] with final gap \(51289/2755350\). Therefore
\((f_{\sigma}(Q)-{c_R})_+\le H_{\sigma}(Q)-c_*\) in each row, and \[
\mathcal P(R,s,\alpha)
\le 2\pi^2 S(R)(E_\rho({\sigma})-c_*)
<\frac{14\pi^2}{11}<4\pi.
\] The last inequality is \(\pi<22/7\). This proves Theorem~\ref{bs:theorem}.
\end{proof}

\subsection{\texorpdfstring{Product estimates for \(3/10\le s\le2/5\)}{Product estimates for 3/10 <= s <=2/5}}\label{sec:central-strip}

We establish the product bound for \(3/10\le s\le2/5\) on the
central weight interval. A convex-hull comparison reduces the proof
to the single weight \(R=7/20\).

\begin{theorem}[Central product estimate]\label{thm:central-product-strip}
For every
\[
 \frac7{20}\le R\le\frac{13}{20},\qquad
 \frac3{10}\le s\le\frac25,\qquad 0\le\alpha\le1,
\]
one has \(\mathcal P(R,s,\alpha)<4\pi\).
\end{theorem}

\subsubsection{Weight reduction}\label{cp:statement-and-weight-reduction}

\begin{proposition}[Reduction to one weight]\label{prop:central-weight-reduction}
For every $0\le s<1$ and $7/20\le R\le13/20$,
\[
 \sup_{0\le\alpha\le1}\mathcal P(R,s,\alpha)
 \le\sup_{0\le\alpha\le1}\mathcal P(7/20,s,\alpha).
\]
\end{proposition}
\begin{proof}
For a fixed scalar function $f$ convex on $[-1,1]$, the function
$\Phi(a,b)=\Av f(a\cos U+b\cos V)$ is convex on $|a|+|b|\le1$,
separately even and symmetric in its arguments. For $0<\rho<1$, let
$\mathscr E_\rho$ be the filled ellipse $a^2/\rho+b^2/(1-\rho)\le1$.
Cauchy--Schwarz places each such ellipse, and hence the convex hulls
below, in $|a|+|b|\le1$. For $1/2\le\rho<1$, the support function in direction $(\xi,\eta)\in\mathbb R^2$ of
$\operatorname{conv}(\mathscr E_\rho\cup \mathscr E_{1-\rho})$ is
\[
 \sqrt{\rho\max(\xi^2,\eta^2)+(1-\rho)\min(\xi^2,\eta^2)}.
\]
The convex hulls are therefore nested as $\rho$ increases. By convexity
and symmetry, the supremum of $\Phi$ on each hull equals its supremum
on the positive boundary arc of $\mathscr E_\rho$. Thus this supremum is
nondecreasing in $\max(R,1-R)$, which is at most $13/20$ on the stated
interval and equals $13/20$ when $R=7/20$.

The area envelope $\Abar(R)$ is nonincreasing and the cutoff
$c_R=L/\Abar(R)$ is nondecreasing for $R\ge1/3$. In the nonnegative
integrand defining $\mathcal P$, first replace $\Abar(R),c_R$ by
$\Abar(7/20),c_{7/20}$; then
apply the preceding convex-hull comparison to the fixed convex
function $f(q)=(\chi_s(q)^2-c_{7/20})_+$. This proves the assertion.
\end{proof}

\subsubsection{Fixed constants and
interpolation}\label{cp:fixed-constants-and-interpolation}

Set \[
 R_0=7/20,\quad c=4487/10000,\quad d=110299/200000,
 \quad {\sigma}=s/(1+s^2),\quad w=4{\sigma}^2,
\] \[
 f_{\sigma}(q)=\frac{1-w}{(1-2{\sigma}q)^2},\qquad
 Q=\sqrt{7x/20}\cos U+\sqrt{13(1-x)/20}\cos V,
 \quad x=\alpha^2\in[0,1].
\] At this weight, \[
 A_0=\pi^2\sqrt{533}/10,\qquad
 c_0=\frac{40}{\sqrt{533}}\left(\frac1{\sqrt3}-\frac1\pi\right).
\] The following rational bounds give
\[
 c<c_0,\qquad d<4\pi/A_0.
\] Explicitly, \(\sqrt{533}<23087/1000\), \(\sqrt3<97/56\), and
\(333/106<\pi<355/113\) give \[
 c_0>334640000/745733187>c,
 \qquad4\pi/A_0>4520000/8195885>d.
\] Since \(Q\) is symmetric, \[
 {\mathbb E}(f_{\sigma}(Q)-c)_+={\mathbb E} g_{\sigma}(Q^2),
\] where \[
 g_{\sigma}(z)=\max\{0,g_B(z),g_A(z)\},\quad
 g_B(z)=\tfrac12(f_{\sigma}(\sqrt z)-c),\quad
 g_A(z)=\frac{(1-w)(1+wz)}{(1-wz)^2}-c.
\] Let \({\widehat H}_{\sigma}(z)=\sum_{j=0}^4{\widehat h}_j({\sigma})z^j\) be the quartic Hermite
polynomial matching \(g_B,g_B'\) at \(a=1/9,b=9/16\), and matching
\(g_B\) at 1. For an entirely rational definition, put \[
 u_2=\frac{g_B(b)-g_B(a)-g_B'(a)(b-a)}{(b-a)^2},
\] \[
 u_3=\frac{g_B'(b)-g_B'(a)-2u_2(b-a)}{(b-a)^2},
\] \[
 u_4=\frac{g_B(1)-g_B(a)-g_B'(a)(1-a)-u_2(1-a)^2-u_3(1-a)^2(1-b)}{(1-a)^2(1-b)^2},
\] \[
 {\widehat H}_{\sigma}(z)=g_B(a)+g_B'(a)(z-a)+u_2(z-a)^2
       +u_3(z-a)^2(z-b)+u_4(z-a)^2(z-b)^2.
\] Here \(g_B'(q^2)={\sigma}(1-4{\sigma}^2)/[q(1-2{\sigma}q)^3]\) at positive rational \(q\).

For every \(0<{\sigma}<1/2\), \({\widehat H}_{\sigma}\ge g_B\) on \([0,1]\). Indeed,
\((1+2{\sigma}\sqrt z)^{-2}\) is completely monotone: repeated application of
\(-(2\sqrt z)^{-1}d/d\sqrt z\) to products of inverse powers of
\(\sqrt z\) and \(1+2{\sigma}\sqrt z\) preserves nonnegative coefficients.
Meanwhile, \[
 (1-2{\sigma}\sqrt z)^{-2}+(1+2{\sigma}\sqrt z)^{-2}
 =\frac{2(1+4{\sigma}^2z)}{(1-4{\sigma}^2z)^2}
\] has positive derivatives of all orders. Thus \(g_B^{(5)}>0\). The
Hermite remainder has the sign of \((z-a)^2(z-b)^2(z-1)\), proving the
assertion, with \(z=0\) by continuity.

Define \[
 \delta({\sigma})=1-4{\sigma}^2-c-{\widehat h}_0({\sigma}),\qquad n_k({\sigma})={\widehat h}_1({\sigma})-k\delta({\sigma})-3w(1-w).
\] The majorant used below is \[
 P_{\sigma}(z)=
 \begin{cases}
 {\widehat H}_{\sigma}(z)+\delta({\sigma})(1-z)^7,&11/40\le {\sigma}<23/80,\\
 {\widehat H}_{\sigma}(z)+\delta({\sigma})(1-z)^8,&23/80\le {\sigma}<1/3,\\
 {\widehat H}_{\sigma}(z)+(1-z)^8/300,&1/3\le {\sigma}\le10/29.
 \end{cases}
\] For convenience each of the following estimates is proved on the
closed interval associated with the corresponding polynomial; at a
shared endpoint the above definition selects the polynomial from the
interval on the right. Since \({\sigma}(3/10)=30/109>11/40\) and
\({\sigma}(2/5)=10/29\), this covers the required range.

\subsubsection{Rational identities and
monotonicities}\label{cp:rational-identities-and-monotonicities}

Put \[
 {D_{\mathrm{prod}}}({\sigma})=6459180000(3-2{\sigma})^3(1-2{\sigma})(2-3{\sigma})^3>0,
 \quad N_j={D_{\mathrm{prod}}}{\widehat h}_j,\quad N_\delta={D_{\mathrm{prod}}}\delta.
\] For direct checking, the numerators are \[
\begin{aligned}
N_0={}&3(10698381172752{\sigma}^7-33010875554264{\sigma}^6+40021442287624{\sigma}^5\\
 &\hspace{21mm}-22448422581950{\sigma}^4+3787025203735{\sigma}^3\\
 &\hspace{21mm}+1860784191522{\sigma}^2-964284945804{\sigma}+128194053624),\\
N_1={}&-2500{\sigma}(2{\sigma}+1)(47356985312{\sigma}^5-130682010590{\sigma}^4+142897698425{\sigma}^3\\
 &\hspace{16mm}-78064690338{\sigma}^2+21400278732{\sigma}-2357020728),\\
N_2={}&1250{\sigma}(2{\sigma}+1)(275157623104{\sigma}^5-692049165454{\sigma}^4+673920615625{\sigma}^3\\
 &\hspace{16mm}-322739797698{\sigma}^2+76882344588{\sigma}-7322128056),\\
N_3={}&-360000{\sigma}(2{\sigma}+1)(1164374848{\sigma}^5-2800985998{\sigma}^4+2587536073{\sigma}^3\\
 &\hspace{16mm}-1181143746{\sigma}^2+273441420{\sigma}-26042040),\\
N_4={}&25920000{\sigma}(2{\sigma}+1)(6889792{\sigma}^5-16178254{\sigma}^4+14545129{\sigma}^3\\
 &\hspace{16mm}-6509826{\sigma}^2+1499148{\sigma}-142776),\\
N_\delta={}&(1-4{\sigma}^2-c){D_{\mathrm{prod}}}-N_0.
\end{aligned}
\]

Write the exact moments as \[
 M_{2j}(x)=\sum_{i=0}^j
 \frac{\binom{2j}{2i}\binom{2i}{i}\binom{2j-2i}{j-i}}{4^j}
 (7x/20)^i(13(1-x)/20)^{j-i}.
\] The first five are \[
\begin{aligned}
M_0={}&1,\\
M_2={}&13/40-3x/20,\\
M_4={}&507/3200+39x/1600-219x^2/1600,\\
M_6={}&2197/25600+507x/3200-897x^2/2560+153x^3/1280,\\
M_8={}&199927/4096000+46137x/204800-188097x^2/409600\\
 &+34671x^3/204800+8169x^4/409600.
\end{aligned}
\] Define \[
 C_k(x)=\sum_{j=0}^k(-1)^j\binom kjM_{2j}(x)
       ={\mathbb E}(1-Q^2)^k=\sum_j c_{k,j}x^j,
\] \[
 B_{\mathrm{prod}}({\sigma},x)=\sum_{j=0}^4{\widehat h}_jM_{2j}=\sum_{j=0}^4b_j({\sigma})x^j,
 \quad F_k({\sigma},x)=B_{\mathrm{prod}}({\sigma},x)+\delta({\sigma})C_k(x)=\sum_{j=0}^ka_{k,j}({\sigma})x^j.
\] All coefficients are thus explicitly defined rational functions.

The identities in Section~\ref{cp:complete-positive-expansions} prove these monotonicities:

\begin{itemize}
\setlength{\itemsep}{0pt}\setlength{\parskip}{0pt}
\item
  On \([11/40,7/20]\), the functions
  \({\widehat h}_2,-{\widehat h}_3,{\widehat h}_4,-\delta,n_7,n_8,a_{8,1},-b_2,a_{8,3}\) are strictly
  increasing.
\item
  On \([8/25,7/20]\), \(a_{8,0}\) is strictly decreasing.
\item
  On \([1/3,7/20]\), \(b_0\) is strictly decreasing.
\item
  The two direct comparisons \[
  a_{7,0}<551/1000\quad(11/40\le {\sigma}\le23/80),
  \qquad
  a_{8,0}<551/1000\quad(23/80\le {\sigma}\le8/25)
  \] hold.
\end{itemize}

Some immediate consequences used repeatedly are worth recording. Since
\(c_{8,2}>0\), \(a_{8,2}=b_2+\delta c_{8,2}\) decreases. Since
\(c_{8,4}<0\) and \(b_4=(8169/409600){\widehat h}_4\), \(a_{8,4}\) increases. Since
\(c_{8,1},c_{8,3}>0\), monotonicity of \(a_{8,1},a_{8,3}\), together
with \(\delta'<0\), implies that \(b_1,b_3\) increase. Also \(b_4\)
increases. Finally, \[
 \delta(1/3)=\frac{428911057}{135642780000}>0,
 \qquad
 \frac1{300}-\delta(1/3)=\frac{23231543}{135642780000}>0.
\]

\subsubsection{Scalar majorization by a cubic Taylor
estimate}\label{cp:scalar-majorization-by-a-cubic-taylor-estimate}

We prove \(P_{\sigma}\ge g_{\sigma}\) on \([0,1]\). Throughout the required range,
\(g_B(z)\ge g_B(0)>0\), since \(1-4(10/29)^2=441/841>c\). Thus the
nonnegative branch is already covered by \({\widehat H}_{\sigma}\ge g_B\). The correction
is positive in all three ranges, by the monotonicity of $-\delta$ and
the value $\delta(1/3)>0$ established in
Section~\ref{cp:rational-identities-and-monotonicities}. It remains to
majorize \(g_A\).

For a row of the following table, let \([a,b]\) be its parameter
interval, let \(k\) be its correction exponent, and put \(Z=q^2\). The
four columns \(n,A,B,C\) are rational lower bounds. The endpoint values
listed after the table are exactly \(n+AZ+BZ^2+CZ^3\).

\[
\begin{array}{ccccrrrr}
a&b&k&q&n&A&B&C\\\hline
11/40&7/25&7&9/20&29/1000&23/100&-69/50&-33/100\\
7/25&23/80&7&3/7&65/1000&1/20&-6/5&-3/10\\
23/80&3/10&8&2/5&59/1000&1/20&-221/100&-1/4\\
3/10&8/25&8&1/3&4/25&-61/100&-8/5&-7/50\\
8/25&1/3&8&1/4&8/25&-89/100&-1/2&0
\end{array}
\] The corresponding endpoint lower bounds, in the same order, are \[
\frac{103975247}{6400000000},\quad
\frac{749227}{23529800},\quad\frac{47}{5000},\quad
\frac{5269}{72900},\quad\frac{3359}{12800}.
\]

Here is a precise derivation and a way to check every table entry using
only rational operations. With \(w_b=4b^2\), put \[
 K_{b,Z}=\frac{(1-w_b)w_b^2(5+w_b Z)}{(1-w_b Z)^4}.
\] The lower bounds in each row satisfy \[
\begin{aligned}
 n&\le n_k(a),\\
 A&\le {\widehat h}_2(a)+\binom k2\delta(b)-K_{b,Z},\\
 B&\le {\widehat h}_3(b)-\binom k3\delta(a),\\
 C&\le\min\{0,{\widehat h}_4(a)\}.
\end{aligned}
\] The exact rational differences are displayed in Section~\ref{sec:explicit-central-arithmetic}. All their differences are positive except the deliberate
equality \(C=0\) in the last row. Also direct substitution gives \[
 f_a(-q)<c
\] for all five rows. Since \(f_{\sigma}(-q)\) decreases in \({\sigma}\) and in
\(q\ge0\), the identity \[
 g_A(z)-g_B(z)=\tfrac12(f_{\sigma}(-\sqrt z)-c)
\] shows that \(g_A\le g_B\) on \(z\ge Z\).

For \(0\le z\le Z\), \[
 g_A'(0)=3w(1-w),\qquad
 \tfrac12g_A''(z)=\frac{(1-w)w^2(5+wz)}{(1-wz)^4}\le K_{b,Z}.
\] The latter inequality uses that \(w^2(1-w)\) increases for
\(w\le2/3\), and each remaining positive factor increases with \(w,z\).
Taylor's theorem gives \[
 g_A(z)\le g_A(0)+3w(1-w)z+K_{b,Z}z^2.
\] The fourth derivative of \((1-z)^k\) is nonnegative on \([0,1]\),
hence \[
 (1-z)^k\ge1-kz+\binom k2z^2-\binom k3z^3.
\] Using \({\widehat h}_0+\delta=g_A(0)\), the monotonicities of Section~\ref{cp:rational-identities-and-monotonicities} give \[
 P_{\sigma}(z)-g_A(z)\ge z(n+Az+Bz^2+Cz^3).
\] The cubic on the right is concave on \([0,Z]\), because \(B,C\le0\).
Its values at 0 and at \(Z\) are positive, as the table shows. This
proves the assertion in the first two majorant ranges.

For \(1/3\le {\sigma}\le10/29\), use \(P_{\sigma}={\widehat H}_{\sigma}+(1-z)^8/300\). We have
\(P_{\sigma}(0)>g_A(0)\) because \(\delta\) decreases and
\(\delta(1/3)<1/300\). Also \(f_{1/3}(-1/5)=125/289<c\), so only
\(z\le1/25\) needs attention. Direct endpoint evaluation and the
monotonicities give \[
 n_8({\sigma})>2/5,\quad \delta({\sigma})>-1/80,\quad
 3w(1-w)\ge20/27,\quad {\widehat h}_2({\sigma})>-4/25,\quad {\widehat h}_3({\sigma}),{\widehat h}_4({\sigma})>0.
\] Consequently \({\widehat h}_1=n_8+8\delta+3w(1-w)>1\), and \[
 P_{\sigma}'(z)\ge1-\frac8{25}z-\frac2{75}
 \ge1-\frac8{625}-\frac2{75}>\frac9{10}.
\] On the other hand, \[
 g_A'(z)=\frac{w(1-w)(3+wz)}{(1-wz)^3}
 \le\frac{3+z}{4(1-z)^3}
 \le\frac{3+1/25}{4(1-1/25)^3}<\frac9{10}.
\] Thus \(P_{\sigma}\ge g_A\) throughout this last interval as well.

\subsubsection{The moment comparison in four
intervals}\label{cp:the-moment-comparison-in-four-intervals}

The comparison needed is \({\mathbb E} P_{\sigma}(Q^2)<d\). All estimates below hold for
every \(0\le x\le1\).

\paragraph{\texorpdfstring{\(11/40\le {\sigma}\le23/80\)}{11/40 <= sigma <= 23/80}}\label{cp:le-tle2380}

We use \(F_7\). The coefficient signs needed are \[
 a_{7,1}<0,\quad a_{7,3},a_{7,5},a_{7,6}>0,
 \quad a_{7,4},a_{7,7}<0,
\] \[
 a_{7,2}+a_{7,3}+a_{7,5}+a_{7,6}<0.
\] Here are endpoint formulas checking the only nonautomatic
comparisons. Writing \(a=11/40,b=23/80\), one has \[
 a_{7,1}({\sigma})
 =a_{8,1}({\sigma})+\delta({\sigma})(c_{7,1}-c_{8,1})
 \le a_{8,1}(b)+\delta(a)(c_{7,1}-c_{8,1})<-1/20.
\] The positive difference \(c_{7,1}-c_{8,1}\) is given explicitly by
the moment definition. Furthermore, \[
\begin{aligned}
 a_{7,2}+a_{7,3}+a_{7,5}+a_{7,6}
 ={}&-\frac{219}{1600}{\widehat h}_2-\frac{591}{2560}{\widehat h}_3-\frac{118755}{409600}{\widehat h}_4\\
 &+\delta(c_{7,2}+c_{7,3}+c_{7,5}+c_{7,6})\\
 \le{}&-\frac{219}{1600}{\widehat h}_2(a)-\frac{591}{2560}{\widehat h}_3(b)
 -\frac{118755}{409600}{\widehat h}_4(a)\\
 &+\delta(a)(c_{7,2}+c_{7,3}+c_{7,5}+c_{7,6})<-1/40.
\end{aligned}
\] Also \[
 a_{7,3}\ge\frac{153}{1280}{\widehat h}_3(b)+\frac{34671}{204800}{\widehat h}_4(a)>0.
\] The signs of \(a_{7,5},a_{7,6},a_{7,7}\) follow from those of the
respective coefficients of \(C_7\), because \(\delta>0\). Finally
\({\widehat h}_4(b)<0\) and \(c_{7,4}<0\) give \(a_{7,4}<0\). Thus \[
 F_7({\sigma},x)\le a_{7,0}+a_{7,1}x
 +(a_{7,2}+a_{7,3}+a_{7,5}+a_{7,6})x^2
 \le a_{7,0}<551/1000<d.
\]

\paragraph{\texorpdfstring{\(23/80\le {\sigma}\le8/25\)}{23/80 <= sigma <= 8/25}}\label{cp:le-tle825}

We use \(F_8\). The coefficient \(a_{8,1}\) increases and \[
 a_{8,1}(8/25)
 =-\frac{16071393982926641546748379699}
 {4701647180412852633600000000000}<0.
\] The coefficients \(a_{8,3},a_{8,6}\) are positive; \(a_{8,4}<0\)
because it increases and \(a_{8,4}(8/25)<0\); and
\(a_{8,5},a_{8,7},a_{8,8}<0\). Since \(a_{8,2}\) decreases, \(a_{8,3}\)
increases, and \(a_{8,6}=\delta c_{8,6}\) decreases, \[
 a_{8,2}({\sigma})+a_{8,3}({\sigma})+a_{8,6}({\sigma})
 \le a_{8,2}(23/80)+a_{8,3}(8/25)+a_{8,6}(23/80)<0.
\] The exact positive differences are recorded in Section~\ref{sec:explicit-central-arithmetic}.
Hence \[
 F_8({\sigma},x)\le a_{8,0}+a_{8,1}x
 +(a_{8,2}+a_{8,3}+a_{8,6})x^2
 \le a_{8,0}<551/1000<d.
\]

\paragraph{\texorpdfstring{\(8/25\le {\sigma}\le1/3\)}{8/25 <= sigma <= 1/3}}\label{cp:le-tle13}

The monotonicities in Section~\ref{cp:rational-identities-and-monotonicities}
and endpoint evaluations give \[
 a_{8,0}<1099/2000,\quad a_{8,1}<1/30,\quad
 a_{8,2}<-171/500,\quad a_{8,3}<81/500,
\] \[
 a_{8,4}<1/125,\qquad a_{8,6}<13/1000,
 \qquad a_{8,5},a_{8,7},a_{8,8}<0.
\] For \(a_{8,0},a_{8,2},a_{8,6}\), use their values at \(8/25\); for
\(a_{8,1},a_{8,3},a_{8,4}\), use their values at \(1/3\). The
coefficients of the positive powers above \(x^2\) can be bounded by
their displayed positive upper bounds times \(x^2\). Their sum with the
bound for \(a_{8,2}\) is \(-159/1000<-3/20\). Therefore \[
 F_8({\sigma},x)<\frac{1099}{2000}+\frac{x}{30}-\frac3{20}x^2
 \le\frac{1099}{2000}+\frac1{540}
 =\frac{29773}{54000}<d.
\] The final gap is exactly \(d-29773/54000=773/5400000\).

\paragraph{\texorpdfstring{\(1/3\le {\sigma}\le10/29\)}{1/3 <= sigma <= 10/29}}\label{cp:le-tle1029}

Now \({\mathbb E} P_{\sigma}(Q^2)=B_{\mathrm{prod}}({\sigma},x)+C_8(x)/300\). The monotone coefficients of \(B_{\mathrm{prod}}\), evaluated at the corresponding endpoints, give \[
 b_0<109/200,\quad b_1<77/1000,\quad b_2<-43/100,
 \quad 0<b_3<187/1000,\quad0<b_4<19/1000.
\] The Bernstein coefficients of the fixed polynomial \(C_8\) on
\([0,1]\) are all below \(1/2\). Their maximum is \[
 \frac{3948897100181}{8388608000000}<\frac12.
\] They can be obtained from the displayed moment formula by
\(\beta_i=\sum_{j=0}^i c_{8,j}\binom ij/\binom8j\). All nine exact
coefficients are printed in Section~\ref{cp:complete-positive-expansions}. Thus \(C_8/300<1/600\).

For \(0\le x\le1/4\), \[
 B_{\mathrm{prod}}({\sigma},x)<109/200+(77/1000)x-(19/50)x^2
 \le109/200+5929/1520000.
\] Indeed, \(-43/100+(187/1000)/4+(19/1000)/16<-19/50\). Adding
\(1/600\) gives a rational number strictly below \(d\). For
\(1/4\le x\le1/2\), \[
 B_{\mathrm{prod}}({\sigma},x)<109/200+(77/1000)x-(33/100)x^2
 \le109/200+77/4000-33/1600,
\] since this quadratic decreases for \(x\ge1/4\). For \(1/2\le x\le1\),
\[
 B_{\mathrm{prod}}({\sigma},x)<109/200+(77/1000)x-(11/50)x^2
 \le109/200+77/2000-11/200.
\] This quadratic decreases for \(x\ge1/2\). In both cases, adding
\(1/600\) again gives a bound strictly below \(d\).

\subsubsection{Conclusion}\label{cp:conclusion}

\begin{proof}[Proof of Theorem~\ref{thm:central-product-strip}]
Sections~\ref{cp:scalar-majorization-by-a-cubic-taylor-estimate} and~\ref{cp:the-moment-comparison-in-four-intervals} give \[
 {\mathbb E}(f_{\sigma}(Q)-c)_+={\mathbb E} g_{\sigma}(Q^2)\le {\mathbb E} P_{\sigma}(Q^2)<d.
\] Consequently \[
 \mathcal P(7/20,s,\alpha)
 \le A_0{\mathbb E}(f_{\sigma}(Q)-c)_+<A_0d<4\pi.
\] Proposition~\ref{prop:central-weight-reduction} proves Theorem~\ref{thm:central-product-strip}.
\end{proof}

\subsection{\texorpdfstring{An eighth-degree majorant for \(2/5\le s\le1/2\)}{An eighth-degree majorant for 2/5 <= s <=1/2}}\label{sec:octic-product-strip}

The product majorant of Proposition~\ref{prop:global-product} admits an
exact reduction in the weight parameter and a further algebraic bound
on a fixed strip. Section~\ref{octic:positive-expansions} supplies the finite
comparisons as thirteen positive polynomial expansions.

Recall that $c_R=L/\Abar(R)$ and
\[
 Q_{R,\alpha}=\sqrt R\,\alpha\cos U+
 \sqrt{1-R}\sqrt{1-\alpha^2}\cos V,
\]
where $U,V$ are independent uniform angles. Then
\[
 \mathcal P(R,s,\alpha)
 =\Abar(R)\Av\bigl(\chi_s(Q_{R,\alpha})^2-c_R\bigr)_+.
\]
Cauchy--Schwarz gives $|Q_{R,\alpha}|\le1$, so the scalar majorants
below are used on their stated domain. All additional notation in this
section is local.

\begin{samepage}
\begin{theorem}[An octic product bound]\label{thm:octic-product-strip}
For every
\[
 \frac7{20}\le R\le\frac{13}{20},\qquad
 \frac25\le s\le\frac12,\qquad0\le\alpha\le1,
\]
one has $\mathcal P(R,s,\alpha)<4\pi$.
\end{theorem}
\end{samepage}

\subsubsection{Rational constants and the symmetrized integrand}
By Proposition~\ref{prop:central-weight-reduction}, it suffices to consider
$R_0=7/20$. At this weight,
\[
 A_0=\frac{\pi^2\sqrt{533}}{10},\qquad
 c_0=\frac{40}{\sqrt{533}}\left(\frac1{\sqrt3}-\frac1\pi\right).
\]
We use
\begin{equation}\label{octic:constants}
 c=\frac{4487}{10000}<c_0,\qquad
 d=\frac{110299}{200000}<\frac{4\pi}{A_0}.
\end{equation}
Indeed $\sqrt{533}<23087/1000$, $\sqrt3<97/56$ and
$333/106<\pi<355/113$ imply
\[
 c_0>\frac{334640000}{745733187}>\frac{4487}{10000},
\]
\[
 \frac{4\pi}{A_0}>\frac{4520000}{8195885}>
 \frac{110299}{200000}.
\]
Set ${\sigma}=s/(1+s^2)\in[10/29,2/5]$ and
\[
 \widehat f_{\sigma}(q)=\frac{1-4{\sigma}^2}{(1-2{\sigma}q)^2}.
\]
For $Q=Q_{R_0,\alpha}$, symmetry gives
\[
 \Av(\widehat f_{\sigma}(Q)-c)_+=\Av g_{\sigma}(Q^2),
\]
where
\begin{equation}\label{octic:branches}
 \begin{aligned}
 g_{\sigma}(z)&=\max\{0,g_B(z),g_A(z)\},\qquad 0\le z\le1,\\
 g_B(z)&=\tfrac12\bigl(\widehat f_{\sigma}(\sqrt z)-c\bigr),\\
 g_A(z)&=\frac{(1-4{\sigma}^2)(1+4{\sigma}^2z)}{(1-4{\sigma}^2z)^2}-c.
 \end{aligned}
\end{equation}

\subsubsection{A structural Hermite majorant}
Put $a=1/9$, $b=9/16$. Let ${\widehat H}_{\sigma}(z)=\sum_{j=0}^4{\widehat h}_j({\sigma})z^j$ be the
quartic satisfying
\[
 {\widehat H}_{\sigma}(a)=g_B(a),\quad {\widehat H}_{\sigma}'(a)=g_B'(a),
\]
\[
 {\widehat H}_{\sigma}(b)=g_B(b),\quad {\widehat H}_{\sigma}'(b)=g_B'(b),\quad {\widehat H}_{\sigma}(1)=g_B(1).
\]
These conditions define the coefficients rationally. Explicitly, set
\begin{align*}
 u_2&=\frac{g_B(b)-g_B(a)-g_B'(a)(b-a)}{(b-a)^2},\\
 u_3&=\frac{g_B'(b)-g_B'(a)-2u_2(b-a)}{(b-a)^2},\\
 u_4&=\frac{g_B(1)-g_B(a)-g_B'(a)(1-a)-u_2(1-a)^2
             -u_3(1-a)^2(1-b)}{(1-a)^2(1-b)^2}.
\end{align*}
Then
\begin{equation}\label{octic:interpolant}
 \begin{aligned}
 {\widehat H}_{\sigma}(z)={}&g_B(a)+g_B'(a)(z-a)+u_2(z-a)^2\\
 &+u_3(z-a)^2(z-b)+u_4(z-a)^2(z-b)^2.
 \end{aligned}
\end{equation}

\begin{lemma}\label{octic:hermite-majorant}
For ${\sigma}\in[10/29,2/5]$ and $0\le z\le1$, one has
${\widehat H}_{\sigma}(z)\ge g_B(z)$.
\end{lemma}
\begin{proof}
For $\beta>0$, the function $(1+\beta\sqrt z)^{-2}$ is completely
monotone on $z>0$. Indeed, application of
$-(2\sqrt z)^{-1}d/d\sqrt z$ to a product of nonnegative inverse
powers of $\sqrt z$ and $1+\beta\sqrt z$ gives a sum of the same
kind with positive coefficients. Induction proves the assertion.
On the other hand,
\[
 (1-\beta\sqrt z)^{-2}+(1+\beta\sqrt z)^{-2}
 =\frac{2(1+\beta^2z)}{(1-\beta^2z)^2}
\]
has positive derivatives of every order on $0<z<\beta^{-2}$.
Every odd derivative of $(1-\beta\sqrt z)^{-2}$ is therefore
positive there. Taking $\beta=2{\sigma}$ gives $g_B^{(5)}>0$ on $(0,1]$.
The Hermite remainder is
\[
 g_B(z)-{\widehat H}_{\sigma}(z)=\frac{g_B^{(5)}(\xi)}{5!}
 (z-a)^2(z-b)^2(z-1)\le0\qquad(z>0).
\]
Continuity gives the endpoint $z=0$.
\end{proof}

The identities in Section~\ref{octic:positive-expansions} give
\begin{equation}\label{octic:coefficient-signs}
 {\widehat h}_0,{\widehat h}_1,{\widehat h}_2,{\widehat h}_4>0,\qquad 32{\widehat h}_2{\widehat h}_4-9{\widehat h}_3^2>0.
\end{equation}
Consequently
\[
 {\widehat H}_{\sigma}'(z)={\widehat h}_1+z(2{\widehat h}_2+3{\widehat h}_3z+4{\widehat h}_4z^2)\ge {\widehat h}_1>0,
\]
and ${\widehat H}_{\sigma}\ge {\widehat h}_0>0$. For $10/29\le {\sigma}\le3/8$, the identities
\eqref{octic:positive-identities} for $P_6,P_7$, defined in
\eqref{octic:polynomials}, give
\[
 {\widehat h}_0\ge1-4{\sigma}^2-c=g_A(0),\qquad {\widehat h}_1>4/5.
\]
For $0\le z\le1/64$,
\[
 \begin{aligned}
 g_A'(z)&=\frac{(1-4{\sigma}^2)4{\sigma}^2(3+4{\sigma}^2z)}{(1-4{\sigma}^2z)^3}\\
 &\le\frac{3+z}{4(1-z)^3}
 \le\frac{197632}{250047}<\frac45.
 \end{aligned}
\]
Thus ${\widehat H}_{\sigma}\ge g_A$ on this interval. For $z\ge1/64$, use the exact
identity $\widehat f_{10/29}(-1/8)=4/9<c$. Since $\widehat f_{\sigma}(q)$ decreases in ${\sigma}$
for $q\le0$, it follows that $\widehat f_{\sigma}(-\sqrt z)<c$ and hence
$g_A\le g_B$. Finally, for ${\sigma}\ge3/8$, $\widehat f_{\sigma}(0)\le7/16<c$, so
$g_A\le g_B$ everywhere. Together with Lemma~\ref{octic:hermite-majorant},
these observations prove
\begin{equation}\label{octic:full-majorant}
 {\widehat H}_{\sigma}(z)\ge g_{\sigma}(z)\qquad(0\le z\le1).
\end{equation}

\subsubsection{The moment comparison}
Write $x=\alpha^2$. The exact moments at $R_0=7/20$ are
\begin{equation}\label{octic:moments}
 \begin{aligned}
 M_{2j}(x)&=\Av Q^{2j}\\
 &=\sum_{i=0}^j
 \frac{\binom{2j}{2i}\binom{2i}{i}\binom{2j-2i}{j-i}}{4^j}
 (7x/20)^i(13(1-x)/20)^{j-i}.
 \end{aligned}
\end{equation}
Define the positive denominator
\begin{equation}\label{octic:denominator}
 {D_{\mathrm{prod}}}({\sigma})=6459180000(3-2{\sigma})^3(1-2{\sigma})(2-3{\sigma})^3,
\end{equation}
the polynomials $N_j={D_{\mathrm{prod}}}{\widehat h}_j$, and
\begin{equation}\label{octic:residual}
 \begin{aligned}
 E({\sigma},x)&={D_{\mathrm{prod}}}({\sigma})\left[d-\sum_{j=0}^4{\widehat h}_j({\sigma})M_{2j}(x)\right]\\
 &=\sum_{j=0}^4E_j({\sigma})x^j.
 \end{aligned}
\end{equation}
Put $K=E_2+E_3/4+E_4/16$. The identities in
Section~\ref{octic:positive-expansions} prove
\begin{equation}\label{octic:moment-signs}
 \begin{gathered}
 E_3,E_4<0,\qquad K>0,\qquad4E_0K-E_1^2>0,\\
 E_x({\sigma},1/4)>0,\qquad E_x({\sigma},1)>0.
 \end{gathered}
\end{equation}

\begin{proof}[Proof of Theorem~\ref{thm:octic-product-strip}]
For $0\le x\le1/4$, \eqref{octic:moment-signs} gives
\[
 E({\sigma},x)\ge E_0+E_1x+Kx^2>0
\]
by completion of the square. For $1/4\le x\le1$, the function
$E_x({\sigma},x)$ is concave, since $E_{xxx}=6E_3+24E_4x<0$.
Its positive endpoint values imply $E_x>0$ throughout this interval;
hence $E({\sigma},x)\ge E({\sigma},1/4)>0$. Thus, by
\eqref{octic:full-majorant} and \eqref{octic:residual},
\[
 \Av(\widehat f_{\sigma}(Q)-c)_+\le\Av {\widehat H}_{\sigma}(Q^2)<d.
\]
Since $c<c_0$ and $A_0d<4\pi$ by \eqref{octic:constants}, this proves
$\mathcal P(R_0,s,\alpha)<4\pi$.
Proposition~\ref{prop:central-weight-reduction} gives the whole weight interval.
\end{proof}

\subsection{\texorpdfstring{Product estimates for \(1/4\le s\le1/2\)}{Product estimates for 1/4 <= s <=1/2}}
\label{sec:whole-low-radius-product}

\begin{theorem}\label{thm:whole-low-radius-product}
For every $9/20\le R\le R_*=1936/2561$, every
$1/4\le s\le1/2$, and every $0\le\alpha\le1$,
\[
 \mathcal P(R,s,\alpha)<4\pi.
\]
\end{theorem}

The proof gives a bound independent of the amplitude parameter. A convexity argument then leaves two elementary
polynomials of degree six.

\subsubsection{A bound for the full conformal kernel}

\begin{lemma}[The complete circular average]\label{lem:full-average-series}
Let \(a,b\ge0\), \(a+b\le1\), \(0\le\beta<1\), and \(u=1-\beta^2\).
For \(Q=a\cos U+b\cos V\), put
\[
 \lambda=\beta a,\qquad \mu=\beta b,\qquad
 w=1-\lambda^2-\mu^2,\qquad
 G=w^2-4\lambda^2\mu^2,\qquad z=\frac{4\lambda^2\mu^2}{w^2}.
\]
Then \(w>0\), \(G>0\), \(0\le z<1\), and
\[
 \mathbb E\frac{u}{(1-\beta Q)^2}
 =\frac{u\sqrt w}{G}\,\mathscr C(z),
 \qquad
 \mathscr C(z)=1-\sum_{n=1}^{\infty}\delta_nz^n.
\]
Here
\[
 \gamma_0=1,\qquad
 \gamma_n=\gamma_{n-1}\left(1-\frac1{16n^2}\right),\qquad
 \delta_n=\frac{\gamma_{n-1}}{16n^2}>0.
\]
In particular, \(\mathscr C(z)\le1-z/16\le1\).
\end{lemma}
\begin{proof}
Since \(\lambda+\mu<1\), both
\(1-(\lambda+\mu)^2\) and \(1-(\lambda-\mu)^2\) are positive,
and their product is \(G\).
Condition on the difference of the two angles and combine the cosines.
The identity
\[
 \frac1{2\pi}\int_0^{2\pi}\frac{d\theta}{(1-r\cos\theta)^2}
 =(1-r^2)^{-3/2}\qquad(|r|<1)
\]
follows from the substitution $t=\tan(\theta/2)$ in the average
of $(1-r\cos\theta)^{-1}$, which gives $(1-r^2)^{-1/2}$;
adding $r$ times its derivative to that average gives the displayed
identity. For fixed difference of the angles, the combined cosine has
amplitude $(\lambda^2+\mu^2+2\lambda\mu\cos\theta)^{1/2}$.
The difference and the remaining angle are independent Haar angles
(the change of variables is an automorphism of the two-torus). Thus
\[
 \mathbb E\frac{u}{(1-\beta Q)^2}
 =u\,\mathbb E_\theta(w-2\lambda\mu\cos\theta)^{-3/2}.
\]
Expanding the last integrand by the binomial series and integrating the
even cosine powers yields
\[
 \mathbb E_\theta(1-\sqrt z\cos\theta)^{-3/2}
 =\sum_{n=0}^{\infty}\gamma_nz^n.
\]
The series converges uniformly on each compact subinterval of \(0\le z<1\).
In fact, writing $(a)_k=a(a+1)\cdots(a+k-1)$ for the rising
factorial, we have $\gamma_n=(3/2)_{2n}/(4^n(n!)^2)$. Hence
\[
 \frac{\gamma_n}{\gamma_{n-1}}
 =\frac{(2n-1/2)(2n+1/2)}{4n^2}=1-\frac1{16n^2}.
\]
Multiplication by \(1-z=G/w^2\) therefore gives
\[
 (1-z)\sum_{n=0}^{\infty}\gamma_nz^n
 =1-\sum_{n=1}^{\infty}(\gamma_{n-1}-\gamma_n)z^n
 =\mathscr C(z),
\]
which proves the formula and its upper bounds.
\end{proof}

\begin{proof}[Proof of Theorem~\ref{thm:whole-low-radius-product}]
Put $\beta=2s/(1+s^2)$, ${\omega}=\beta^2$, $u=1-{\omega}$, $h=1-R$,
$d=2R-1$, ${\varkappa}=4Rh$, and $x=\alpha^2$. Set
\[
 \lambda=\beta\sqrt{Rx},\qquad
 \mu=\beta\sqrt{h(1-x)},\qquad
 w=1-\lambda^2-\mu^2.
\]
Conditioning on $U-V$, combining the two cosines, and evaluating the
remaining circular average gives
\[
 {\mathfrak F}:=\mathbb E \chi_s(Q_{R,\alpha})^2
 =u\,\mathbb E_\theta(w-2\lambda\mu\cos\theta)^{-3/2}.
\]

Lemma~\ref{lem:full-average-series} gives
\[
 {\mathfrak F}\le\frac{u\sqrt w}{G},\qquad G=w^2-4\lambda^2\mu^2.
\]
With
\[
 y={\omega}(x-h)-d,\qquad w_0=\frac{{\varkappa}(1+u)}2,
\]
direct expansion gives the exact identities
\begin{equation}\label{eq:low-radius-square}
 G=y^2+{\varkappa}u,\qquad w=w_0-dy.
\end{equation}
The tangent inequality for the square root gives
$\sqrt w\le\sqrt{w_0}(1-{\eta_{\mathrm{tan}}} y)$, where ${\eta_{\mathrm{tan}}}=d/(2w_0)$.
For $K>0$ and real ${\eta_{\mathrm{tan}}}$,
\[
 \sup_{y\in\mathbb R}\frac{1-{\eta_{\mathrm{tan}}} y}{y^2+K}
 =\frac{1+\sqrt{1+{\eta_{\mathrm{tan}}}^2K}}{2K}.
\]
Indeed, subtracting the left rational function from the claimed upper
bound and clearing its positive denominator gives a quadratic with
zero discriminant. Taking $K={\varkappa}u$ and using
$\sqrt{1+t}\le1+t/2$ proves
\[
 {\mathfrak F}\le\sqrt{\frac{1+u}{2{\varkappa}}}
       \left(1+\frac{d^2u}{4{\varkappa}(1+u)^2}\right).
\]
In terms of \(s\), this becomes
\begin{align}
 {\mathfrak F}&\le\frac{\sqrt{1+s^4}}{\sqrt {\varkappa}(1+s^2)}
 \left[1+\frac{d^2}{16{\varkappa}}
       \left(\frac{1-s^4}{1+s^4}\right)^2\right]\notag\\
 &\le\frac{1+e(R)}{\sqrt {\varkappa}}\,\frac{1+s^4/2}{1+s^2},
 \qquad e(R)=\frac{(2R-1)^2}{64R(1-R)}.
 \label{eq:low-radius-alpha-free}
\end{align}
Both inequalities are uniform in $\alpha$.

\subsubsection{The symmetric cutoff and convexity in the weight}
\label{sec:product-weight-convexity}

Write $c_R=L/\Abar(R)$ and
$d_s=((1-s)/(1+s))^2$. The bounds in Section~\ref{sec:global-majorants} give
$4/9<c_R<3/5$ on the actual weight range. For $1/4\le s\le1/2$,
$d_s\le9/25<c_R$.
Pair $q,-q$ in the symmetric distribution of $Q_{R,\alpha}$. If both
kernel values are below $c_R$, their sum is at least $2u$ by convexity.
If both are above it, their truncated sum is $2c_R$. In the remaining
case this sum is at least $c_R+d_s$.

The function $2u-d_s$ decreases on $[1/4,1/2]$. Its derivative has the
sign of
\[
 (1+s^2)^3-4s(1+s)^4<0,
\]
since the first term is at most $125/64$ and the second is at least
$625/256$. Its minimum is
$18/25-1/9=137/225>3/5$. Thus all three cases give
\[
 \mathbb E\min\{\chi_s(Q_{R,\alpha})^2,c_R\}
 \ge\frac{c_R+d_s}{2}.
\]
Together with \eqref{eq:low-radius-alpha-free}, this proves
\begin{equation}\label{eq:low-radius-H}
 \mathcal P(R,s,\alpha)\le\pi^2 H_s(R)-\frac L2,
\end{equation}
where
\[
 H_s(R)={G_R}(1+e(R))\mathscr H(s)-\sqrt{D(R)}\,d_s,
\]
\[
 {G_R}=\sqrt{3+1/R},\quad D(R)=(1-R)(1+3R),\quad
 \mathscr H(s)=\frac{1+s^4/2}{1+s^2}.
\]
For every fixed $s$ in the range of this estimate, the function $H_s$
is convex on $(0,1)$. In fact, the calculation below also proves that
$G_R(1+e(R))$ is convex as a function of $R$ throughout $(0,1)$.
Indeed,
\[
 {G_R}'=-\frac1{2R^2{G_R}},\quad
 {G_R}''=\frac{3(4R+1)}{4R^4{G_R}^3}>0,
\]
\[
 e'=\frac{2R-1}{64R^2(1-R)^2},\qquad
 e''=\frac{1-3R(1-R)}{32R^3(1-R)^3}>0.
\]
For $R\le1/2$, the cross term $2e'{G_R}'$ in $(e{G_R})''$ is nonnegative.
For $R\ge1/2$, the inequality $e''{G_R}+2e'{G_R}'>0$ is equivalent to
\[
 2[1-3R(1-R)](1+3R)>(2R-1)(1-R).
\]
The left side is at least $5/4$, and the right side is at most $1/2$.
Since $e(R)\ge0$, these inequalities show $(eG_R)''\ge0$.
Thus $G_R(1+e(R))$ is convex on $(0,1)$. Finally $\sqrt D$ is concave,
$\mathscr H(s)>0$ and $d_s\ge0$, proving the assertion for $H_s$.

Because $R_*<4/5$, it remains to bound the two endpoint values of
$H_s$. The endpoint $4/5$ is used only for this elementary convex
comparison; no cutoff estimate beyond $R_*$ is asserted.
We shall prove
\[
 H_s(9/20),\ H_s(4/5)<K_0:=\frac{1791}{1000}
 <\frac2{\sqrt3}+\frac2\pi.
\]
For the final comparison use
$\sqrt3<97/56$, $\pi<22/7$, and
$112/97+7/11=1911/1067>1791/1000$.

\subsubsection{Two short polynomial comparisons}
\label{sec:low-radius-endpoints}

At $R=9/20$, one has
${G_R}=\sqrt{47}/3$, $e=1/1584$, and $\sqrt D=\sqrt{517}/20$.
The rational square comparisons give
\[
 {G_R}(1+e)<\frac{22867}{10000},\qquad
 \sqrt D>\frac{11368}{10000}.
\]
For example, use $\sqrt{47}<68557/10000$ and
$(68557/10000)(1585/4752)<22867/10000$ for the first;
$(22736/1000)^2<517$ proves the second.
After multiplication by $(1+s^2)(1+s)^2$, it suffices to prove
positivity of the polynomial ${P_{\mathrm{end}}}$ defined by
\[
 20000{P_{\mathrm{end}}}(s)=12822-65300s+71378s^2+26168s^3
                 +35689s^4-45734s^5-22867s^6.
\]
For $1/4\le s\le1/2$,
\begin{align*}
 {P_{\mathrm{end}}}''(s)&\ge\frac1{20000}\left(
 142756+\frac{157008}{4}+\frac{428268}{16}
              -\frac{914680}{8}-\frac{686010}{16}\right)\\
 &=\frac{412513}{160000}>2.
\end{align*}
At $a=3/8$,
\[
 {P_{\mathrm{end}}}(a)=\frac{14419029}{5242880000}>\frac1{400},\qquad
 |{P_{\mathrm{end}}}'(a)|=\frac{330767056}{5242880000}<\frac1{10}.
\]
Writing $t=s-a$, Taylor's formula and completion of the square give
\[
 {P_{\mathrm{end}}}(s)\ge {P_{\mathrm{end}}}(a)+{P_{\mathrm{end}}}'(a)t+t^2
 \ge {P_{\mathrm{end}}}(a)-\frac{{P_{\mathrm{end}}}'(a)^2}{4}
 >\frac1{400}-\frac{(1/10)^2}{4}=0.
\]

At $R=4/5$, one has
${G_R}=\sqrt{17}/2$, $e=9/256$ and $\sqrt D=\sqrt{17}/5$.
Since $4<\sqrt{17}<33/8$,
\[
 {G_R}(1+e)=\frac{265\sqrt{17}}{512}<\frac{107}{50},\qquad
 \sqrt D>\frac45.
\]
The corresponding cleared polynomial ${Q_{\mathrm{end}}}$ is
\[
 2000{Q_{\mathrm{end}}}(s)=902-4596s+6084s^2+3964s^3
                 +3042s^4-4280s^5-2140s^6.
\]
Its second derivative satisfies
\[
 {Q_{\mathrm{end}}}''(s)=\frac{1521}{250}
 +s\left(\frac{2973}{250}-\frac{214}{5}s^2\right)
 +s^2\left(\frac{4563}{250}-\frac{321}{10}s^2\right)>6,
\]
because both brackets are positive for $s\le1/2$.
At $a=3/8$,
\[
 {Q_{\mathrm{end}}}(a)=\frac{17404241}{131072000}>\frac18,\qquad
 |{Q_{\mathrm{end}}}'(a)|=\frac{7219497}{8192000}<1.
\]
Therefore
\[
 {Q_{\mathrm{end}}}(s)\ge {Q_{\mathrm{end}}}(a)+{Q_{\mathrm{end}}}'(a)t+3t^2
 \ge {Q_{\mathrm{end}}}(a)-\frac{{Q_{\mathrm{end}}}'(a)^2}{12}
 >\frac18-\frac1{12}=\frac1{24}>0.
\]
This establishes both endpoint comparisons. Convexity and
\eqref{eq:low-radius-H} now yield
\[
 \mathcal P(R,s,\alpha)
 <\pi^2\left(\frac2{\sqrt3}+\frac2\pi\right)-\frac L2
 =4\pi,
\]
proving Theorem~\ref{thm:whole-low-radius-product}.
\end{proof}

\subsection{\texorpdfstring{Product estimates for \(1/2\le s\le3/5\)}{Product estimates for 1/2 <= s <= 3/5}}\label{sec:product-half-three-fifths}

\begin{theorem}[Moderate-radius product bound]\label{thm:wide-product-three-fifths}
For $4/9\le R\le R_*$, $1/2\le s\le3/5$ and $0\le\alpha\le1$,
$\mathcal P(R,s,\alpha)<4\pi$.
\end{theorem}
\begin{proof}
Use the parameters
\[
 u=\left(\frac{1-s^2}{1+s^2}\right)^2,\quad {\omega}=1-u,
 \quad d_s=\left(\frac{1-s}{1+s}\right)^2,
 \quad \eta=\min\{R,1-R\},
\]
\[
 D=(1-R)(1+3R),\quad S=\sqrt D,\quad
 {G_R}=\sqrt{3+1/R},\quad e=\frac{(2R-1)^2}{64R(1-R)},
\]
\[
 \mathscr H(s)=\frac{1+s^4/2}{1+s^2},\qquad
 C_R=\frac{1+e}{\sqrt{4R(1-R)}}.
\]
Estimate~\eqref{eq:low-radius-alpha-free} gives
${\mathfrak F}:=\mathbb E \chi_s(Q)^2\le C_R\mathscr H(s)$.

\paragraph{A quadratic lower bound for the removed part.}
Set $c=c_R=L/\Abar(R)$ and $B=c+d_s-2u$. On this strip,
$u\le9/25<4/9<c$, so for $q\ge0$ one has $f(-q)<c$, where
$f(q)=u/(1-\sqrt {\omega}\,q)^2$.
Before $f(q)$ reaches $c$,
\[
 f(q)+f(-q)
 =2u\sum_{j\ge0}(2j+1){\omega}^jq^{2j}
 \ge2u+6u{\omega}q^2.
\]
After it reaches $c$, the sum of the two minima is at least
$c+d_s$. If $B\ge0$, then $B<1<6u{\omega}$, because
\[
 6u{\omega}\ge\frac{86400}{83521}>1,
 \qquad c+d_s<\frac35+\frac19<1.
\]
Consequently the paired minimum is at least $2u+Bq^2$.
Since $\mathbb E Q^2\ge\eta/2$, this gives
\[
 \mathbb E\min\{f(Q),c\}\ge u+\frac\eta4 B
 \quad(B\ge0).
\]
If $B<0$, the paired minimum is at least $c+d_s$, giving instead
$\mathbb E\min\{f(Q),c\}\ge(c+d_s)/2$.
Thus, throughout the strip,
\[
 \mathcal P\le\max\{\mathcal U_1(R,s),\mathcal U_2(R,s)\},
\]
where
\begin{align*}
 \mathcal U_1&=\Abar C_R\mathscr H(s)-\Abar u
       -\frac\eta4\bigl(L+\Abar d_s-2\Abar u\bigr),\\
 \mathcal U_2&=\Abar C_R\mathscr H(s)-\frac L2-\frac{\Abar d_s}2.
\end{align*}

\paragraph{Monotonicity in the centre parameter.}
For $4/9\le R\le19/25$, one has $1\le C_R<6/5$.
Indeed, $4R(1-R)\ge456/625>(64/75)^2$ and
$e\le169/7296<3/125$.
For $1/2\le s\le3/5$,
\[
 \mathscr H'(s)=-\frac{s(2-2s^2-s^4)}{(1+s^2)^2},\qquad
 u'(s)=-\frac{8s(1-s^2)}{(1+s^2)^3},\qquad
 d_s'=-\frac{4(1-s)}{(1+s)^3}.
\]
Elementary endpoint bounds yield
\[
 \frac{719}{2312}\le|\mathscr H'|\le\frac{69}{125},\qquad
 |u'|\ge\frac{5000}{4913},\qquad
 -\frac{d_s'}2\le\frac8{27}.
\]
Therefore
\[
 \frac{\partial_s\mathcal U_2}{\Abar}
 =C_R\mathscr H'-d_s'/2
 \le-\frac{719}{2312}+\frac8{27}<0,
\]
and, since $\eta\le1/2$,
\[
 \frac{\partial_s\mathcal U_1}{\Abar}
 =C_R\mathscr H'-(1-\eta/2)u'-\eta d_s'/4
 \ge-\frac{414}{625}+\frac34\frac{5000}{4913}
 =\frac{309768}{3070625}>0.
\]
\paragraph{Convexity in the weight.}
At a fixed $s$ write
\[
 \mathcal U_1=\pi^2\left[{G_R}(1+e)\mathscr H(s)-S W\right]-\frac{\eta L}4,
 \qquad W=(2-\eta)u+\frac\eta2d_s.
\]
The function ${G_R}(1+e)$ is convex by the calculation in
Section~\ref{sec:product-weight-convexity}.
For $R\ge1/2$, $W=(1+R)u+(1-R)d_s/2$ is positive and affine,
with $W'=u-d_s/2>0$; since $S'<0$ and $S''<0$, $(SW)''<0$.
For $4/9\le R\le1/2$, put $t=u-d_s/2>0$. Then
$W'= -t$, $W\ge3u/2$, and
\[
 (SW)''<\frac{-3W+2t(3R-1)}S<0,
\]
using $S''<-3/S$ and $3R-1\le1/2$.
Thus $\mathcal U_1(\cdot,s)$ is convex on either side of $R=1/2$.
At $s=3/5$, its first derivative has the positive jump
\[
 \partial_R\mathcal U_1(1/2+,s)-\partial_R\mathcal U_1(1/2-,s)
 =\frac{\Abar}2(c+d_s-2u)>0,
\]
because $u=64/289$, $d_s=1/16$ and $c>4/9>128/289$.
Hence $\mathcal U_1(\cdot,3/5)$ is convex on $[4/9,19/25]$.
Also $\mathcal U_2(\cdot,s)$ is convex on this interval because $G_R(1+e)$ is convex and $S$ is concave. It therefore suffices to check $\mathcal U_1(\cdot,3/5)$ and $\mathcal U_2(\cdot,1/2)$ at $R=4/9,19/25$. The latter endpoint is used solely for the auxiliary upper functions.

\paragraph{Endpoint comparisons.}
The following bounds follow by squaring (all quantities are positive):
\[
\begin{array}{c|c|c|c|c}
 R&\eta&e&G_R(1+e)\text{ upper bound}&S\text{ lower bound}\\\hline
 4/9&4/9&1/1280&22931/10000&2277/2000\\
 19/25&6/25&169/7296&1063/500&887/1000
\end{array}
\]
For example, the exact radical expressions in the first row are
$1281\sqrt{21}/2560$ and $\sqrt{105}/9$; in the second row they are
$(7465/7296)\sqrt{82/19}$ and $\sqrt{492}/25$.
At $s=3/5$, $\mathscr H=1331/1700$, $u=64/289$, $d_s=1/16$, and
$W=7457/20808$ or $45923/115600$, respectively. At $s=1/2$,
$\mathscr H=33/40$ and $d_s=1/9$.
The preceding bounds, $\pi<22/7$, and $1/\sqrt3>56/97$ give the following
strict rational comparisons, using
\begin{align*}
 \mathcal U_1/\pi&=\pi[G_R(1+e)\mathscr H-SW-\eta/\sqrt3]+\eta,\\
 \mathcal U_2/\pi&=\pi[G_R(1+e)\mathscr H-Sd_s-2/\sqrt3]+2.
\end{align*}
\[
\begin{array}{c|c|c}
 R& \mathcal U_1(R,3/5)/\pi\text{ upper bound}&\mathcal U_2(R,1/2)/\pi\text{ upper bound}\\\hline
 4/9&7061299425947/1766079000000&532233241/135800000\\
 19/25&154178073569/39246200000&218395057/61110000
\end{array}
\]
For each entry, replace $G_R(1+e)$ by its upper bound, $S$ and
$1/\sqrt3$ by their lower bounds and $\pi$ by $22/7$ in the
positive upper bracket, then add $\eta$ or $2$ as appropriate.
Every displayed fraction is less than $4$. If an original bracket
were nonpositive, the required inequality would be immediate.
Convexity in $R$ and monotonicity in $s$ complete the proof.
\end{proof}

\subsection{A fixed two-cosine estimate}\label{sec:fixed-two-cosine}

This fixed-amplitude reduction supplies the low-weight product range.
The logarithmic argument proving Theorem~\ref{thm:fixed-product-high}
also supplies the endpoint used in
Theorem~\ref{thm:fixed-product-mid}.

\subsubsection{Reduction to one fixed pair of amplitudes}

Write $\Aeq=4\pi^2/\sqrt3$ and
$c_{\rm eq}=1-\sqrt3/\pi$. The normalized height of a fundamental
embedding can, after an integral change of angular variables and phase
shifts, be written
\[
 q=a_1\cos U+a_2\cos V+a_3\cos(U+V+\phi),\qquad
 a_i=\sqrt{r_i}|p_i|\ge0.
\]
Cauchy--Schwarz gives $M:=a_1+a_2+a_3\le1$.

\begin{proposition}\label{prop:fixed-amplitudes}
If $R=\max r_i\le4/9$, then, for every $0\le s<1$,
\[
 J_s(\weights,p)\le \Aeq\,\mathbb E
 \left(\chi_s\left(\frac{\cos U+2\cos V}{3}\right)^2
                 -c_{\rm eq}\right)_+.
\]
Here $U,V$ are independent uniform angles.
\end{proposition}
\begin{proof}
Set $f(q)=(\chi_s(q)^2-c_{\rm eq})_+$, a convex function. Replacing the
actual area and cutoff by $\Aeq$ and $c_{\rm eq}$ increases the
upper bound for $J$. Separate a largest amplitude $a=\max a_i$.
For each fixed angle in its mode, the other two modes combine into a
cosine of amplitude at most $b=M-a$. Circular averaging of a convex
function is nondecreasing in the nonnegative amplitude, so the resulting
average is at most
\[
 \Phi(a,b):=\mathbb E f(a\cos U+b\cos V).
\]
Now $a\le\sqrt R\le2/3$, $b\le2M/3\le2/3$, and $a+b\le1$.
The function $\Phi$ is convex, symmetric, separately even, and
nondecreasing in each nonnegative variable. The vertices of the polygon
\[
 a,b\ge0,\qquad a,b\le2/3,\qquad a+b\le1
\]
are $(0,0),(2/3,0),(2/3,1/3),(1/3,2/3),(0,2/3)$.
Convexity and monotonicity therefore give
$\Phi(a,b)\le\Phi(1/3,2/3)$.
\end{proof}

\begin{theorem}\label{thm:fixed-product-high}
For every fundamental embedding with $R\le4/9$ and every
$3/5\le s\le41/64$, one has $J_s(\weights,p)<4\pi$.
\end{theorem}

\begin{proof}[Proof of Theorem~\ref{thm:fixed-product-high}]
By Proposition~\ref{prop:fixed-amplitudes}, it suffices to estimate the
fixed two-cosine distribution. We establish its density bound, monotonicity
and the required endpoint comparisons.
\subsubsection{Density and an elementary logarithmic bound}

Put $Q=(\cos U+2\cos V)/3$, $\delta=1/3$, $\kappa=9/8$, and
$\varrho_*=3/(2\pi\sqrt2)$. For $\delta<q<1$, the density of $Q$ is
\begin{equation}\label{eq:fixed-density}
 \varrho(q)=\varrho_*\sum_{n=0}^{\infty}a_nz(q)^n,
 \qquad z(q)=\kappa(1-q^2),\qquad
 a_n=\frac{\binom{2n}{n}^2}{16^n}.
\end{equation}
Indeed, the angular substitution
$\theta=(U+V)/2$, $\varphi=(U-V)/2$ shows that $Q$ has the
same distribution as
$\cos\theta\cos\varphi-\delta\sin\theta\sin\varphi$.
Conditioning on $\varphi$ gives the amplitude
$\sqrt{1-(1-\delta^2)\sin^2\varphi}$. Substituting
$\sin\varphi=\sqrt{z}\sin\psi$ in its density yields
\[
 \varrho(q)=\frac{2}{\pi^2\sqrt{1-\delta^2}}
 \int_0^{\pi/2}\frac{d\psi}{\sqrt{1-z\sin^2\psi}},
\]
and the binomial series gives \eqref{eq:fixed-density}.

For every $n\ge1$, one has $a_n<1/(\pi n)$. To see this, set
$b_n=\binom{2n}{n}/4^n$ and
$I_j=\int_0^{\pi/2}\sin^j x\,dx$. Integration by parts gives
$I_{2n}=(\pi/2)b_n$ and $I_{2n-1}=1/(2nb_n)$.
The strict inequality $I_{2n}<I_{2n-1}$ proves the claim.
Consequently, with $C=106/333>1/\pi$ and $N\ge0$,
\begin{equation}\label{eq:fixed-RN}
 \varrho(q)\le \varrho_*\mathcal R_N(q),\qquad
 \mathcal R_N(q)=1-C\log(1-z(q))+
       \sum_{n=1}^N\left(a_n-\frac Cn\right)z(q)^n.
\end{equation}
Only $N=1,4$ are needed. Their coefficients are determined by
\[
 a_1=\frac14,\quad a_2=\frac9{64},\quad
 a_3=\frac{25}{256},\quad a_4=\frac{1225}{16384}.
\]

\subsubsection{Monotonicity and the endpoint}
\label{sec:fixed-high-monotonicity}

Set $c=c_{\rm eq}$, $\beta=2s/(1+s^2)$, $u=1-\beta^2$, and
\[
 I(\beta)=\int_{q_0}^1
  \left(\frac{u}{(1-\beta q)^2}-c\right)\varrho(q)\,dq,
 \qquad q_0=q_0(s,c)=\frac{1-\sqrt{u/c}}{\beta}.
\]
Here we use the truncation height \eqref{eq:truncation-height}.
The two rows below ensure $\delta<q_0<\beta$ on the stated range.
Differentiation, using the vanishing of the bracket at $q_0$, gives
\[
 I'(\beta)=2\int_{q_0}^1
           \frac{q-\beta}{(1-\beta q)^3}\varrho(q)\,dq.
\]
By \eqref{eq:fixed-density}, $\varrho\ge \varrho_*$ on this integration range,
so the positive portion $q\ge\beta$ is at least $\varrho_*/u$. On $q\le\beta$,
\[
 \partial_\beta\frac{\beta-q}{(1-\beta q)^3}
 =\frac{1+2\beta q-3q^2}{(1-\beta q)^4}>0.
\]
Define
\[
 \mathcal N(l,b)=2\int_l^b
      \frac{b-q}{(1-bq)^3}\mathcal R_1(q)\,dq.
\]
If $\beta_a\le\beta\le\beta_b$ and $q_0\ge l>\delta$, then
\begin{equation}\label{eq:fixed-derivative-lower}
 I'(\beta)\ge \varrho_*
 \left(\frac1{1-\beta_a^2}-\mathcal N(l,\beta_b)\right).
\end{equation}
The required constants are
\begingroup
\renewcommand{\arraystretch}{2.2}
\[
\begin{array}{c|c|c|c|c|c}
 s\text{-interval}&\beta_a&\beta_b&l&
 \mathcal N(l,\beta_b)<&1/(1-\beta_a^2)\ge\\\hline
 {[3/5,31/50]}&15/17&3100/3461&337/1000&89/20&289/64\\
 {}[31/50,41/64]&3100/3461&5248/5777&3/8&49/10&5
\end{array}
\]
\endgroup
Both differences in \eqref{eq:fixed-derivative-lower} are positive.
The elementary bounds $\sqrt3<97/56$ and $\pi>333/106$ give
$c>4183/9324$. At the left endpoints the values of
$u/(1-\beta l)^2$ are, respectively,
\[
 \frac{2560000}{5707321}<\frac{4183}{9324},\qquad
 \frac{9474084}{21132409}<\frac{4183}{9324}.
\]
Thus $q_0>l$ there. Rewriting \eqref{eq:truncation-height} as
\[
 q_0=\frac{(1+c^{-1/2})s-(c^{-1/2}-1)/s}{2}
\]
shows that $q_0$ increases with $s$, proving the condition throughout
each row. Also $q_0<\beta$, since $cu<1$.

At $s=41/64$, put
\[
 b=\frac{5248}{5777},\qquad u=1-b^2,\qquad l=\frac{331}{800},
 \qquad c_*:=\frac{u}{(1-bl)^2}
 =\frac{405015625}{902822209}<\frac{4183}{9324}<c.
\]
Replacing $c$ by $c_*$ and using \eqref{eq:fixed-RN} with $N=4$ gives
\[
 I(b)\le \varrho_*\mathcal T,\qquad
 \mathcal T=\int_l^1
       \left(\frac{u}{(1-bq)^2}-c_*\right)\mathcal R_4(q)\,dq
       <\frac{16327}{10000}.
\]
The explicit elementary evaluations are given in
Section~\ref{sec:fixed-log-evaluations}. Since
$3\cdot16327^2=799712787<800000000$, the last fraction is smaller than
$4/\sqrt6$. Thus
\[
 \Aeq{}I(b)<\Aeq\varrho_*\frac4{\sqrt6}=4\pi.
\]
Strict monotonicity and Proposition~\ref{prop:fixed-amplitudes} prove
Theorem~\ref{thm:fixed-product-high}.

\subsubsection{Explicit evaluation of the three constants}
\label{sec:fixed-log-evaluations}

All the evaluations reduce to logarithms of positive rational numbers.
For rational $l<r<1/b$, put $w_0=1-bl$, $w_1=1-br$, and
\[
 M_n(l,r;b)=b^{-n-1}\sum_{j=0}^n(-1)^j\binom njD_j,
\]
where
\[
 D_0=\frac1{w_1}-\frac1{w_0},\quad
 D_1=\log\frac{w_0}{w_1},\quad
 D_j=\frac{w_0^{j-1}-w_1^{j-1}}{j-1}\quad(j\ge2).
\]
Substitution $w=1-bq$ gives
$M_n=\int_l^r q^n/(1-bq)^2\,dq$. Also put
\[
 \mathcal J_a(l,r;b)=
 \frac1{(1-ba)^2}\log\frac{(r-a)w_0}{(l-a)w_1}
 +\frac{1/w_1-1/w_0}{1-ba}.
\]
Partial fractions give
$\mathcal J_a=\int_l^r[(q-a)(1-bq)^2]^{-1}\,dq$.
Writing $d_1=1/4-C$ and $W(q)=(1-q^2)/(1-bq)^2$, integration by parts
now yields, with all $M,\mathcal J$ on $[l,b]$,
\begin{align*}
 \mathcal N(l,b)={}&W(b)\mathcal R_1(b)-W(l)\mathcal R_1(l)\\
 &+2C\left[-M_1+\frac{1-\delta^2}{2}(\mathcal J_\delta+\mathcal J_{-\delta})\right]
 +\frac94d_1(M_1-M_3).
\end{align*}

For the endpoint use the parameters $b,l,u,c_*$ specified in
Section~\ref{sec:fixed-high-monotonicity}, set $r=1$, and define
\[
 A_n=\kappa^n\sum_{j=0}^n(-1)^j\binom nj
 \left[uM_{2j}-c_*\frac{1-l^{2j+1}}{2j+1}\right]
 \quad(0\le n\le4).
\]
Put $\mathscr V(q)=1-z(q)=\kappa(q^2-\delta^2)$ and
\[
 R_- =\frac{(1-\delta)w_0}{(l-\delta)w_1},\qquad
 R_+ =\frac{(1+\delta)w_0}{(l+\delta)w_1}.
\]
A further integration by parts gives
\begin{align*}
 H&:=\int_l^1\left(\frac{u}{(1-bq)^2}-c_*\right)\log \mathscr V(q)\,dq\\
 &=\left(c_*l-\frac{u}{bw_0}\right)\log \mathscr V(l)
 +\left(-\frac{u}{b(1-b\delta)}+c_*\delta\right)\log R_-\\
 &\quad+\left(-\frac{u}{b(1+b\delta)}-c_*\delta\right)\log R_+
 +2c_*(1-l),
\end{align*}
and consequently
\[
 \mathcal T=A_0-CH+
             \sum_{n=1}^4\left(a_n-\frac Cn\right)A_n.
\]

For completeness the following table supplies all logarithm bounds.
An integer $n$ in its last column means
$n/10^{12}<\log v<(n+1)/10^{12}$. Every entry follows from the same
finite positive series: for $1\le y\le2$, with $x=(y-1)/(y+1)$,
\[
 2\sum_{j=0}^{20}\frac{x^{2j+1}}{2j+1}
 \le\log y\le
 2\sum_{j=0}^{20}\frac{x^{2j+1}}{2j+1}
       +\frac{2x^{43}}{43(1-x^2)}.
\]
For arbitrary $v>0$, write $v=2^ky$ with $1\le y\le2$ and use the
same bounds for $\log2$, reversing endpoints when $k<0$.

\begin{center}\small
\begin{tabular}{lll}
Case&$v$&$n$\\\hline
First derivative&$22121/8000000$&$-5890669437074$\\
&$74511479/95828168$&$-251603476934$\\
&$83628143/23685210$&$1261529299488$\\
&$14108775700/26053731$&$6294391290320$\\
&$30834404300/4763095731$&$1867733268797$\\[2pt]
Second derivative&$17/512$&$-3405111280984$\\
&$214499807/266989832$&$-218901736678$\\
&$22004593/5832225$&$1327852630086$\\
&$303714424/5832225$&$3952689290952$\\
&$655787912/99147825$&$1889225513526$\\[2pt]
Endpoint&$90141/13225$&$1919266133310$\\
&$346049/5120000$&$-2694329334506$\\
&$5769024/102097$&$4034334852633$\\
&$11538048/948497$&$2498526748493$\\
\end{tabular}
\end{center}
The five entries for each derivative are, in order,
$\mathscr V(l),\mathscr V(b),w_0/w_1,(b-\delta)w_0/((l-\delta)w_1)$ and
$(b+\delta)w_0/((l+\delta)w_1)$.
The endpoint entries are $w_0/w_1,\mathscr V(l),R_-,R_+$.
Substitution gives the deliberately coarser rational bounds
\[
 \mathcal N(337/1000,3100/3461)<\frac{4405}{1000}<\frac{89}{20},
 \qquad
 \mathcal N(3/8,5248/5777)<\frac{4868}{1000}<\frac{49}{10},
\]
\[
 \begin{gathered}
 A_0<\frac{1529585}{10^6},\quad
 A_1>\frac{300364}{10^6},\quad
 A_2>\frac{113249}{10^6},\\
 A_3>\frac{56946}{10^6},\quad
 A_4>\frac{33397}{10^6},\quad
 -H<\frac{396096}{10^6}.
 \end{gathered}
\]
As $a_n-C/n<0$, these imply the exact rational comparison
\[
 \mathcal T<\frac{1781242254949}{1091174400000}
 <\frac{16327}{10000},
\]
whose last difference is $318187931/1091174400000>0$.
This completes all finite arithmetic used in the strip proof.
\end{proof}

\subsubsection{The fixed-amplitude estimate for \texorpdfstring{$1/2\le s\le3/5$}{1/2 <= s <= 3/5}}
\label{sec:fixed-mid}

\begin{theorem}[The fixed-amplitude bound]\label{thm:fixed-product-mid}
Let $U,V$ be independent uniform angles, and put
$Q=(\cos U+2\cos V)/3$, $c=1-\sqrt3/\pi$, and
$\Aeq=4\pi^2/\sqrt3$. Then
\[
 \Aeq\,\mathbb E(\chi_s(Q)^2-c)_+<4\pi
 \qquad\left(\frac12\le s\le\frac35\right).
\]
Consequently, for every fundamental embedding with $R\le4/9$,
every unit direction $p$, and every $1/2\le s\le3/5$, one has
$J_s(\weights,p)<4\pi$.
\end{theorem}

\begin{proof}
The assertion about $J$ follows from Proposition~\ref{prop:fixed-amplitudes}.
We prove the fixed-amplitude estimate.

\paragraph{A short rational majorant of the full average.}
Put
\[
 \beta=\frac{2s}{1+s^2},\quad u=1-\beta^2,
 \quad {\mathfrak F}=\mathbb E\frac{u}{(1-\beta Q)^2}.
\]
For the fixed amplitudes, Lemma~\ref{lem:full-average-series} gives
\[
 w=\frac{4+5u}{9},\qquad G=\frac{u(8+u)}9,
 \qquad z=\frac{16(1-u)^2}{(4+5u)^2}.
\]
Since $\mathscr C(z)=1-\sum_{n\ge1}\delta_nz^n$ with $\delta_n>0$ and $\delta_1=1/16$,
\[
 {\mathfrak F}\le\frac{3\sqrt{4+5u}}{8+u}\left(1-\frac z{16}\right)
 =\frac{9(1+2u)(5+4u)}{(8+u)(4+5u)^{3/2}}.
\]
For $1/2\le s\le23/40$, this implies the simpler bound
\begin{equation}\label{eq:fixed-rational-F}
 {\mathfrak F}\le {\mathfrak F}_*(s):=1-\frac23s^2+\frac29s^4
                    +\frac{27s^6}{4(10+21s^2)}.
\end{equation}
Here is an explicit algebraic proof. Set $t=s^2$ and
\[
 B=9-2t+9t^2,\quad C_0=3-2t+3t^2,\quad
 D_0=9+2t+9t^2,\quad E=9+14t+9t^2.
\]
Since $1/4\le t\le529/1600<1/3$,
\[
 \sqrt B\ge P:=3-t/3+40t^2/27,
 \qquad B-P^2=\frac{80}{729}t^3(9-20t)>0.
\]
Thus the preceding bound is at most $9(1+t)C_0D_0/(EBP)$.
The difference between ${\mathfrak F}_*$ and this last rational function is
\[
 \frac{t^3\,T(t)}
 {36(21t+10)(9t^2-2t+9)(9t^2+14t+9)(40t^2-9t+81)},
\]
where
\begin{align*}
 T(t)={}&1331640t^6+102141t^5+4649295t^4-1707854t^3\\
       &+2188518t^2-498591t+32643.
\end{align*}
On $[1/4,1/3]$,
\[
 T(t)\ge\frac{4857700}{3}t^2-498591t+32643
 \ge\frac{27592}{3}>0.
\]
The last quadratic is increasing for $t\ge1/4$.
This proves \eqref{eq:fixed-rational-F}.

\paragraph{A uniform density component for the removed part.}
The density $\varrho$ of $Q$ satisfies $\varrho(q)\ge1/\pi$ on $(-1,1)$.
Indeed, its two branches are
\[
 \varrho(q)=\frac{2}{\pi^2\sqrt{8/9}}
 K\!\left(\frac{1-q^2}{8/9}\right)
 \quad(|q|>1/3),
\]
\[
 \varrho(q)=\frac{2}{\pi^2\sqrt{1-q^2}}
 K\!\left(\frac{8/9}{1-q^2}\right)
 \quad(|q|<1/3),
\]
where $K(z)=\int_0^{\pi/2}(1-z\sin^2\theta)^{-1/2}\,d\theta$.
Both bounds follow directly from $K(z)\ge\pi/2$.

Use the slightly larger cutoff $c_0=449/1000>c$.
For example, $\sqrt3>433/250$ and $\pi<22/7$ imply this comparison.
It suffices to prove $\mathbb E\max\{f(Q),c_0\}<1$, where
$f(q)=u/(1-\beta q)^2$, because the maximum is nondecreasing in the cutoff.
Put
\[
 d_s=\left(\frac{1-s}{1+s}\right)^2,
 \qquad \ell=\min\{u,(c_0+d_s)/2\}.
\]
The symmetric pair satisfies
$\min\{f(q),c_0\}+\min\{f(-q),c_0\}\ge2\ell$.
Consequently
\[
 \mathbb E\min\{f(Q),c_0\}
 \ge\frac I\pi+(1-2/\pi)\ell,
\]
where direct integration gives
\[
 I=\int_{-1}^1\min\{f(q),c_0\}\,dq
 =1+c_0-\frac{(1-\sqrt{c_0})^2+(1+\sqrt{c_0})^2s^2}{2s}.
\]
The lower bound $\sqrt{c_0}>67/100$ gives
\[
 I>I_0:=\frac{1449}{1000}-\frac{109}{2000s}-\frac{2789s}{2000}.
\]
Using $7/22<1/\pi<106/333$ yields two rational lower bounds
\begin{align*}
 H_1(s)&=\frac{24061}{22000}+\frac2{11}d_s
                         -\frac{763}{44000s}-\frac{19523s}{44000},\\
 H_2(s)&=\frac{4453}{4400}+\frac4{11}u
                 -\frac{763}{44000s}-\frac{19523s}{44000}-\frac1{3663},
\end{align*}
such that
\[
 1-c_0+\mathbb E\min\{f(Q),c_0\}\ge\min\{H_1(s),H_2(s)\}.
\]
For $H_1$, note that $I_0-c_0-d_s$ is positive: it is concave in $s$
and its endpoint values are $119/1440$ and
$223219211/7302960000$.
For $H_2$, use $I_0-2u>-2$ and
$0<1/\pi-7/22<1/7326$; the final constant allows for its possible
negative sign.

\paragraph{Two elementary derivative comparisons.}
We show ${\mathfrak F}_*<H_1,H_2$ on $[1/2,23/40]$.
First,
\[
 H_1''\ge\frac{857496919}{14440214250}>0,
 \qquad H_1'(1/2)=-\frac{572717}{1188000}>-\frac{49}{100}.
\]
Thus $H_1'>-49/100$. To control ${\mathfrak F}_*'$, put $a=27/40$ and $b=21/10$.
Polynomial division of $as^6/(1+bs^2)$ gives
\[
 {\mathfrak F}_*'''(s)=24s\left[
 \frac29+\frac ab-\frac ab\frac{1-bs^2}{(1+bs^2)^4}\right]>0.
\]
Therefore ${\mathfrak F}_*'$ is convex, and its maximum is attained at an endpoint.
The two rational endpoint derivatives satisfy
\[
 {\mathfrak F}_*'(1/2)=-\frac{129157}{267912}<-\frac{23}{50},
 \qquad
 {\mathfrak F}_*'(23/40)=-\frac{198087307472543}{423301179456000}
                  <-\frac{23}{50}.
\]
It follows that $(H_1-{\mathfrak F}_*)'>-3/100$.
Since
\[
 H_1(1/2)-{\mathfrak F}_*(1/2)=\frac{38789}{12078000}>\frac3{1000},
\]
we obtain $H_1-{\mathfrak F}_*>3/4000>0$ throughout the interval.

For $H_2$, the bound $|u'|\ge5000/4913$ gives
\[
 H_2'\le-\frac{160922023}{216172000}<-\frac7{10}.
\]
On the other hand, differentiating ${\mathfrak F}_*$ with respect to $t=s^2$
shows $\partial_t{\mathfrak F}_*\ge-2/3+4t/9\ge-5/9$, hence
${\mathfrak F}_*'\ge-23/36$. Thus $H_2-{\mathfrak F}_*$ is strictly decreasing.
At its right endpoint, the following rational comparisons suffice:
\[
 {\mathfrak F}_*(23/40)=\frac{2044346455979}{2498365440000}
                    <\frac{8183}{10000},
\]
\[
 H_2(23/40)=\frac{50010180257224249}{61099275457440000}
                    >\frac{8185}{10000}.
\]
This proves the fixed-amplitude estimate for $1/2\le s\le23/40$.

\paragraph{Extension of the logarithmic density bound through the interior singularity.}
We finish $23/40\le s\le3/5$ by monotonicity, using the bound for the
fixed two-cosine integral at $s=3/5$ established in
Section~\ref{sec:fixed-high-monotonicity}.
Retain $\delta=1/3$, $\kappa=9/8$, $\varrho_*=3/(2\pi\sqrt2)$,
$C=106/333$, and $d_1=1/4-C$. The density bound~\eqref{eq:fixed-RN} with $N=1$ extends to
all $0\le q<1$, $q\ne\delta$, in the form
\begin{equation}\label{eq:fixed-central-log}
 \varrho(q)\le \varrho_*\widetilde{\mathcal R}_1(q),\qquad
 \widetilde{\mathcal R}_1(q)=1-C\log|\kappa(q^2-\delta^2)|
                         +d_1\kappa(1-q^2).
\end{equation}
It remains to establish the bound on the branch $q<\delta$. Set $z=\kappa(1-q^2)\in[1,9/8]$.
The central density divided by $\varrho_*$ is
$z^{-1/2}\sum_{n\ge0}a_nz^{-n}$.
Applying the same coefficient bound $a_n<C/n$ to this series, its
upper bound minus $\widetilde{\mathcal R}_1(q)$ is at most
\[
 D(z)=z^{-1/2}-1+Cz^{-1/2}\log z+d_1(z^{-3/2}-z).
\]
Here the discarded term is nonpositive because
$\log(z-1)<0$ and $z^{-1/2}\le1$.
Now $D(1)=0$ and
\[
 z^{5/2}D'(z)
 \le(C-1/2)z+(C-1/4)(z^{5/2}+3/2)
 \le(4C-5/4)z-\tfrac12(C-1/4)
 \le-\frac{85}{10656}<0.
\]
We used $z^{5/2}\le3z-2$ on $[1,9/8]$, which follows by differentiating
and observing $(5/2)(9/8)^{3/2}<3$.
This proves \eqref{eq:fixed-central-log}.

Differentiation under the expectation is justified by dominated
convergence on each compact subinterval of \(0<\beta<1\):
the kernel and its first parameter derivative are uniformly bounded,
and the level set \(u/(1-\beta Q)^2=c\) has probability zero.
This argument remains valid when the truncation height crosses
the logarithmic singularity of the density.

Write $I(\beta)=\mathbb E(f(Q)-c)_+$ with the cutoff $c$.
As in Section~\ref{sec:fixed-high-monotonicity},
\[
 I'(\beta)=2\int_{q_0}^1\frac{q-\beta}{(1-\beta q)^3}\varrho(q)\,dq,
 \qquad q_0=q_0(s,c)=\frac{1-\sqrt{(1-\beta^2)/c}}\beta.
\]
The expression for $q_0$ follows from \eqref{eq:truncation-height}.
For $\beta_a\le\beta\le\beta_b$ and $q_0\ge l>0$, define
\[
 \widetilde{\mathcal N}(l,b)
 =2\int_l^b\frac{b-q}{(1-bq)^3}\widetilde{\mathcal R}_1(q)\,dq.
\]
By \eqref{eq:fixed-density}, the positive part $q\ge\beta$ is at least
$\varrho_*/(1-\beta^2)$; the negative kernel is increasing in $\beta$ for
$q\le\beta$. Hence
\[
 I'(\beta)\ge \varrho_*\left[
       \frac1{1-\beta_a^2}-\widetilde{\mathcal N}(l,\beta_b)\right].
\]
The following two rows make the bracket strictly positive:
\[
\begin{array}{c|c|c|c|c|c}
 s\text{-interval}&\beta_a&\beta_b&l&\widetilde{\mathcal N}<
   &1/(1-\beta_a^2)>\\\hline
 {[23/40,47/80]}&1840/2129&7520/8609&36/125&377/100&79/20\\
 {[47/80,3/5]}&7520/8609&15/17&5/16&199/50&21/5
\end{array}
\]
At the respective left endpoints,
\[
 \frac{1-\beta_a^2}{(1-\beta_a l)^2}
 =\frac{14630625}{32615521},\quad
   \frac{145161}{323761},
 \qquad\text{both}<\frac{4183}{9324}<c.
\]
Since $q_0$ increases in $s$, the lower cutoffs hold on both rows.
Thus $I$ is strictly increasing on $[23/40,3/5]$.
The bound at $s=3/5$ established in
Section~\ref{sec:fixed-high-monotonicity} completes the claim.

\paragraph{Explicit constants in the two logarithmic comparisons.}
The primitive formulas for $\mathcal N$ in
Section~\ref{sec:fixed-log-evaluations} apply with $\log|\mathscr V(q)|$ and
absolute values in $\mathcal J_\delta$.
If an interval crosses $q=\delta$, apply integration by parts after
removing $(\delta-\varepsilon,\delta+\varepsilon)$ and let
$\varepsilon\downarrow0$. The logarithmic boundary terms cancel, with
error $O(\varepsilon|\log\varepsilon|)$; the remaining $\mathcal J_\delta$ is its
ordinary principal-value logarithmic primitive. The integral
$\widetilde{\mathcal N}$ is an absolutely convergent improper integral.

For completeness, an integer $n$ below means
$n/10^6<\log v<(n+1)/10^6$.
Every row follows from the positive finite atanh-series bounds in
Section~\ref{sec:fixed-log-evaluations}, using twenty-one terms.
\[
\begin{array}{c|r|r}
\text{row}&v&n\\\hline
1&3961/125000&-3451818\\
 &434838719/592919048&-310083\\
 &1386746329/439112025&1149960\\
 &11236205155/298596177&3627804\\
 &25103668445/4092524073&1813852\\\hline
2&31/2048&-4190632\\
 &217/289&-286530\\
 &3349/1024&1184945\\
 &1379/16&4456525\\
 &197/32&1817467
\end{array}
\]
In each row these five quantities are $|\mathscr V(l)|$, $|\mathscr V(b)|$,
$(1-bl)/(1-b^2)$,
$|(b-\delta)(1-bl)/((l-\delta)(1-b^2))|$,
and $(b+\delta)(1-bl)/((l+\delta)(1-b^2))$.
Substitution in the displayed elementary primitive gives respectively
$\widetilde{\mathcal N}<377/100$ and $\widetilde{\mathcal N}<199/50$.
This completes the proof of Theorem~\ref{thm:fixed-product-mid}.
\end{proof}

\section{Large conformal centres}\label{app:large-centres}

We establish the radial cap criteria and the product estimates needed
for $3/5\le s<1$. The product direction is divided into three
amplitude ranges after the reduction in Lemma~\ref{lem:product-direction}.

\subsection{Radial cap inequalities}
\label{sec:oscillator-caps}

The radial density in Lemma~\ref{lem:radial-cap-density} can be retained
more accurately than in its rational bound $h/(q-e)$ after normalization
by $4\pi$. Three finite radial majorants yield direct-cap criteria. They apply only once the truncation height
lies above the critical-height ceiling; that geometric condition is
stated explicitly throughout.

Fix $0<s<1$, put $h=1-R$, and suppose that $m=\max q>\sqrt R$. Define
\[
 e=(2R-1)m+2\sqrt{R(1-R)}\sqrt{1-m^2},\qquad
 M=\frac{m-e}{2}>0.
\]
Choose $0<c\le L/A(\weights)$ and write
\[
 P=1-s^2,\quad B=1+s^2,\quad
 q_0=q_0(s,c)=\frac{B-P/\sqrt c}{2s},\quad {\mathfrak a}=hc.
\]
The expression for $q_0$ is \eqref{eq:truncation-height} in this notation.
The case $m\le q_0$ has $J=0$. In the remaining direct-cap case
$e\le q_0<m$, set
\[
 D=B-2se,\quad b=B-2sm,\quad y=\frac bD,\quad
 z=\frac{P}{\sqrt c\,D}.
\]
Then $0<y<z\le1$. Identity~\eqref{radial:conformal-square} gives
\begin{equation}\label{eq:oscillator-square}
 y\ge {\mathfrak a}z^2.
\end{equation}
Lemma~\ref{lem:radial-cap-density} supplies the following bound for the
height density $\rho(q)=-\frac{d}{dq}\Area\{p\cdot F_{\weights}>q\}$,
with area measured using $g_{\weights}$
\begin{equation}\label{eq:oscillator-density}
 \rho(q)\le\frac{2\pi h}{M}\frac{\theta}{\sin\theta},
 \qquad
 \cos\theta=\frac{q-(m+e)/2}{M},\quad e<q<m.
\end{equation}

\begin{theorem}[Direct-cap criteria]\label{thm:oscillator-direct-caps}
Under the preceding geometric hypotheses, assume ${\mathfrak a}\le2/3$. Each of the
following alternatives implies $J_s(\weights,p)<4\pi$:
\[
 \begin{array}{c|c}
 \text{lower bound for }{\mathfrak a}&\text{range of }z\\ \hline
 9/40&0<z\le1\\
 1/5&0<z\le1\\
 9/50&0<z\le1\\
 7/50&0<z\le19/20.
 \end{array}
\]
The third row includes the parameter ranges of the first two rows.
The first and fourth criteria suffice for the applications below.
The proofs of the four criteria use finite radial majorants and share
several auxiliary identities.
\end{theorem}

\begin{proof}[Proof of Theorem~\ref{thm:oscillator-direct-caps}]
We apply finite polynomial majorants to the density bound
\eqref{eq:oscillator-density} and verify the four stated parameter ranges.
\subsubsection{Finite oscillator inequalities and monotonicity}
\label{sec:finite-oscillator-weights}

Let
\[
 \mathcal B(w)=\arccos(1-w)\sqrt{\frac{2-w}{w}},\qquad 0<w<2.
\]
It solves
\[
 w(2-w)\mathcal B'(w)+\mathcal B(w)=2-w.
\]
Set $c_1=1/3$ and $c_n=(n-1)c_{n-1}/(2n+1)$ for $n\ge2$, and
\[
 \mathcal B_N(w)=2-\sum_{n=1}^N c_nw^n.
\]
Direct substitution gives
\[
 w(2-w)\mathcal B_N'(w)+\mathcal B_N(w)-(2-w)
   =Nc_Nw^{N+1}>0.
\]
The difference $\mathcal B_N-\mathcal B$ vanishes at zero. Multiplication
of its first-order equation by the positive integrating factor
$\sqrt{w/(2-w)}$ expresses that difference as a positive integral.
Consequently $\mathcal B_N>\mathcal B$ on $(0,2)$. This is a finite
polynomial calculation and does not require a convergence argument at
$w=2$.

Use the cap coordinate $x=(q-e)/(m-e)$. Here
$w=1-\cos\theta=2(1-x)$. Equation~\eqref{eq:oscillator-density} and the
strict oscillator inequalities give
\[
 \frac{J}{4\pi}<U_N({\mathfrak a};y,z)
 :={\mathfrak a}\int_{(1-z)/(1-y)}^1
 \left(\frac{z^2}{[1-(1-y)x]^2}-1\right)w_N(x)\dd x,
\]
where $w_N(x)=\mathcal B_N(2(1-x))/(2x)$. The first three weights are
\begin{align*}
 w_1(x)&=\frac2{3x}+\frac13,\\
 w_2(x)&=\frac8{15x}+\frac7{15}+\frac2{15}(1-x),\\
 w_3(x)&=\frac{16}{35x}+\frac{19}{35}
          +\frac{22}{105}(1-x)+\frac8{105}(1-x)^2.
\end{align*}
They are positive. Increasing $y$ decreases the integrand and increases
the lower endpoint; hence $U_N$ decreases in $y$. By
\eqref{eq:oscillator-square}, it is enough to bound $U_N({\mathfrak a};{\mathfrak a}z^2,z)$.

For a positive weight $w(x)\ge a_0/x+(1-a_0)$, let
\[
 \mathcal F_w({\mathfrak a},z)=
 \int_{(1-z)/(1-{\mathfrak a}z^2)}^1
 \left(\frac{z^2}{[1-(1-{\mathfrak a}z^2)x]^2}-1\right)w(x)\dd x-\frac1{\mathfrak a}.
\]
The moving-endpoint term vanishes when differentiating in ${\mathfrak a}$. Thus
\[
 \partial_{\mathfrak a}\mathcal F_w=
 \frac1{{\mathfrak a}^2}-2z^4\int_{(1-z)/(1-{\mathfrak a}z^2)}^1
 \frac{xw(x)}{[1-(1-{\mathfrak a}z^2)x]^3}\dd x.
\]
The weights $1/x$ and $1$ give, respectively,
\begin{align*}
 F_{\mathfrak a}(z)&=-z+z^2\log\frac1{{\mathfrak a}z}
 +(1-z^2)\log\frac{1-z}{1-{\mathfrak a}z^2},\\
 T_{\mathfrak a}(z)&=\frac{-2z+({\mathfrak a}+1)z^2}{1-{\mathfrak a}z^2}.
\end{align*}
The scalar formulas below are understood by continuous extension at
$z=0$ and $z=1$ whenever a logarithmic term has a removable endpoint
singularity. In particular, $F_{\mathfrak a}(0)=T_{\mathfrak a}(0)=0$,
$F_{\mathfrak a}(1)=\log(1/{\mathfrak a})-1$, and $T_{\mathfrak a}(1)=-1$.
Consequently
\[
 \partial_{\mathfrak a}\mathcal F_w
 \le a_0\partial_{\mathfrak a}F_{\mathfrak a}+(1-a_0)\partial_{\mathfrak a}T_{\mathfrak a}.
\]
After multiplication by the positive factor
${\mathfrak a}(1-{\mathfrak a}z^2)^2/z^2$, the right-hand side becomes
\[
 a_0({\mathfrak a}-1)(1-{\mathfrak a}z^2)+(1-a_0){\mathfrak a}(1-z)^2.
\]
This is convex in $z$ and has endpoint values ${\mathfrak a}-a_0$ and
$-a_0(1-{\mathfrak a})^2$. Therefore
\begin{equation}\label{eq:oscillator-A-monotone}
 \partial_{\mathfrak a}\mathcal F_w\le0\qquad(0<{\mathfrak a}\le a_0<1).
\end{equation}
For $w_1,w_2,w_3$, the corresponding values of $a_0$ are
$2/3,8/15,16/35$.

\subsubsection{The first and second corrections}

Write $\mathcal F_N=\mathcal F_{w_N}$, so that
$U_N({\mathfrak a};{\mathfrak a}z^2,z)=1+{\mathfrak a}\mathcal F_N({\mathfrak a},z)$. Exact integration gives
\[
 \mathcal F_1=\frac23F_{\mathfrak a}+\frac13T_{\mathfrak a},\qquad
 \mathcal F_2=\frac8{15}F_{\mathfrak a}+\frac7{15}T_{\mathfrak a}+\frac2{15}\mathcal J_1,
\]
where
\[
 \mathcal J_1({\mathfrak a},z)=\frac{z^2[-\log({\mathfrak a}z)-3/2+2{\mathfrak a}z-({\mathfrak a}z)^2/2]}
                   {(1-{\mathfrak a}z^2)^2}.
\]
We first establish the convexity needed below. For \(0<z<1\), put
\[
 g_{\mathfrak a}(z)=\log\frac{1-{\mathfrak a}z^2}{{\mathfrak a}z(1-z)}.
\]
Direct differentiation gives
\[
 F_{\mathfrak a}''(z)-2g_{\mathfrak a}(z)+4
 =\frac{2N_{\mathfrak a}^*(z)}{(1-z)(1-{\mathfrak a}z^2)^2},
\]
where
\[
 N_{\mathfrak a}^*(z)=(1-{\mathfrak a}z^2)^2-(1-z)
 \bigl[(1-{\mathfrak a})+{\mathfrak a}(3-{\mathfrak a})z^2-2{\mathfrak a}^2z^4\bigr].
\]
At \({\mathfrak a}=1/4\), the inequalities \(z(1-z)\le1/4\) and
\(1-z^2/4\ge3/4\) imply \(g_{1/4}(z)\ge\log12>2\).
For the last comparison, the elementary exponential series gives
\(\exp(1)<3\), hence \(\exp(2)<9<12\).
Moreover,
\[
 16N_{1/4}^*(z)
 =4+z(11z^2-19z+12)+z^4(3-2z)>0.
\]
Indeed the quadratic in parentheses has positive leading coefficient
and discriminant \(-167\). Therefore \(F_{1/4}''(z)>0\).
The identity
\[
 \partial_{\mathfrak a} F_{\mathfrak a}''(z)
 =\frac{2({\mathfrak a}-1)(1+3{\mathfrak a}z^2)}{{\mathfrak a}(1-{\mathfrak a}z^2)^3}<0
\]
then proves \(F_{\mathfrak a}''(z)>0\) for every \(0<{\mathfrak a}\le1/4\).

Also
\[
 T_{\mathfrak a}''=\frac{2N_{\mathfrak a}(z)}{(1-{\mathfrak a}z^2)^3},\quad
 N_{\mathfrak a}={\mathfrak a}+1-6{\mathfrak a}z+3{\mathfrak a}({\mathfrak a}+1)z^2-2{\mathfrak a}^2z^3.
\]
Here $N_{\mathfrak a}'=-6{\mathfrak a}(1-z)(1-{\mathfrak a}z)\le0$ and $N_{\mathfrak a}(1)=(1-{\mathfrak a})^2>0$.
Thus $T_{\mathfrak a}''>0$.

At ${\mathfrak a}=9/40$, the function $\mathcal F_1$ is strictly convex, vanishes
at zero, and satisfies
\[
 \mathcal F_1(9/40,1)=\frac23\log\frac{40}{9}-1<0.
\]
The last inequality follows from
$\sum_{k=0}^5(3/2)^k/k!=5711/1280>40/9$.
Convexity and~\eqref{eq:oscillator-A-monotone} prove the first row of
Theorem~\ref{thm:oscillator-direct-caps}.

For the second row, it remains to check $\mathcal J_1''>0$ at ${\mathfrak a}=1/5$.
Put $\ell=-\log({\mathfrak a}z)$ and $D_0=1-{\mathfrak a}z^2$. Direct differentiation gives
\begin{align}
 D_0^4\mathcal J_1''={}&2(1+8{\mathfrak a}z^2+3{\mathfrak a}^2z^4)\ell-6+12{\mathfrak a}z
                -(26{\mathfrak a}+6{\mathfrak a}^2)z^2\notag\\
 &+32{\mathfrak a}^2z^3-(4{\mathfrak a}^2+6{\mathfrak a}^3)z^4+4{\mathfrak a}^3z^5.
 \label{eq:oscillator-J1dd}
\end{align}
The elementary logarithm estimates
\begin{equation}\label{eq:oscillator-log5}
 \frac85<\log5<\frac{13}{8}
\end{equation}
follow, respectively, from
\[
 2\sum_{j=0}^3\frac{(2/3)^{2j+1}}{2j+1}
 =\frac{122492}{76545}=\frac85+\frac4{15309}
\]
and
$\sum_{k=0}^5(13/8)^k/k!=19839493/3932160>5$.
Hence at ${\mathfrak a}=1/5$,
\[
 \ell>\frac85+(1-z)+\frac{(1-z)^2}{2}+\frac{(1-z)^3}{3}.
\]
Substitution into~\eqref{eq:oscillator-J1dd} gives a polynomial $P_{20}$.
Writing $b_{j,n}(z)=\binom njz^j(1-z)^{n-j}$, its full positive identity is
\begin{align*}
 39375P_{20}={}&34125b_{0,7}+13875b_{1,7}+9650b_{2,7}
 +11340b_{3,7}\\
 &+14928b_{4,7}+19200b_{5,7}+24480b_{6,7}+30240b_{7,7}.
\end{align*}
Thus $\mathcal F_2(1/5,\cdot)$ is strictly convex. Its endpoint values are
\[
 \mathcal F_2(1/5,0)=0,\qquad
 \mathcal F_2(1/5,1)=\frac{89\log5-148}{120}<0,
\]
using~\eqref{eq:oscillator-log5}. This proves the second row for
$1/5\le {\mathfrak a}\le8/15$; the first row handles the larger values up to $2/3$.

\subsubsection{The third correction}

The last weight has the exact scalar expression
\begin{equation}\label{eq:oscillator-F3}
 \mathcal F_3=\frac{16}{35}F_{\mathfrak a}+\frac{19}{35}T_{\mathfrak a}
                  +\frac{22}{105}\mathcal J_1+\frac8{105}\mathcal J_2,
\end{equation}
where
\[
 \mathcal J_2({\mathfrak a},z)=
 \frac{\frac23z^3+{\mathfrak a}z^4(1+2\log({\mathfrak a}z))-2{\mathfrak a}^2z^5+\frac13{\mathfrak a}^3z^6}
 {(1-{\mathfrak a}z^2)^3}.
\]
Equivalently, for $k=1,2$,
\[
 \mathcal J_k({\mathfrak a},z)=\int_{(1-z)/(1-{\mathfrak a}z^2)}^1
 \left(\frac{z^2}{[1-(1-{\mathfrak a}z^2)x]^2}-1\right)(1-x)^k\dd x.
\]
We prove strict convexity in $z$ at ${\mathfrak a}=9/50$ and ${\mathfrak a}=7/50$.
The only additional derivative identity required is
\begin{align}
 D_0^5\mathcal J_2''={}&4z+26{\mathfrak a}z^2+(20{\mathfrak a}-40{\mathfrak a}^2)z^3
                 +(26{\mathfrak a}^2+10{\mathfrak a}^3)z^4\notag\\
 &+(8{\mathfrak a}^2-52{\mathfrak a}^3)z^5+(6{\mathfrak a}^4-4{\mathfrak a}^3)z^6-4{\mathfrak a}^4z^7\notag\\
 &-12{\mathfrak a}z^2(2+5{\mathfrak a}z^2+{\mathfrak a}^2z^4)\ell.
 \label{eq:oscillator-J2dd}
\end{align}

For ${\mathfrak a}\le1/5$, the inequality $\log u\ge1-1/u$ and
\eqref{eq:oscillator-log5} give $\log(1/{\mathfrak a})>13/5-5{\mathfrak a}$.
Thus in~\eqref{eq:oscillator-J1dd} one may substitute
\[
 \ell>\begin{cases}
 17/10+(1-z)+(1-z)^2/2,&{\mathfrak a}=9/50,\\
 19/10+(1-z)+(1-z)^2/2,&{\mathfrak a}=7/50.
 \end{cases}
\]
The resulting polynomials $P_{18}$ and $P_{14}$ have the complete
positive expansions
\begin{align*}
 937500P_{18}={}&375000b_{0,6}+87500b_{1,6}+133850b_{2,6}
                   +292650b_{3,6}\\
              &+461093b_{4,6}+597860b_{5,6}+757140b_{6,6},\\
 937500P_{14}={}&750000b_{0,6}+387500b_{1,6}+328650b_{2,6}
                   +392850b_{3,6}\\
              &+488561b_{4,6}+578270b_{5,6}+698580b_{6,6}.
\end{align*}
This proves $\mathcal J_1''>0$ at both values of ${\mathfrak a}$.

For $0<z\le1$, elementary differentiation gives
\[
 -\log z\le\frac{1-z^2}{2z}.
\]
We also have
\[
 \log\frac{50}{9}<\frac{26}{15}<\frac74,
 \qquad \log\frac{50}{7}<2.
\]
Indeed,
\[
 \sum_{k=0}^5\frac{(26/15)^k}{k!}-\frac{50}{9}
   =\frac{620197}{11390625}>0,
 \qquad
 \sum_{k=0}^5\frac{2^k}{k!}=\frac{109}{15}>\frac{50}{7}.
\]
The coefficient of $\ell$ in~\eqref{eq:oscillator-J2dd} is negative.
Substitute $\ell<7/4+(1-z^2)/(2z)$ at ${\mathfrak a}=9/50$, and
$\ell<2+(1-z^2)/(2z)$ at ${\mathfrak a}=7/50$. The lower bounds are
$D_0^5\mathcal J_2''>zQ_{18}(z)$ and $D_0^5\mathcal J_2''>zQ_{14}(z)$, where, with
$t=1-z$,
\begin{align*}
 1562500Q_{18}={}&1148368+532521t+189630t^2+455220zt^2\\
                &+956367t^4+70713zt^4+48114t^6,\\
 3125000Q_{14}={}&3309651+779247t+1796060t^2+252840zt^2\\
                &+1318394t^4+46991zt^4+46648t^6.
\end{align*}
Every displayed coefficient is positive. Thus $\mathcal J_2''>0$, and
\eqref{eq:oscillator-F3} is strictly convex at both specified values of
${\mathfrak a}$.

At ${\mathfrak a}=9/50$ the endpoint calculation is
\begin{align*}
 \mathcal F_3(9/50,1)
 &=\frac{5203208}{7236705}\log\frac{50}{9}
                         -\frac{662936}{529515}\\
 &<-\frac{618472}{108550575}<0.
\end{align*}
Since $\mathcal F_3(9/50,0)=0$, convexity proves the third row of
Theorem~\ref{thm:oscillator-direct-caps}, first up to ${\mathfrak a}=16/35$ by
\eqref{eq:oscillator-A-monotone}, and then up to $2/3$ by the first row.

For the fourth row, at ${\mathfrak a}=7/50$ and $z=19/20$ set
\[
 \ell_0=\log\frac{1000}{133},\qquad
 n_0=\log\frac{1000}{17473}.
\]
Exact substitution in~\eqref{eq:oscillator-F3} gives
\[
 \mathcal F_3(7/50,19/20)=a\ell_0+\frac{39}{875}n_0+b,
\]
where
\[
 a=\frac{422940225430991}{666825875352125}<\frac{127}{200},
 \qquad
 b=-\frac{11170787340922}{9243130945475}<-\frac{241}{200}.
\]
The positive atanh series at $1/3$ gives
\[
 \frac23<\log2
 <2\left[\frac13+\frac{1/27}{3(1-1/9)}\right]=\frac{25}{36}.
\]
Since $1000/133<8$ and $17473/1000>16$, it follows that
$\ell_0<25/12$ and $n_0<-8/3$. Therefore
\[
 \mathcal F_3(7/50,19/20)
 <\frac{127}{200}\frac{25}{12}
   -\frac{39}{875}\frac83-\frac{241}{200}
 =-\frac{79}{84000}<0.
\]
Strict convexity between $0$ and $19/20$, followed by monotonicity in
${\mathfrak a}$ and the same overlapping range as above, proves the final row.
This completes the proof of Theorem~\ref{thm:oscillator-direct-caps}.
\end{proof}

\subsubsection{Consequences and geometric scope}

For the envelope cutoff $c=c_R=L/\Abar(R)$,
\[
 {\mathfrak a}(R)=2\left(\frac1{\sqrt3}-\frac1\pi\right)
          \sqrt{\frac{1-R}{1+3R}}
\]
is decreasing. The bounds $\sqrt3<97/56$ and $\pi>333/106$ give
\[
 2\left(\frac1{\sqrt3}-\frac1\pi\right)>
 \frac{16732}{32301}.
\]
For $R\le R_*$, the bounds in Section~\ref{sec:global-majorants} give ${\mathfrak a}(R)>7/50$,
so the fourth criterion in Theorem~\ref{thm:oscillator-direct-caps}
applies whenever $z\le19/20$.

The hypotheses for applying these criteria are
$q_0\ge e_R(m)$ and, for the fourth row when needed,
$\chi_s(e_R(m))\le(19/20)\sqrt c$. These will be verified in the
particular radius and weight ranges in which the criteria are used.

\subsection{\texorpdfstring{Large centres with \(R\le1/2\)}{Large centres with R <=1/2}}
\label{sec:central-small-weights}

\begin{theorem}\label{thm:central-small-weights}
For every fundamental embedding and direction, if
\[
 \frac13\le R\le\frac12,\qquad \frac35\le s<1,
\]
then $J_s(\weights,p)<4\pi$.
\end{theorem}

\begin{proof}
Put
\[
 c_0=\frac49,\quad a=\frac{15}{17},\quad {t_{\mathrm h}}=\frac{13}{20},\quad
 b_0=\frac13,\quad H_* =\frac{943}{1000},
\]
\[
 f(q)=\chi_{3/5}(q)^2=\frac{64}{(17-15q)^2},\quad
 C=\frac{f({t_{\mathrm h}})-c_0}{({t_{\mathrm h}}+b_0)^2}
   =\frac{2340800}{2927521},\quad K=\frac{127}{200}.
\]
We use the envelope cutoff ${c_R}=L/\Abar(R)>4/9$.
For a nonempty cap with $q_0\ge e_R(m)$, the first criterion in
Theorem~\ref{thm:oscillator-direct-caps} applies throughout $R\le1/2$,
since the bounds in Section~\ref{sec:global-majorants} give
$9/40<(1-R){c_R}<2/3$.

\paragraph{A uniform envelope.}
For $q\le a$ and $s\ge3/5$, we have $\chi_s(q)^2\le f(q)$. For $q>a$,
maximizing in $s$ gives $\chi_s(q)^2\le1/(1-q^2)$, and
\[
 \frac1{1-q^2}-f(q)
 =\frac{(q-a)^2}{(1-q^2)(1-aq)^2}.
\]
The denominator decreases on $[a,H_*]$, and the exact comparison
\[
 324(1-H_*^2)(1-aH_*)^2-1
   =\frac{34858710071}{2890000000000}>0
\]
therefore proves
\begin{equation}\label{eq:central-envelope}
 \chi_s(q)^2\le f(q)+324(q-a)_+^2\qquad(q\le H_*,\ s\ge3/5).
\end{equation}
The derivative argument in Lemma~\ref{lem:quadratic-cap-majorant} applies here: indeed $(1+(3/5)^2)/(1-(3/5)^2)=17/8>\sqrt{8/3}$. Hence $(f-c_0)/(q+b_0)^2$ increases
on its positive support. Thus ${Q_{\mathrm{bulk}}}(q)=C(q+b_0)^2$ majorizes
$(f-c_0)_+$ for $q\le{t_{\mathrm h}}$, while
\[
 \psi(q)=f(q)-c_0-{Q_{\mathrm{bulk}}}(q)\ge0\qquad(q\ge{t_{\mathrm h}}).
\]
Consequently ${Q_{\mathrm{bulk}}}$ plus the cap excess $\psi+324(q-a)_+^2$ is a
pointwise majorant of the actual truncated integrand on $q\le H_*$.

\paragraph{Quadratic base and density.}
The moment contribution is at most
\[
 C\frac{\Abar(R)}{4\pi}\left(\frac R2+\frac19\right)<KC.
\]
Indeed the function before $C$ increases on $[1/3,1/2]$: its derivative
has the sign of $-12R^2+(14/3)R+22/9>0$. At $R=1/2$ its value is
$13\pi\sqrt5/144$. Using $\pi<22/7$ and $\sqrt5<161/72$ gives
\[
 \frac{13\pi\sqrt5}{144}<\frac{23023}{36288}<\frac{127}{200},
 \qquad \frac{127}{200}-\frac{23023}{36288}=\frac{71}{129600}.
\]
For fixed $m>\sqrt R$ and $q\le m$ with $q>e_R(m)$, direct differentiation
gives
\[
 \partial_R\frac{1-R}{q-e_R(m)}
 =\frac{m+\sqrt{(1-R)/R}\sqrt{1-m^2}-q}{[q-e_R(m)]^2}\ge0.
\]
The same statement holds with $q=m$. Combining
\eqref{eq:oscillator-density} with the weight $w_1$ in
Section~\ref{sec:finite-oscillator-weights} gives
\[
 \frac{\rho(q)}{4\pi}
 <\frac{2(1-R)}{3(q-e_R(m))}
    +\frac{1-R}{3(m-e_R(m))}
 \qquad(e_R(m)<q<m).
\]
In the applications below, $m>39/50>1/\sqrt2$. The two terms therefore
increase with $R$ up to $1/2$ whenever $q>\sqrt{1-m^2}$.
At $R=1/2$, put $e=\sqrt{1-m^2}$. For $R\le1/2$ and $e<q<m$ we obtain
\begin{equation}\label{eq:central-density}
 \frac{\rho(q)}{4\pi}<\frac1{3(q-e)}+\frac1{6(m-e)}.
\end{equation}

\paragraph{The extreme maximum ranges.}
If $m\le39/50$, Lemma~\ref{radial:bulk-bound} and $\Abar/(4\pi)<11/6$ give
\[
 \frac{J}{4\pi}
 <\frac{11}{6}\left(\frac{6400}{2809}-\frac49\right)
       \frac{1/2}{2(39/50)^2+1/2}
 =\frac{79688125}{81379539}<1.
\]
Here $\chi_s(39/50)$ decreases for $s\ge3/5$, because $39/50<15/17$.
If $m\ge\sqrt{8/9}$, then
\[
 e_R(m)\le\sqrt{1-m^2}\le\frac13,
 \qquad
 q_0(s,4/9)=\frac{5s}{4}-\frac1{4s}\ge\frac13.
\]
By \eqref{eq:truncation-height}, since ${c_R}>4/9$, the actual
truncation height is strictly larger than $q_0(s,4/9)$. The first criterion in
Theorem~\ref{thm:oscillator-direct-caps} closes the direct cap, or
$J=0$ if its positive support is empty.

It remains to treat $39/50<m<\sqrt{8/9}<H_*$. Here $e_R(m)<{t_{\mathrm h}}$, so
\eqref{eq:central-density} applies above ${t_{\mathrm h}}$.

\paragraph{Four explicit integral bounds.}
The following rational constants bound $e=\sqrt{1-m^2}$ from above
and $m-e$ from below. Each ceiling comparison follows by squaring;
each $d$ is the left endpoint minus the stated ceiling.
\begin{center}
\renewcommand{\arraystretch}{1.5}
\begin{tabular}{c|ccc}
$m$ interval & $H$ & ceiling for $e$ & lower bound $d$\\ \hline
$[39/50,17/20]$&$17/20$&$313/500$&$77/500$\\
$[17/20,9/10]$&$9/10$&$527/1000$&$323/1000$\\
$[9/10,93/100]$&$93/100$&$109/250$&$58/125$\\
$[93/100,943/1000]$&$943/1000$&$46/125$&$281/500$
\end{tabular}
\end{center}
In each row use the kernel $1/[3(q-e)]+1/(6d)$ and extend the integral
to $H$. The correction in~\eqref{eq:central-envelope} satisfies
\[
 \int_a^H324(q-a)^2\left(\frac1{3(q-e)}+\frac1{6d}\right)\dd q
 \le\left(\frac{36}{a-e}+\frac{18}{d}\right)(H-a)_+^3.
\]
Put
\[
 k=17-15e,\quad D_H=\frac1{17-15H}-\frac1{17-15{t_{\mathrm h}}},
\]
\[
 A_e=\frac{64}{3k^2}-\frac{c_0+C(e+b_0)^2}{3},\qquad
 B_e=\frac{64}{3k^2}.
\]
Elementary partial fractions bound the normalized cap contribution plus
its quadratic base by
\begin{equation}\label{eq:central-E}
 E=E_0+A_e\log\frac{H-e}{{t_{\mathrm h}}-e}
          +B_e\log\frac{17-15{t_{\mathrm h}}}{17-15H},
\end{equation}
where every term of $E_0$ is rational:
\begin{align*}
 E_0={}&KC+\left(\frac{64}{3k}+\frac{64}{90d}\right)D_H\\
 &-\frac C3\left[\frac{H^2-{t_{\mathrm h}}^2}{2}+(e+2b_0)(H-{t_{\mathrm h}})\right]\\
 &-\frac{c_0(H-{t_{\mathrm h}})+(C/3)[(H+b_0)^3-({t_{\mathrm h}}+b_0)^3]}{6d}\\
 &+\left(\frac{36}{a-e}+\frac{18}{d}\right)(H-a)_+^3.
\end{align*}
For example, the only non-polynomial primitive needed is
\[
 \int\frac{64\dd q}{(17-15q)^2(q-e)}
 =\frac{64}{k^2}\log\frac{q-e}{17-15q}
      +\frac{64}{k(17-15q)}.
\]
For all four displayed $e$ values, $A_e<0$ and $B_e>0$.
The following eight logarithm bounds, substituted in
\eqref{eq:central-E}, yield the stated rational upper bounds:
\begin{center}
\small
\renewcommand{\arraystretch}{1.5}
\begin{tabular}{c|cc|cc|c}
row & first argument & log lower bound & second argument & log upper bound & $E<$\\ \hline
1&$28/3$&$2233/1000$&$29/17$&$107/200$&$4813/5000$\\
2&$373/123$&$1109/1000$&$29/14$&$729/1000$&$4779/5000$\\
3&$247/107$&$209/250$&$145/61$&$433/500$&$4937/5000$\\
4&$575/282$&$89/125$&$1450/571$&$233/250$&$9903/10000$
\end{tabular}
\end{center}
Every entry in the final column is less than one.
For explicit verification, write each log argument as
$x=2^j(1+v)/(1-v)$. In the order of the table, the eight pairs $(j,v)$ are
\[
 (3,1/13),\ (0,6/23),\ (1,127/619),\ (1,1/57),
\]
\[
 (1,33/461),\ (1,23/267),\ (1,11/1139),\ (1,77/648).
\]
Define the finite rational functions
\[
 \mathcal L(v)=2\sum_{i=0}^4\frac{v^{2i+1}}{2i+1},\qquad
 \mathcal U(v)=\mathcal L(v)+\frac{2v^{11}}{11(1-v^2)}.
\]
The positive atanh expansion and its geometric tail give
\[
 j\mathcal L(1/3)+\mathcal L(v)<\log x
 <j\mathcal U(1/3)+\mathcal U(v).
\]
Substitution of the eight rational pairs proves the printed log bounds;
substitution in~\eqref{eq:central-E} proves the final column by rational
arithmetic. The bulk range, these four intervals and the direct
extreme cover every maximum $m$, proving the theorem.
\end{proof}

\subsection{\texorpdfstring{The cap range \(s\ge(1+R)/2\)}{The cap range s >=(1+R)/2}}
\label{sec:global-highcap-linear}

\begin{theorem}\label{thm:global-highcap-linear}
For every fundamental embedding and every direction $p$, if
\[
 \frac9{20}\le R\le R_*\qquad\text{and}\qquad
 \frac{1+R}{2}\le s<1,
\]
then $J_s(\weights,p)<4\pi$.
\end{theorem}

\begin{proof}
Write ${c_R}=L/\Abar(R)$ and use the bounds of Section~\ref{sec:global-majorants}
${c_R}>4/9$ and ${\mathfrak a}=(1-R){c_R}>7/50$. Also ${\mathfrak a}<2/3$, since ${\mathfrak a}(R)$
decreases and ${\mathfrak a}(1/3)=(2/3)(1-\sqrt3/\pi)<2/3$.
Choose
\[
 s_0=\frac{1+R}{2},\qquad M=\frac{73}{100}+\frac R4.
\]
Since $R_*<19/25$, we have $s_0<22/25$, $M<23/25$ and
\[
 M^2-R\ge\left(\frac{23}{25}\right)^2-\frac{19}{25}
       =\frac{54}{625}>0.
\]
Here $M^2-R$ decreases on the relevant interval. Moreover, for $s\ge s_0$,
\[
 \frac{2s}{1+s^2}\ge\frac{2(29/40)}{1+(29/40)^2}
 =\frac{2320}{2441}>\frac{23}{25}>M.
\]
Thus $\chi_s(M)$ decreases with $s\ge s_0$.

\paragraph{\texorpdfstring{Bulk maxima $m\le M$.}{Bulk maxima m <= M .}}
The identity $2M-s_0=24/25$ gives
\[
 \chi_{s_0}(M)=\frac{1-s_0^2}{1-24s_0/25}\le\frac{25}{16},
\]
since the last inequality is equivalent to $(4s_0-3)^2\ge0$.
Lemma~\ref{radial:bulk-bound} therefore gives
\[
 \frac{J}{4\pi}\le\phi(R)\left(\frac{625}{256}-\frac49\right),
 \qquad
 \phi(R)=\frac{\Abar(R)}{4\pi}\frac{R}{2M^2+R}.
\]
We claim $\phi<1/2$. Since $\pi<22/7$, it suffices to prove
\[
 T(R)=49(2M^2+R)^2-484R^2(1-R)(1+3R)>0.
\]
With $x=R-3/5$, the exact expansion is
\begin{align*}
 T={}&\frac{92977}{64}x^4+\frac{253983}{100}x^3
       +\frac{2773571}{2500}x^2\\
    &+\frac{377858}{15625}x+\frac{12148801}{390625}.
\end{align*}
For $x\ge0$ all coefficients are positive. If $-3/20\le x\le0$, dropping
the positive even powers gives
\[
 T>31-25\frac3{20}-2540\left(\frac3{20}\right)^3
   =\frac{7471}{400}>0.
\]
It follows that
\[
 \frac{J}{4\pi}<\frac12\left(\frac{625}{256}-\frac49\right)
              =\frac{4601}{4608}<1.
\]

\paragraph{\texorpdfstring{Cap maxima $m\ge M$.}{Cap maxima m >= M .}}
For $0<\nu<1$, let
\[
 q_\nu(s)=\frac{\nu+1}{2\nu}s-\frac{1-\nu}{2\nu s},
 \qquad \chi_s(q_\nu(s))=\nu.
\]
This function increases with $s$. Since $e_R(m)$ decreases for
$m>\sqrt R$, it is enough to prove
\[
 q_\nu(s_0)>e_R(M),\qquad
 \nu\le\frac{19}{20}\sqrt{{c_R}}.
\]
These inequalities imply $q_0>e_R(m)$ and
$z=\chi_s(e_R(m))/\sqrt{{c_R}}<19/20$. If the positive cap is nonempty,
the fourth criterion in Theorem~\ref{thm:oscillator-direct-caps}
applies; otherwise $J=0$.

For $9/20\le R\le3/5$, take $\nu=3/5$. Because ${c_R}>4/9$,
\[
 \nu=\frac35<\frac{19}{30}<\frac{19}{20}\sqrt{{c_R}}.
\]
Set $W=q_\nu(s_0)-(2R-1)M$. Exact simplification gives
\[
 W=\frac{-150R^3-313R^2+256R+219}{300(R+1)}>0.
\]
For the sign, use
$q_\nu(s_0)\ge q_\nu(29/40)=147/290>1/2$, while
$(2R-1)M\le22/125<1/2$. The squared difference is
\[
 W^2-4R(1-R)(1-M^2)=\frac{S_1(R)}{90000(R+1)^2}.
\]
With $x=3/5-R\in[0,3/20]$, its numerator is
\begin{align*}
 S_1={}&60000x^5-99575x^4-73120x^3+108784x^2\\
       &+\frac{258048}{5}x+\frac{46656}{25}.
\end{align*}
The two negative terms are absorbed into the quadratic term, since
\[
 108784-73120\frac3{20}-99575\left(\frac3{20}\right)^2
 =\frac{1529209}{16}>0.
\]
Thus $S_1>0$. Taking positive square roots proves
$q_\nu(s_0)>e_R(M)$.

For $3/5\le R\le R_*$, take $\nu=33/50$. The cutoff ${c_R}$ increases
for $R\ge1/3$. Using $\sqrt3<97/56$, $\pi>333/106$ and
$\sqrt{28/25}<53/50$, we obtain
\[
 {c_R}\ge {c_{3/5}}>
 \frac{100}{53}\left(\frac{56}{97}-\frac{106}{333}\right)
 =\frac{836600}{1711953}>\left(\frac{66}{95}\right)^2.
\]
The last difference is $93047732/15450375825>0$.
Hence $\nu<(19/20)\sqrt{{c_R}}$.
In this case
\[
 W=q_\nu(s_0)-(2R-1)M
   =\frac{-825R^3-1784R^2+1283R+1392}{1650(R+1)}>0.
\]
Indeed $q_\nu(s_0)\ge q_\nu(4/5)=301/440>2/3$, while
$(2R-1)M<(13/25)(23/25)=299/625<2/3$.
Now
\[
 W^2-4R(1-R)(1-M^2)=\frac{S_2(R)}{2722500(R+1)^2}.
\]
With $x=19/25-R\in[0,4/25]$, the numerator is
\begin{align*}
 S_2={}&1711875x^5-3646925x^4-1668788x^3
                  +\frac{112411588}{25}x^2\\
      &+\frac{699971712}{625}x+\frac{70972452}{15625}.
\end{align*}
Again the negative terms are absorbed, because
\[
 \frac{112411588}{25}-1668788\frac4{25}
       -3646925\left(\frac4{25}\right)^2
 =\frac{103402404}{25}>0.
\]
Thus $S_2>0$, and positive square roots give the required ceiling
comparison. Both cap ranges are therefore covered, completing the proof.
\end{proof}

\subsection{A sixth-moment estimate for the removed part}\label{sec:sixth-moment-removed}

A sixth-moment lower bound for the part removed by truncation gives
the following product estimate.

\begin{theorem}[The sixth-moment strip]\label{thm:sixth-moment-product}
For
\[
 \frac12\le R\le R_*,\qquad
 \frac35\le s\le\frac{13}{20},\qquad0\le\alpha\le1,
\]
one has $\mathcal P(R,s,\alpha)<4\pi$.
\end{theorem}

\begin{proof}
We use the full-average estimate~\eqref{eq:low-radius-alpha-free},
\[
 {\mathfrak F}:=\mathbb E f_s(Q)\le C_R\mathscr H(s),\qquad
 C_R=\frac{1+e}{\sqrt{4R(1-R)}},\quad
 e=\frac{(2R-1)^2}{64R(1-R)},\quad
 \mathscr H(s)=\frac{1+s^4/2}{1+s^2}.
\]
Here
\[
 Q=\sqrt R\,\alpha\cos U+\sqrt{1-R}\sqrt{1-\alpha^2}\cos V,
 \qquad M_{2j}=\mathbb E Q^{2j},\qquad \eta=1-R,
\]
and $U,V$ are independent uniform angles. Put
\[
 u=\left(\frac{1-s^2}{1+s^2}\right)^2,\quad {\omega}=1-u,\quad
 d=\left(\frac{1-s}{1+s}\right)^2,\quad
 f_s(q)=\frac{u}{(1-\sqrt {\omega}\,q)^2}.
\]
The area and cutoff are $\Abar=2\pi^2S$, $c=L/\Abar$, where
\[
 S=\sqrt{(1-R)(1+3R)},\qquad L=4\pi^2/\sqrt3-4\pi.
\]

\paragraph{Two moment inequalities.}
For $1/2\le R\le19/25$, one has
\begin{equation}\label{eq:sixth-moment-ratio}
 M_4\le\frac35M_2.
\end{equation}
To verify this directly, write $x=\alpha^2$. The polynomial
$(3/5)M_2-M_4$ has coefficients
\[
 \frac{3(6R(1-R)-1)}8x^2
 +\frac{3(15R^2-16R+3)}{20}x
 +\frac{3(1-R)(5R-1)}{40}.
\]
Its leading coefficient is positive. Writing $a=2R-1$, its
quadratic discriminant is
\[
 \frac9{6400}(75a^4-6a^2-5)<0.
\]
Indeed, $0\le a^2\le169/625$; the convex quadratic
$75y^2-6y-5$ has negative values at both endpoints of this interval,
namely $-5$ and $-17792/15625$. This proves
\eqref{eq:sixth-moment-ratio} even for unrestricted real $x$.

We also have
\begin{equation}\label{eq:sixth-moment-h2}
 2M_2-M_4\ge
 \begin{cases}
 9/16-13R/32,&1/2\le R\le2/3,\\
 \eta-3\eta^2/8,&2/3\le R\le19/25.
 \end{cases}
\end{equation}
In fact, the exact polynomial $2M_2-M_4$ is
\[
 \eta-\frac38\eta^2
 +\frac{9R^2-4R-1}{4}x
 +\frac{3(6R(1-R)-1)}8x^2.
\]
For $1/2\le R\le5/8$, subtracting $9/16-13R/32$ leaves a
quadratic with positive leading coefficient and discriminant
\[
 \frac{2R-1}{64}(54R^3+36R^2-31R-10)\le0.
\]
The cubic in parentheses increases on this interval and has value
$-545/256$ at $5/8$. For $5/8\le R\le2/3$, the coefficient of $x$
is at least $1/256$, and the constant after that subtraction is
\[
 -\frac{(3R-2)(4R+1)}{32}\ge0.
\]
All coefficients are therefore nonnegative in this second interval.
For $R\ge2/3$, the displayed exact polynomial directly gives the
second line of \eqref{eq:sixth-moment-h2}, because its two nonconstant
coefficients are positive.

By Cauchy--Schwarz, $M_6M_2\ge M_4^2$. Consequently
\[
 \mathbb E[Q^2(1-Q^2)^2]
 =M_2-2M_4+M_6
 \ge M_2\left(1-\frac{M_4}{M_2}\right)^2
 \ge\frac4{25}M_2\ge\frac2{25}\eta.
\]
Combining this with \eqref{eq:sixth-moment-h2}, and defining
\[
 h(z)=\frac52z-2z^2+\frac12z^3
     =(2z-z^2)+\frac12z(1-z)^2,
\]
we obtain $\mathbb E h(Q^2)\ge \mathfrak m(R)$, where
\begin{equation}\label{eq:sixth-moment-m}
 \mathfrak m(R)=
 \begin{cases}
 \displaystyle\frac9{16}-\frac{13R}{32}+\frac{1-R}{25},
       &1/2\le R\le2/3,\\[2mm]
 \displaystyle\frac{26}{25}(1-R)-\frac38(1-R)^2,
       &2/3\le R\le19/25.
 \end{cases}
\end{equation}
The function $\mathfrak m$ is positive and decreasing. It is continuous at
$2/3$, where its value is $61/200$, and
$\mathfrak m(1/2)=607/1600<2/5$.

\paragraph{The paired minimum.}
Let $B=c+d-2u$. On the present physical strip,
\[
 \frac49<c<\frac{59}{100},\qquad
 \frac{53361}{323761}\le u\le\frac{64}{289},\qquad
 \frac{49}{1089}\le d\le\frac1{16}.
\]
The upper bound on \(c\) follows directly from \(R\le19/25\):
\[
 S(R)^2\ge\frac{492}{625}>
 \left(\frac{22}{25}\right)^2,\qquad
 {c_R}=\frac{2w_*}{S(R)}
 <\frac{5181}{8800}<\frac{59}{100}.
\]
Here \(w_*<25905/100000\) is the elementary constant bound established
in Section~\ref{sec:analytic-highcaps}. This estimate is also valid at the auxiliary
endpoint \(R=189/250\).
The displayed inequalities imply
\[
 B>\frac49+\frac{49}{1089}-\frac{128}{289}
   =\frac{14645}{314721}>0.
\]
Put $u_0=53361/323761$. Since $u(1-u)$ increases for $u<1/2$,
\[
\begin{aligned}
 6u{\omega}-\frac52 B
 &>6u_0(1-u_0)-\frac52\left(\frac{59}{100}
                          +\frac1{16}-2u_0\right)\\
 &=\frac{314301084219}{16771389619360}>0.
\end{aligned}
\]
For $0\le z\le1$, the polynomial $h$ satisfies
\[
 0\le h(z)\le1,\qquad h(z)\le\frac52z,
 \qquad1-h(z)=(1-z)^2(1-z/2).
\]
For $q\ge0$, $f_s(-q)<c$. Before $f_s(q)$ reaches $c$, its paired
sum satisfies
\[
 f_s(q)+f_s(-q)\ge2u+6u{\omega}q^2\ge2u+B h(q^2).
\]
After it reaches $c$, the sum of the two minima is at least
$c+d=2u+B\ge2u+B h(q^2)$. Symmetry of $Q$ therefore gives
\begin{equation}\label{eq:sixth-moment-min}
 \mathbb E\min\{f_s(Q),c\}\ge u+\frac B2\mathfrak m(R).
\end{equation}
It follows that $\mathcal P\le\mathcal U(R,s)$, where
\begin{equation}\label{eq:sixth-moment-upper}
 \mathcal U=\Abar C_R\mathscr H(s)-\Abar u
             -\frac{\mathfrak m(R)}2(L+\Abar d-2\Abar u).
\end{equation}

\paragraph{\texorpdfstring{Monotonicity in $s$.}{Monotonicity in s .}}
We use the bound $1\le C_R<6/5$ and the derivative formulas from
Section~\ref{sec:product-half-three-fifths}. On
$3/5\le s\le13/20$, those formulas give
\[
 |\mathscr H'(s)|\le\frac{9347}{23120}<\frac{41}{100},\qquad
 |u'(s)|\ge\frac{177408000}{184220009}>\frac{19}{20},\qquad
 d'(s)<0.
\]
For example, the first estimate uses $s\le13/20$ in the factor $s$,
$s\ge3/5$ in $2-2s^2-s^4$ and in $(1+s^2)^2$;
the second uses $s\ge3/5$, $1-s^2\ge231/400$, and
$1+s^2\le569/400$.
Since $0<\mathfrak m<2/5$, differentiation of \eqref{eq:sixth-moment-upper} yields
\[
 \frac{\partial_s\mathcal U}{\Abar}
 =C_R\mathscr H'-(1-\mathfrak m)u'-\frac{\mathfrak m}2d'
 >-\frac65\frac{41}{100}+\frac35\frac{19}{20}
 =\frac{39}{500}>0.
\]
It is enough to evaluate the upper bound at $s=13/20$.

\paragraph{\texorpdfstring{Convexity in $R$.}{Convexity in R .}}
Let $G=\sqrt{3+1/R}$, and put $t_0=u-d/2>0$ and
$W=2u-2t_0\mathfrak m(R)$. Then
\[
 \mathcal U=\pi^2[G(1+e)\mathscr H(s)-SW]-\frac L2\mathfrak m(R).
\]
The elementary convexity proof for $G(1+e)$ in
Section~\ref{sec:product-weight-convexity} holds on $[1/2,19/25]$. Also
\[
 S'<0,\qquad S''<-3/S,\qquad S^2\le4/3.
\]
Because $\mathfrak m'<0$ and $\mathfrak m<2/5$, we have $W'>0$ and $W\ge u$.
On the first interval in \eqref{eq:sixth-moment-m}, $W''=0$, so
$(SW)''<0$. On its second interval, $\mathfrak m''=-3/4$ and
$W''=3t_0/2$, whence
\[
 (SW)''<\frac{-3W+(3/2)t_0S^2}{S}
 \le\frac{-3u+2t_0}{S}<0.
\]
Finally $-Lm/2$ is convex on each interval. Thus $\mathcal U$
is convex on each side of $2/3$.
At the joining point,
\[
 \mathfrak m'(2/3+)-\mathfrak m'(2/3-)=-\frac{11}{32}.
\]
Since $\mathfrak m$ is continuous, the derivative jump of the complete upper
function is
\[
 \partial_R\mathcal U(2/3+,s)-\partial_R\mathcal U(2/3-,s)
 =\frac{11\Abar B}{64}>0.
\]
This proves convexity on the full interval $[1/2,189/250]$.
This auxiliary upper endpoint is legitimate since
\[
 \frac{189}{250}-R_*=\frac{29}{640250}>0.
\]
Only the upper function is evaluated beyond $R_*$.

\paragraph{Two rational endpoint comparisons.}
At $s=13/20$,
\[
 u=\frac{53361}{323761},\qquad d=\frac{49}{1089},\qquad
 \mathscr H(s)=\frac{348561}{455200}.
\]
For $R=1/2$, let $\mathfrak m_0=607/1600$. Since $G=\sqrt5$, $e=0$,
and $S=\sqrt5/2$, define the positive rational number
\[
 K=\mathscr H(s)-(1-\mathfrak m_0)u-\frac{\mathfrak m_0d}{2}
  =\frac{738894406787}{1128242332800}.
\]
Using $\sqrt5<161/72$, $1/\sqrt3>56/97$, and $\pi<22/7$ gives
\[
\begin{aligned}
 \frac{\mathcal U(1/2,13/20)}\pi
 &<\frac{22}7\left[\frac{161}{72}K
                     -2\mathfrak m_0\frac{56}{97}\right]+2\mathfrak m_0\\
 &=\frac{1427144584604149}{358165656921600}<4,
\end{aligned}
\]
with gap $5518043082251/358165656921600$.

For $R_1=189/250$, one has
\[
 \mathfrak m_1=\frac{115717}{500000},\qquad
 W_1=2(1-\mathfrak m_1)u+\mathfrak m_1d
     =\frac{46497237335227}{176287864500000},
\]
\[
 G(1+e)=\frac{11785}{11529}\sqrt{\frac{817}{189}}
        <\frac{21253}{10000},\qquad
 S=\frac{\sqrt{49837}}{250}>\frac{8929}{10000}.
\]
The radical comparisons follow by squaring, with respective positive
differences
\[
 \frac{113308237057541}{2512147194900000000},\qquad
 \frac{12159}{100000000}.
\]
Consequently
\[
\begin{aligned}
 \frac{\mathcal U(R_1,13/20)}\pi
 &<\frac{22}7\left[
   \frac{21253}{10000}\frac{348561}{455200}
   -\frac{8929}{10000}W_1-2\mathfrak m_1\frac{56}{97}\right]+2\mathfrak m_1\\
 &=\frac{108751546385993117737}{27204422726250000000}<4,
\end{aligned}
\]
with gap $66144519006882263/27204422726250000000$.
In these uses of $\pi<22/7$, the displayed rational brackets are
positive; if the original bracket were nonpositive, the desired
comparison would be immediate.
Convexity in $R$, monotonicity in $s$, and \eqref{eq:sixth-moment-min}
give $\mathcal P<4\pi$ throughout the stated strip.
\end{proof}

\subsection{Direction reduction and a density lower bound}\label{sec:product-direction-density}

\begin{lemma}[Restriction of the product direction]\label{lem:product-direction}
Let $f:[-1,1]\to\mathbb R$ be convex, let $1/2\le R<1$, and let
$U,V$ be independent uniform angles. Then
\[
 \sup_{0\le x\le1}\mathbb E f\bigl(\sqrt{Rx}\cos U+
                    \sqrt{(1-R)(1-x)}\cos V\bigr)
 =\sup_{R\le x\le1}\mathbb E f\bigl(\sqrt{Rx}\cos U+
                    \sqrt{(1-R)(1-x)}\cos V\bigr).
\]
\end{lemma}
\begin{proof}
Put $\Phi(a,b)=\mathbb E f(a\cos U+b\cos V)$ on the diamond
$|a|+|b|\le1$, so every argument of $f$ lies in its stated domain.
It is jointly convex, symmetric, separately even, and consequently
nondecreasing in each nonnegative amplitude.
If $x<1/2$, replace it by $1-x$. For the new amplitudes $a',b'$,
\[
 \max\{a',b'\}\ge\max\{a,b\},\qquad
 a'+b'\ge a+b.
\]
Indeed $a'=\sqrt R\sqrt{1-x}$ dominates both old amplitudes, while
\[
 a'+b'-(a+b)
 =\bigl(\sqrt R-\sqrt{1-R}\bigr)
           \bigl(\sqrt{1-x}-\sqrt x\bigr)\ge0.
\]
For a convex symmetric coordinatewise nondecreasing function, these two
comparisons imply $\Phi(a,b)\le\Phi(a',b')$.
To see this directly, put $A=\max\{a',b'\}$, $B=\min\{a',b'\}$.
Increase the old coordinates, without exceeding $A$, until their sum is
$A+B$. The resulting pair is on the segment joining $(A,B)$ and $(B,A)$;
monotonicity followed by convexity proves the assertion.
Thus only $x\ge1/2$ need be considered.

For $1/2\le x\le R$, the amplitudes satisfy $a\ge b$. Set
$m=a+b$ and $\delta=a-b$. Cauchy--Schwarz gives $m\le1$, and
$\delta$ is increasing in $x$, so $0\le\delta\le2R-1$.
The integer torus map $(\theta,\varphi)\mapsto
(U,V)=(\theta+\varphi,\theta-\varphi)$ sends Haar probability to Haar
probability. Hence
\[
 \Phi(a,b)=\Psi(m,\delta):=
 \mathbb E f(m\cos\theta\cos\varphi-\delta\sin\theta\sin\varphi).
\]
For $|m|,|\delta|\le1$, Cauchy--Schwarz gives
\[
 |m\cos\theta\cos\varphi-\delta\sin\theta\sin\varphi|
 \le\max\{|m|,|\delta|\}
 \bigl(|\cos\theta\cos\varphi|+|\sin\theta\sin\varphi|\bigr)
 \le1.
\]
Thus every argument of $f$ also lies in $[-1,1]$ throughout the square
on which this second comparison is made.
The function $\Psi$ is jointly convex and separately even:
$\theta\mapsto\pi-\theta$ reverses $m$ alone, and
$\theta\mapsto-\theta$ reverses $\delta$ alone.
It is therefore nondecreasing in the two nonnegative variables.
Consequently $\Psi(m,\delta)\le\Psi(1,2R-1)$, which is exactly the
value at $x=R$. This proves the claimed restriction.
\end{proof}

If $R<1/2$, interchange the two cosine terms first, replacing
$(R,x)$ by $(1-R,1-x)$. The same conclusion therefore holds with
$\rho=\max\{R,1-R\}$ after this interchange.

On the restricted branch $R\le x\le1$, one has
\[
 \sqrt R\le m\le1,\qquad 2R-1\le\delta\le\sqrt R,
\]
and, with $d=2R-1$ and ${\varkappa}=4R(1-R)$,
\[
 \delta=dm+\sqrt {\varkappa}\sqrt{1-m^2}.
\]
The latter follows by writing
$\sqrt x=\sqrt R\,m+\sqrt{1-R}\sqrt{1-m^2}$ and
$\sqrt{1-x}=\sqrt{1-R}\,m-\sqrt R\sqrt{1-m^2}$.

\subsubsection{Geometry and the removed part in the inner amplitude region}
\label{sec:inner-geometry-minimum}

\begin{lemma}\label{lem:inner-total-height}
Let $1/2\le R\le R_*$, $13/20\le s<1$, and $R\le x\le1$. Put
\[
 h=1-R,\quad d=2R-1,\quad {\varkappa}=4Rh,\quad
 \beta=\frac{2s}{1+s^2},\quad u=1-\beta^2,\quad {\omega}=\beta^2,
\]
\[
 m=\sqrt{Rx}+\sqrt{h(1-x)},\qquad y={\omega}(x-h)-d.
\]
If $y\le(3/10)\sqrt{{\varkappa}u}$, then $m>\beta$.
\end{lemma}
\begin{proof}
The hypothesis $x\ge R$ is essential. On this branch, $m$ decreases
from $1$ to $\sqrt R$. Moreover
$\beta\ge520/569>9/10$ and $R<R_+ :=19/25<81/100$, so
$\beta>\sqrt R$. The unique point with $m=\beta$ is
\[
 x_\beta=(\sqrt R\,\beta+\sqrt h\sqrt u)^2.
\]
Indeed $\sqrt{1-x_\beta}=\sqrt h\,\beta-\sqrt R\sqrt u>0$,
and direct substitution gives $m=\beta$. At this point,
\[
 \frac{y_\beta}{\sqrt{{\varkappa}u}}
 =(1-u)^{3/2}-\frac d{\sqrt {\varkappa}}\sqrt u(2-u).
\]
The ratio $d/\sqrt {\varkappa}$ increases with $R$, its derivative being
$1/[4(R(1-R))^{3/2}]$. Since $\sqrt{114}>32/3$,
\[
 \frac d{\sqrt {\varkappa}}<\frac{13}{2\sqrt{114}}<\frac{39}{64}.
\]
For $s\ge13/20$, one has $\beta\ge520/569$ and
$\sqrt u\le231/569$. The function $t(2-t^2)$ increases up to this
last endpoint. Therefore
\begin{align*}
 \frac{y_\beta}{\sqrt{{\varkappa}u}}
 &>\left(\frac{520}{569}\right)^3
  -\frac{39}{64}\frac{231}{569}
       \left[2-\left(\frac{231}{569}\right)^2\right]\\
 &=\frac{3646115551}{11790080576}
  =\frac3{10}+\frac{545456891}{58950402880}>\frac3{10}.
\end{align*}
The variable $y$ increases strictly with $x$, so $y\le(3/10)\sqrt{{\varkappa}u}$
implies $x<x_\beta$. Monotonicity of $m$ proves the claim.
\end{proof}

\begin{lemma}\label{lem:removed-part-monotone}
Let $Q=a\cos U+b\cos V$ with $a,b\ge0$ and $m=a+b\le1$.
Suppose $0<c<1$, $0<\beta<1$, $u=1-\beta^2\le4c/9$, and $m\ge\beta$.
For $f(q)=u/(1-\beta q)^2$, one has
\[
 \mathbb E\min\{f(Q),c\}\ge L_0(m)\ge L_0(\beta),
\]
where
\[
 L_0(m)=\frac1{\pi m}\int_{-m}^m\min\{f(q),c\}\dd q
                 +\left(1-\frac2\pi\right)u,
\]
\begin{equation}\label{eq:removed-part-at-beta}
 L_0(\beta)=
 \frac{2\sqrt{cu}-cu-u/(2-u)}{\pi(1-u)}
     +\left(1-\frac2\pi\right)u.
\end{equation}
\end{lemma}
\begin{proof}
The random variable $Q$ has an even density ${\varpi}$. We claim that
\[
 {\varpi}(q)\ge\frac1{\pi m}\qquad\text{for almost every }-m<q<m.
\]
Suppose first that $a,b>0$, and put $\delta=|a-b|$ and $\gamma=4ab$.
The angles $\theta=U-V$ and $V$ are independent and uniform. Conditional
on $\theta$, the variable $Q$ has the same distribution as
\[
 \sqrt{a^2+b^2+2ab\cos\theta}\,\cos V.
\]
Since $a^2+b^2+2ab\cos\theta=m^2-\gamma\sin^2(\theta/2)$,
averaging the conditional arcsine density and substituting
$z=\sin(\theta/2)$ gives, with $\Delta_q=m^2-q^2$,
\[
 {\varpi}(q)=\frac2{\pi^2}
 \int_0^{\min\{1,\sqrt{\Delta_q/\gamma}\}}
 \frac{\dd z}{\sqrt{1-z^2}\sqrt{\Delta_q-\gamma z^2}}.
\]
For $\delta<|q|<m$, one has $0<\Delta_q<\gamma$. Dropping the
factor $(1-z^2)^{-1/2}\ge1$ gives
\[
 {\varpi}(q)\ge\frac2{\pi^2}
  \int_0^{\sqrt{\Delta_q/\gamma}}
   \frac{\dd z}{\sqrt{\Delta_q-\gamma z^2}}
 =\frac1{\pi\sqrt\gamma}\ge\frac1{\pi m}.
\]
For $|q|<\delta$, one has $\Delta_q>\gamma$, and
$(\Delta_q-\gamma z^2)^{-1/2}\ge\Delta_q^{-1/2}$. Hence
\[
 {\varpi}(q)\ge\frac2{\pi^2\sqrt{\Delta_q}}
  \int_0^1\frac{\dd z}{\sqrt{1-z^2}}
 =\frac1{\pi\sqrt{m^2-q^2}}\ge\frac1{\pi m}.
\]
The exceptional levels $|q|=\delta$ do not affect the integral.
If one amplitude vanishes, then $m>0$ by the hypothesis $m\ge\beta$,
and the arcsine density $1/(\pi\sqrt{m^2-q^2})$ proves the same bound.

Each symmetric pair of truncated values has sum at least $2u$.
If both values lie below $c$, this is convexity of $f$. If both lie
above, use $c\ge9u/4$. In the mixed case, use
\[
 c+\frac{u}{(1+\beta m)^2}\ge c+\frac u4
 \ge\frac{5u}{2}>2u.
\]
Subtract the uniform density $1/(\pi m)$ from ${\varpi}$. The remaining even,
nonnegative density has mass $1-2/\pi$. The pair estimate proves
$\mathbb E\min\{f(Q),c\}\ge L_0(m)$.

Put $r=\sqrt{u/c}\le2/3$. The truncation point
$q_0=q_0(s,c)=(1-r)/\beta$ from \eqref{eq:truncation-height}
lies in $(0,\beta)$: indeed $uc<1$ implies $r>u$.
Thus, for $m\ge\beta$, direct integration yields
\[
 \int_{-m}^m\min\{f(q),c\}\dd q
 =cm+\frac{2\sqrt{cu}-c-u/(1+\beta m)}{\beta}.
\]
Differentiating gives
\[
 L_0'(m)=
 \frac{c-2\sqrt{cu}+u(1+2\beta m)/(1+\beta m)^2}
      {\pi\beta m^2}.
\]
For $t=\beta m\le1$, one has $(1+2t)/(1+t)^2\ge3/4$. Hence
\[
 L_0'(m)\ge
 \frac{c(1-r/2)(1-3r/2)}{\pi\beta m^2}\ge0.
\]
Substitution of $m=\beta$ gives~\eqref{eq:removed-part-at-beta}.
\end{proof}

For the envelope cutoff, $4/9<{c_R}<3/5$. On $s\ge13/20$,
\[
 u\le\left(\frac{231}{569}\right)^2<\frac{16}{81}
       <\frac49{c_R}.
\]
Consequently, once Lemma~\ref{lem:product-direction} has reduced the
amplitude range to $x\ge R$, the inner condition in
Lemma~\ref{lem:inner-total-height}
provides all hypotheses of Lemma~\ref{lem:removed-part-monotone}.

\subsection{The inner product directions}
\label{sec:high-inner}

We use the amplitude reduction to \(x=\alpha^2\ge R\) from
Lemma~\ref{lem:product-direction}, the density lower bound from the
proof of Lemma~\ref{lem:removed-part-monotone}, and the geometric
estimate of Lemma~\ref{lem:inner-total-height}.

\begin{proposition}\label{prop:inner-product}
Suppose
\[
 \frac12\le R\le R_*=\frac{1936}{2561},\qquad
 \frac{13}{20}\le s\le\frac{1+R}{2},\qquad x=\alpha^2\ge R.
\]
Put
\[
 \tau=\frac{1-s^2}{1+s^2},\quad u=\tau^2,\quad {\omega}=1-u,\quad
 d=2R-1,\quad {\varkappa}=4R(1-R),\quad
 y={\omega}(x-1+R)-d.
\]
If \(y\le3\sqrt{{\varkappa}u}/10\), then \(\mathcal P(R,s,\alpha)<4\pi\).
\end{proposition}

\begin{proof}[Proof of Proposition~\ref{prop:inner-product}]
We bound the full conformal integral, retain a lower bound for the removed
part, and then use separate convexity and endpoint comparisons.
\subsubsection{A tangent bound for the complete conformal integral}

Set
\[
A=2\pi^2\sqrt{(1-R)(1+3R)},\qquad c=L/A,
\]
and write
\[
 \lambda=\sqrt{{\omega}Rx},\quad \mu=\sqrt{{\omega}(1-R)(1-x)},\quad
 w=1-\lambda^2-\mu^2,\quad G_0=w^2-4\lambda^2\mu^2.
\]
The elementary conditional circle integration gives
\[
 {\mathfrak F}=\mathbb E \chi_s(Q)^2
   =u\,\mathbb E_\theta(w-2\lambda\mu\cos\theta)^{-3/2}.
\]
Set \(M=u\sqrt w/G_0\), \(T=u/w^{3/2}\), and
\(z=1-G_0/w^2\in[0,1)\). Lemma~\ref{lem:full-average-series}
gives \({\mathfrak F}=M\mathscr C(z)\). Positivity of the coefficients \(\delta_n\)
implies
\[
 \mathscr C(z)\le \mathscr C_4(z)=1-\frac z{16}-\frac{15z^2}{1024}
                     -\frac{105z^3}{16384}
                     -\frac{15015z^4}{4194304}.
\]
This polynomial is concave on \([0,1]\). Its tangent at \(z=3/4\) is
bounded above, coefficient by coefficient when written in \(1-z\), by
\[
 a+b(1-z),\qquad a=\frac{229}{250},\quad b=\frac{203}{2000}.
\]
The two comparisons are obtained simply by substituting \(3/4\) into
\(\mathscr C_4+(1-z)\mathscr C_4'\) and \(-\mathscr C_4'\); both coefficients are positive.
Consequently
\begin{equation}\label{eq:inner-tangent}
 {\mathfrak F}\le aM+bT.
\end{equation}

The exact identities
\[
 G_0=y^2+{\varkappa}u,\qquad w=\frac{{\varkappa}(1+u)}2-dy
\]
yield a convenient bound independent of \(y\):
\begin{equation}\label{eq:inner-M}
 M\le\overline M
 :=\frac{(1+u/2)(1+d^2u/(4{\varkappa}))}{\sqrt{2{\varkappa}}}.
\end{equation}
For completeness, apply the tangent bound for the square root at
\(w_0={\varkappa}(1+u)/2\) and use
\[
 \sup_{y\in\mathbb R}\frac{1-\eta y}{y^2+K}
       =\frac{1+\sqrt{1+\eta^2K}}{2K},
 \quad \eta=\frac d{{\varkappa}(1+u)},\quad K={\varkappa}u.
\]
First this gives
\(\sqrt{(1+u)/(2{\varkappa})}[1+d^2u/(4{\varkappa}(1+u)^2)]\), and then
\eqref{eq:inner-M}.

Since \(R\le19/25\), one has \({\varkappa}\ge456/625>(17/20)^2\).
Under the inner hypothesis,
\[
 w\ge\frac {\varkappa}2(1-X),\qquad
 X=\frac{12}{17}d\tau-\tau^2.
\]
Here \(0\le u\le(231/569)^2<1/6\) and
\(-1/6\le X\le(6d/17)^2<1/25\).
Taylor's theorem gives
\[
 (1-X)^{-3/2}\le1+\frac32X+3X^2.
\]
For example the second derivative, divided by \(2\), is at most
\((15/8)(25/24)^{7/2}<3\) between \(0\) and \(X\).
It follows that
\begin{equation}\label{eq:inner-T}
 T\le\overline T
 :=\frac{2\sqrt2\,u}{{\varkappa}^{3/2}}\left(1+\frac32X+3X^2\right).
\end{equation}

\subsubsection{Retaining the removed part of the truncated integral}

Lemma~\ref{lem:inner-total-height} gives
\(m=\max Q\ge\sqrt {\omega}\). Lemma~\ref{lem:removed-part-monotone}
therefore gives
\[
 \mathbb E\min(\chi_s(Q)^2,c)\ge L_0(\sqrt {\omega})
 =\frac{2\sqrt c\,\tau-cu-u/(2-u)}{\pi(1-u)}
       +\left(1-\frac2\pi\right)u.
\]
The hypotheses of Lemma~\ref{lem:removed-part-monotone} hold here:
\(c>4/9\),
\(u<1/6<4c/9\), and \(c>2u\).
Combining this inequality with \eqref{eq:inner-tangent}--\eqref{eq:inner-T}
proves \(\mathcal P\le\Phi(R,\tau)\), where, with
\({G_R}=\sqrt{3+1/R}\),
\begin{align}
 \Phi={}&\pi^2{G_R}\left[
 \frac a{\sqrt2}(1+u/2)\left(1+\frac{d^2u}{4{\varkappa}}\right)
 +\frac{2\sqrt2\,b\,u}{{\varkappa}}\left(1+\frac32X+3X^2\right)\right]\notag\\
 &-\frac{2\sqrt{AL}\,\tau}{\pi(1-u)}
 +\frac{Lu}{\pi(1-u)}
 +Au\left[\frac1{\pi(1-u)(2-u)}-1+\frac2\pi\right].
 \label{eq:inner-Phi}
\end{align}
We next prove two convexity statements on the full rectangle
\[
 \frac12\le R\le\frac{19}{25},\qquad 0\le\tau\le\frac1{\sqrt6}.
\]

\subsubsection{Convexity in the weight}

The functions \(A(R)\) and \(\sqrt{A(R)}\) are concave.
The coefficient multiplying \(A\) on the second line of
\eqref{eq:inner-Phi} is nonpositive, since
\[
 \frac1{(1-u)(2-u)}+2\le\frac{36}{55}+2<3<\pi.
\]
The coefficient multiplying \(\sqrt A\) is also nonpositive. These
terms are therefore convex in \(R\).

For the first line, put \({\xi_R}=R-1/2\), and let
\[
 D_R=16R^3(1-R)^3(3R+1)^2.
\]
Direct differentiation gives
\[
 ({G_R}d^j/{\varkappa})''=\frac{{G_R} N_j({\xi_R})}{D_R}\quad(j=0,1,2),
\]
where
\begin{align*}
 N_0={}&216{\xi_R}^4+396{\xi_R}^3+173{\xi_R}^2+14{\xi_R}+59/4,\\
 N_1={}&144{\xi_R}^5+288{\xi_R}^4+210{\xi_R}^3+164{\xi_R}^2+(163/2){\xi_R}-5/2,\\
 N_2={}&48{\xi_R}^5+204{\xi_R}^4+360{\xi_R}^3+185{\xi_R}^2+20{\xi_R}+25/2.
\end{align*}
In particular \({G_R}\) and \({G_R}d^2/{\varkappa}\) are convex. Also
\[
 1+\frac32X+3X^2=P_0+P_1d+P_2d^2,
\]
\[
 P_0=1-\frac32u+3u^2\ge\frac56,\quad
 P_1=\frac{18}{17}\tau(1-4u)\in[0,1/2],\quad
 P_2=\frac{432}{289}u\ge0.
\]
Every nonconstant coefficient of \(P_0N_0+P_1N_1+P_2N_2\)
is nonnegative, while its constant coefficient is at least
\((5/6)(59/4)-(1/2)(5/2)>0\).
Thus all of \(\Phi\) is convex in \(R\).

\subsubsection{\texorpdfstring{Convexity in the radius variable \(\tau\)}{Convexity in the radius variable tau}}

All derivatives in this subsection are with respect to \(\tau\), at fixed \(R\).
Let \({\mathfrak f}=a\overline M+b\overline T\), so that \(\Phi=A[{\mathfrak f}-L_0(\sqrt {\omega})]\).
The elementary estimates
\[
 \overline M''\ge1/\sqrt2,\qquad
 \overline T''\ge2\sqrt2(23/16-3u-10u^2)
\]
are proved as follows. The first is immediate from \eqref{eq:inner-M}.
For the second, the bracket obtained by differentiating
\(\tau^2(1+3X/2+3X^2)\) twice is
\[
 Q_0+Q_1d+Q_2d^2,
\]
\[
 Q_0=2-18u+90u^2,\quad
 Q_1=\frac{108}{17}\tau-\frac{1440}{17}\tau^3,\quad
 Q_2=\frac{5184}{289}u.
\]
Completing the square in \(d\) gives
\[
 Q_0-\frac{Q_1^2}{4Q_2}=\frac{23}{16}-3u-10u^2
 \ge\frac{95}{144}>0.
\]
The assertion at \(u=0\) follows by continuity. Since \({\varkappa}\le1\),
the factor \({\varkappa}^{-3/2}\) can be discarded for this lower bound. Hence
\[
 {\mathfrak f}''\ge\frac1{\sqrt2}
   \left(a+\frac{23b}{4}-12bu-40bu^2\right).
\]

Put \(H(\tau)=\tau^2/[(1-\tau^2)(2-\tau^2)]\).
Differentiating \eqref{eq:removed-part-at-beta} with
\(\beta=\sqrt {\omega}\), and completing the square in \(\sqrt c\), yields
\begin{align*}
 \frac{\pi\Phi_{\tau\tau}}A
 &\ge \frac\pi{\sqrt2}
       \left(a+\frac{23b}{4}-12bu-40bu^2\right)
       +H''+4-2\pi
       -\frac{2u(3+u)^2}{(1-u)^3(1+3u)}.
\end{align*}
Indeed the two terms depending on \(c\) before completion of the square are
\([2c(1+3u)-4\sqrt c\,\tau(3+u)]/(1-u)^3\).
The remaining rational terms simplify completely. Since
\(H=1/(1-\tau^2)-2/(2-\tau^2)\), direct differentiation gives
\[
 H''-\frac{2u(3+u)^2}{(1-u)^3(1+3u)}
 =\frac2{1+3u}-\frac{4(2+3u)}{(2-u)^3}.
\]
For \(0\le u\le1/6\), the three elementary endpoint bounds are
\[
 a+\frac{23b}{4}-12bu-40bu^2\ge a+\frac{95b}{36},\qquad
 \frac2{1+3u}\ge\frac43,\qquad
 \frac{4(2+3u)}{(2-u)^3}\le\frac{2160}{1331}.
\]
Use \(\pi>333/106\), \(1/\sqrt2>70/99\), and \(\pi<22/7\).
The preceding lower bound therefore yields the single rational comparison
\begin{align*}
 \frac{\pi\Phi_{\tau\tau}}A
 &>\frac{333}{106}\frac{70}{99}
       \left(a+\frac{95b}{36}\right)
       +4-\frac{44}{7}+\frac43-\frac{2160}{1331}\\
 &=\frac{386984401}{7110734400}>0.
\end{align*}
This proves strict convexity in \(\tau\).

\subsubsection{Twelve elementary endpoint comparisons}

Separate convexity reduces the desired comparison to the corners of five
rectangles. Put
\[
 (R_0,R_1,R_2,R_3,R_4,R_5)
 =\left(\frac12,\frac{11}{20},\frac35,\frac{13}{20},\frac7{10},\frac{19}{25}\right).
\]
On \(R_i\le R\le R_{i+1}\) the desired radius range is contained in
\[
 \frac{13}{20}\le s\le s_i:=\frac{1+R_{i+1}}2.
\]
These become rectangles in \((R,\tau)\). The corners at \(s=13/20\)
can in turn be reduced to \(R=R_0,R_5\), by weight convexity.

Here are fully rational bounds for the remaining comparisons.
Write \(w_* =1/\sqrt3-1/\pi\). The elementary bounds
\[
 \frac{333}{106}<\pi<\frac{355}{113},\quad
 \frac{265}{153}<\sqrt3<\frac{1351}{780}
\]
give \(25903/100000<w_*<25905/100000\).
The following table specifies numbers \(S_-,G_+,N_-\);
all three columns are integers divided by \(10^6\).
\[
 \begin{array}{c|r|r|r}
 R&S_-&G_+&N_-\\ \hline
 1/2&1118033&2236068&761057\\
 11/20&1092016&2195036&752150\\
 3/5&1058300&2160247&740447\\
 13/20&1016120&2130367&725541\\
 7/10&964365&2104418&706823\\
 19/25&887242&2077448&677970
 \end{array}
\]
Their defining comparisons, verified just by squaring, are
\[
 S_-^2<(1-R)(1+3R),\quad G_+^2>3+1/R,\quad
 N_-^2<2S_-\frac{25903}{100000}.
\]
In the following formula \(u=\tau^2\), \(X=12d\tau/17-u\), and
\[
 \begin{split}
 \Upsilon(R,\tau)={}&\frac{355}{113}G_+
 \left[\frac{99a}{140}(1+u/2)\left(1+\frac{d^2u}{4{\varkappa}}\right)
       +\frac{99b}{35}\frac u{\varkappa}(1+3X/2+3X^2)\right]\\
 &-\frac{4N_-\tau}{1-u}
 +\frac{4(25905/100000)u}{1-u}\\
 &+2S_-u\left[\frac1{(1-u)(2-u)}+2-\frac{333}{106}\right].
 \end{split}
\]
Since the last bracket is negative and all the other replacements have the
indicated sign, \(\Phi(R,\tau)/\pi<\Upsilon(R,\tau)\).
Thus every following check is a rational substitution:
\[
 \Upsilon(R_0,231/569)<381/100,\qquad
 \Upsilon(R_5,231/569)<381/100.
\]
For the upper \(s\) endpoints, the strict upper bounds in the last two
columns are integers divided by \(10^4\):
\[
 \begin{array}{c|c|r|r}
 {[R_i,R_{i+1}]}&s_i&
 10^4\,\Upsilon(R_i,\tau(s_i))\ \hbox{is less than}&
 10^4\,\Upsilon(R_{i+1},\tau(s_i))\ \hbox{is less than}\\ \hline
 {[1/2,11/20]}&31/40&39838&39114\\
 {[11/20,3/5]}&4/5&39621&39074\\
 {[3/5,13/20]}&33/40&39574&39158\\
 {[13/20,7/10]}&17/20&39646&39327\\
 {[7/10,19/25]}&22/25&39900&39601
 \end{array}
\]
All are at most \(399/100<4\). Separate convexity proves
\(\Phi<4\pi\) on the five rectangles, hence on the entire region of the
proposition.
\end{proof}

\subsection{The intermediate product directions}
\label{sec:high-middle}

\begin{proposition}[Intermediate product directions]
\label{prop:middle-product}
Let \(R_*=1936/2561\), and suppose
\[
 \frac12\le R\le R_*,\qquad
 \frac{13}{20}\le s\le\frac{1+R}{2},\qquad
 0\le\alpha\le1,\qquad x=\alpha^2\ge R.
\]
Put
\[
 \tau=\frac{1-s^2}{1+s^2},\qquad u=\tau^2,\qquad {\omega}=1-u,
 \qquad {\varkappa}=4R(1-R),
\]
\[
 y={\omega}(x-(1-R))-(2R-1).
\]
If
\[
 \frac3{10}\sqrt{{\varkappa}u}\le y\le\frac12\sqrt{{\varkappa}u},
\]
then \(\mathcal P(R,s,\alpha)<4\pi\).
\end{proposition}

\begin{proof}[Proof of Proposition~\ref{prop:middle-product}]
Set
\[
 \beta=\sqrt {\omega},\qquad d=2R-1,\qquad
 w_0=\frac{{\varkappa}(1+u)}2,\qquad w=w_0-dy,\qquad G=y^2+{\varkappa}u.
\]
Thus $0<\tau\le231/569<41/100$, $u<1/6$, and $R_*<19/25$.
Throughout this proof write
\[
 S(R)=\sqrt{(1-R)(1+3R)},\qquad G_R=\sqrt{3+1/R}.
\]
The symbol $G_R$ is distinct from the denominator $G=y^2+{\varkappa}u$.

\paragraph{A linear bound for the elliptic correction.}
Lemma~\ref{lem:full-average-series} gives
\[
 {\mathfrak F}:=\mathbb E \chi_s(Q)^2=M\,\mathscr C(z),\qquad
 M=\frac{u\sqrt w}{G},\quad z=1-G/w^2,
\]
\[
 \mathscr C(z)\le \mathscr C_4(z)=1-\frac z{16}-\frac{15z^2}{1024}
                    -\frac{105z^3}{16384}-\frac{15015z^4}{4194304}.
\]
The polynomial $\mathscr C_4$ is concave on $[0,1]$. Its tangent at $z=3/4$
therefore yields
\[
 {\mathfrak F}\le aM+bT,\qquad T=\frac{u}{w^{3/2}},
 \qquad a=\frac{229}{250},\quad b=\frac{203}{2000}.
\]
Indeed, the exact tangent coefficients are
$983243445/1073741824<a$ and $6800029/67108864<b$.

On the middle band,
\[
 M\le\frac{1+u/2}{(109/100)\sqrt{2{\varkappa}}}.
\]
Also put $X=(20/17)d\tau-\tau^2$. Since $\sqrt {\varkappa}>17/20$,
\[
 w\ge\frac {\varkappa}2(1-X),\qquad
 -\frac16<X<\frac1{10}.
\]
For the upper bound on $X$, put $a_X=(20/17)d\le52/85$.
The maximum of $a_X\tau-\tau^2$ over real $\tau$ is
$a_X^2/4\le676/7225<1/10$.
Taylor's formula gives
\[
 (1-X)^{-3/2}\le1+\frac32X+3X^2
 \quad(-1/6\le X\le1/10),
\]
because the second derivative $15(1-X)^{-7/2}/4$ is smaller than $6$
on this interval. Consequently
\begin{equation}\label{eq:middle-full-upper}
 {\mathfrak F}\le \frac{a(1+\tau^2/2)}{(109/100)\sqrt{2{\varkappa}}}
 +\frac{2\sqrt2\,b\tau^2}{{\varkappa}^{3/2}}
          \left(1+\frac32X+3X^2\right).
\end{equation}

\paragraph{A polynomial lower bound for the removed part.}
The density estimate established in the proof of
Lemma~\ref{lem:removed-part-monotone} gives a symmetric density at least
$1/(\pi m)$ on $(-m,m)$ for a two-cosine height of amplitude sum $m$.
On the restricted branch supplied by Lemma~\ref{lem:product-direction},
one has $m\ge\sqrt R$. For any cutoff $c_0>2u$, the paired minimum is at least
$2u$, so
\[
 \mathbb E\min\{\chi_s(Q)^2,c_0\}\ge \mathscr L(m,c_0),
\]
where
\[
 \mathscr L(m,c_0):=\frac{c_0}\pi+
 \frac{2\sqrt{c_0}\tau-c_0-u/(1+\beta m)}{\pi\beta m}
 +(1-2/\pi)u.
\]
The displayed formula is a lower bound even if the formal cutoff point
$q_0=q_0(s,c_0)=(1-\tau/\sqrt{c_0})/\beta$ from
\eqref{eq:truncation-height} exceeds $m$: extending the low branch
to $q_0$ adds the nonpositive integral $\int_m^{q_0}(\chi_s(q)^2-c_0)\,\mathrm dq$.
In the present range $q_0>0$, so no issue occurs at the left endpoint.

If $\tau/\sqrt{c_0}\le2/3$, then $\mathscr L(m,c_0)$ is nondecreasing in
$0<m\le1$. To verify this, put $t=\beta m$ and
$z=\tau/\sqrt{c_0}$. Its derivative has the sign of
\[
 z^2(1+2t)-(2z-1)(1+t)^2
 =(1+2t)\left(z-\frac{1+t}{1+2t}\right)(z-(1+t))\ge0.
\]
We may therefore replace $m$ by any fixed smaller positive number.

For $u\le1/6$ and $m\ge7/10$,
\[
 1+\frac u2\le\frac1\beta\le1+\frac35u,
 \qquad
 \frac1{\beta(1+\beta m)}\le\frac1{1+m}+\frac35u.
\]
The first pair follows by squaring or convexity.
For the last inequality, differentiate in $u$:
\[
 \frac{d}{du}\frac1{\beta(1+\beta m)}
 =\frac{1+2\beta m}{2\beta^3(1+\beta m)^2}
 \le\frac{500}{729}\frac{22600}{26569}<\frac35,
\]
using $\beta>9/10$, $\beta m\ge63/100$, and the fact that
$(1+2t)/(1+t)^2$ decreases for $t\ge0$.
Using these inequalities together with
$7/22<1/\pi<106/333$ and a rational $a_0\le\sqrt{c_0}$, we obtain
\begin{equation}\label{eq:middle-L-poly}
 \mathscr L(m,c_0)\ge -B_0+B_1\tau+B_2\tau^2+B_3\tau^3-B_4\tau^4,
\end{equation}
where
\[
 B_0=\frac{106c_0}{333}\left(\frac1m-1\right),\quad
 B_1=\frac{7a_0}{11m},\quad B_3=\frac{7a_0}{22m},
\]
\[
 B_2=1-\frac{106}{333}\left(2+\frac{3c_0}{5m}
                                      +\frac1{m(1+m)}\right),
 \qquad B_4=\frac{106}{333}\frac3{5m}.
\]

\paragraph{Only two weight intervals.}
All constants below are rational. In the table, decimals are exact
terminating decimals.
\[
\begin{array}{c|c|c|c|c|c|c|c|c}
 R_0&R_1&S_+&S_-&m&c_0&a_0&C_0&D_0\\\hline
 1/2&5/8&1.119&1.038&.707&.462&.679&.665&.343\\
 5/8&19/25&1.039&.887&.790&.498&.705&.638&.422
\end{array}
\]
They satisfy
\[
 S(R_1)>S_-,\quad S(R_0)<S_+,\quad m^2<R_0,
 \quad a_0^2<c_0<\frac2{S_+}\left(\frac{56}{97}-\frac{106}{333}\right).
\]
The last comparison implies $c_0<c_R=L/\Abar(R)$ throughout the row.
In particular $c_0>2u$ and $u<4c_0/9$.
Writing ${\varkappa}_1=4R_1(1-R_1)$, the two remaining constants satisfy
\[
 C_0>\frac{a}{109/100}\frac{\sqrt{3+1/R_0}}{2\sqrt2},
 \qquad
 D_0>\frac{\sqrt2\,b\sqrt{3+1/R_0}}{{\varkappa}_1}.
\]
All these comparisons follow by squaring and integer arithmetic.
They preserve the useful identity $S(R)/\sqrt {\varkappa}=G_R/2$.

Put $A_0=(20/17)(2R_1-1)$ and $X_0=A_0\tau-\tau^2$.
The factor $1+3X/2+3X^2$ is increasing for $X\ge-1/6$.
Since $G_R$ decreases with $R$ and ${\varkappa}\ge {\varkappa}_1$, multiplying
\eqref{eq:middle-full-upper} by $S(R)$ gives
\[
 S(R){\mathfrak F}\le C_0(1+\tau^2/2)
    +D_0\tau^2\left(1+\frac32X_0+3X_0^2\right).
\]
Let $\ell_{\mathrm{mid}}(\tau)$ denote the polynomial on the right of
\eqref{eq:middle-L-poly}, using the constants of the appropriate row.
It is positive on the required interval. Indeed
\[
 B_1>\frac12,\quad B_2>-\frac1{20},\quad B_4<\frac3{10},
\]
so $\ell_{\mathrm{mid}}'>1/2-\tau/10-6\tau^3/5>0$ for $\tau\le41/100$.
At the left endpoints listed below, direct rational substitution gives
$\ell_{\mathrm{mid}}>1/50$.
Thus $S(R)\mathscr L(m,c_0)\ge S_-\ell_{\mathrm{mid}}(\tau)$, and
\[
 \frac{\mathcal P}{2\pi^2}\le \Upsilon(\tau)
 :=C_0(1+\tau^2/2)
 +D_0\tau^2\left(1+\frac32X_0+3X_0^2\right)-S_-\ell_{\mathrm{mid}}(\tau).
\]
This is a degree-six polynomial. The required intervals for $\tau$ are
\[
 [87/425,231/569]\quad\hbox{and}\quad[141/1109,231/569],
\]
respectively, since $s\le(1+R_1)/2$.

\paragraph{Convexity and four rational endpoint comparisons.}
Write
\[
 \Upsilon(\tau)=Q_0+Q_1\tau+Q_2\tau^2+Q_3\tau^3+Q_4\tau^4
               -6A_0D_0\tau^5+3D_0\tau^6,
\]
where
\[
\begin{aligned}
 Q_0&=C_0+S_-B_0,& Q_1&=-S_-B_1,\\
 Q_2&=C_0/2+D_0-S_-B_2,& Q_3&=3A_0D_0/2-S_-B_3,\\
 Q_4&=D_0(-3/2+3A_0^2)+S_-B_4.
\end{aligned}
\]
For the first row, direct substitution gives
\[
 Q_2>13/20,\qquad Q_3>-1/5,\qquad Q_4>-1/5,
 \qquad 30D_0A_0^4<2/25.
\]
The minimum of $\tau^3(3\tau-4A_0)$ for $\tau\ge0$ is $-A_0^4$.
Consequently, throughout $0\le\tau\le41/100$,
\[
 \Upsilon''>\frac{13}{10}-\frac65\tau-\frac{12}{5}\tau^2-\frac2{25}>0.
\]
For the second row, $Q_2>13/20$, $Q_3,Q_4>0$, and
$A_0=52/85>41/100$. The function
$\tau^3(3\tau-4A_0)$ is decreasing on this interval, and
\[
 -30D_0(41/100)^3\bigl(3(41/100)-4A_0\bigr)<11/10.
\]
Thus $\Upsilon''>13/10-11/10>0$ in the second row as well.
Convexity reduces the comparison to the following four substitutions:
\[
\begin{array}{c|c|c|c}
 \text{row}&\tau_{\rm left}&\Upsilon(\tau_{\rm left})<
                       &\Upsilon(231/569)<\\\hline
 1&87/425&627/1000&570/1000\\
 2&141/1109&624/1000&590/1000
\end{array}
\]
These are direct rational inequalities in the explicitly given
polynomial and constants.
They prove $\Upsilon(\tau)<627/1000<7/11$ on both intervals. Finally
\[
 \mathcal P<2\pi^2\frac7{11}<4\pi,
\]
using $\pi<22/7$. This proves Proposition~\ref{prop:middle-product}.
\end{proof}

\subsection{The outer product directions}\label{sec:outer-product}

\begin{proposition}[Outer directions]\label{prop:outer-product}
In the notation of Sections~\ref{sec:high-inner} and
\ref{sec:high-middle}, if
\[
 \frac12\le R\le R_*,\quad
 \frac{13}{20}\le s\le\frac{1+R}{2},\quad
 y\ge\frac12\sqrt{{\varkappa}u},
\]
then \(\mathcal P(R,s,\alpha)<4\pi\).
\end{proposition}
\begin{proof}
Lemma~\ref{lem:full-average-series}, together with the identities
\(w={\varkappa}(1+u)/2-dy\), \(G_0=y^2+{\varkappa}u\), and \(d\ge0\), gives
\[
 {\mathfrak F}\le\frac{u\sqrt w}{G_0}
 \le\frac45\sqrt{\frac{1+u}{2{\varkappa}}}.
\]
Here \(c_R>4/9>2u\). Each symmetric pair of the removed minima has sum
at least \(2u\): before either entry is capped this follows from convexity;
after one entry is capped that entry alone is at least \(2u\).
Consequently \(\mathbb E\min\{\chi_s(Q_{R,\alpha})^2,c_R\}\ge u\), and
\[
 \mathcal P\le\Abar(R)
 \left[\frac45\sqrt{\frac{1+u}{2{\varkappa}}}-u\right].
\]
The bracket decreases in \(u\ge0\), since its derivative is
\((2/5)/\sqrt{2{\varkappa}(1+u)}-1<0\) for \({\varkappa}\ge456/625\).
Setting \(u=0\) enlarges the bound, and gives
\[
 \mathcal P\le\frac45\pi^2\sqrt{\frac{3+1/R}{2}}
 \le\frac45\pi^2\sqrt{\frac52}<4\pi.
\]
The last comparison is equivalent to \(\pi^2<10\), which follows from
\(\pi<22/7\).
\end{proof}

\section{Explicit rational arithmetic and polynomial identities}\label{app:explicit-arithmetic}

\subsection{An elementary identity for interval positivity}\label{app:interval-identity}
Let $P(X)=\sum_{k=0}^n p_kX^k$, with real coefficients, and let $a<b$.
For $j=0,\ldots,n$ put
\begin{equation}\label{eq:positive-expansion-rule}
 c_j=\sum_{k=0}^n p_k\sum_{i=0}^k
 \binom{k}{i}\binom{n-k}{j-i}a^{k-i}b^i,
\end{equation}
where a binomial coefficient outside its defining range is zero. Then
\begin{equation}\label{eq:positive-expansion}
 (1+v)^nP\!\left(\frac{a+bv}{1+v}\right)=\sum_{j=0}^n c_jv^j.
\end{equation}
Indeed, multiply each monomial by $(1+v)^n$ and expand
$(a+bv)^k(1+v)^{n-k}$ by the binomial theorem; collecting the
coefficient of $v^j$ gives \eqref{eq:positive-expansion-rule}.
Equivalently, writing $t=(X-a)/(b-a)$, the same polynomial identity is
\[
 P(X)=\sum_{j=0}^n c_j t^j(1-t)^{n-j},\qquad a\le X\le b.
\]
Thus positivity of all $c_j$ proves strict positivity on the entire
closed interval, including both endpoints. In the applications below,
$c_j=\gamma k_j$ with the displayed positive rational scale $\gamma$
and integer coefficients $k_j>0$. Formula~\eqref{eq:positive-expansion-rule}
provides a direct algebraic derivation of each coefficient from the
specified polynomial. If a positive denominator was cleared or a
positive factor was divided out in defining $P$, its sign is recorded
in the corresponding argument. No numerical approximation is required.

For the Bernstein coefficients of $C_8$, the
identity
\[
 x^j=\sum_{i=j}^8\frac{\binom{i}{j}}{\binom{8}{j}}
          \binom{8}{i}x^i(1-x)^{8-i}
\]
follows by the binomial theorem after factoring out $x^j$. Consequently,
if $C_8(x)=\sum_{j=0}^8 c_{8,j}x^j$, its Bernstein coefficients are
$\beta_i=\sum_{j=0}^i c_{8,j}\binom{i}{j}/\binom{8}{j}$.
The basis functions are nonnegative on $[0,1]$ and sum to one, so
$C_8(x)\le\max_i\beta_i$. The nine exact values are listed in Section~\ref{cp:complete-positive-expansions};
in particular each is strictly smaller than $1/2$.

\subsection{Complete positive
expansions}\label{cp:complete-positive-expansions}

The following identities prove the monotonicities and the two inequalities
on the axis stated in Section~\ref{cp:rational-identities-and-monotonicities}. For a rational function \(\phi\) in the first eleven rows
below, put \[
 P_\phi({\sigma})=\frac{{D_{\mathrm{prod}}}({\sigma})({D_{\mathrm{prod}}}\phi)'({\sigma})-{D_{\mathrm{prod}}}'({\sigma})({D_{\mathrm{prod}}}\phi)({\sigma})}{(3-2{\sigma})^2(2-3{\sigma})^2}.
\] This is a polynomial, by the formulas in Section~\ref{cp:rational-identities-and-monotonicities}. Its sign
equals the sign of \(\phi'\). For the final two rows put \[
 P_{\mathrm{axis}k}={D_{\mathrm{prod}}}({\sigma})(551/1000-a_{k,0}({\sigma})).
\] For each row, with nonnegative fractional-linear variable $v$, write $P$ for the indicated polynomial and $n=\deg P$. The identity has the form \[
 (1+v)^n P\!\left(\frac{a+b v}{1+v}\right)
 =\gamma_{\mathrm{cen}}\sum_{j=0}^n k_j v^j,
 \qquad\gamma_{\mathrm{cen}}>0,\quad k_j>0.
\] It proves positivity on \([a,b]\), including the endpoint \(b\) by
the leading coefficient. Every polynomial on the left is defined here or in Section~\ref{cp:rational-identities-and-monotonicities}. Identity~\eqref{eq:positive-expansion-rule} gives the coefficient formula for each row. The positive integer vectors on
the right are given in increasing order of powers of \(v\).

\paragraph{$\widehat h_2$}\label{cp:h2}

Interval \([11/40,7/20]\), degree \(8\).
\[
 \gamma_{\mathrm{cen}}=\frac{968877}{16384}.
\]

\begin{center}\small
\setlength{\tabcolsep}{6pt}\renewcommand{\arraystretch}{1.05}
\begin{tabular}{@{}rr@{}}\toprule $j$ & Coefficient\\\midrule
0 & $52149274383505619873$\\
1 & $408234760291339100056$\\
2 & $1372183957931726886464$\\
3 & $2615017606198316429152$\\
4 & $3112796310523988972960$\\
5 & $2379777225149683472512$\\
6 & $1142834436323126511104$\\
7 & $315020242713274794496$\\
8 & $38090354832347979008$\\
\bottomrule\end{tabular}
\end{center}

\paragraph{$-\widehat h_3$}\label{cp:minus_h3}

Interval \([11/40,7/20]\), degree \(8\).
\[
 \gamma_{\mathrm{cen}}=\frac{8719893}{512}.
\]

\begin{center}\small
\setlength{\tabcolsep}{6pt}\renewcommand{\arraystretch}{1.05}
\begin{tabular}{@{}rr@{}}\toprule $j$ & Coefficient\\\midrule
0 & $40235247056713121$\\
1 & $609090054164038552$\\
2 & $3021255279671709248$\\
3 & $7658178857787490144$\\
4 & $11418683523323311520$\\
5 & $10466237524656183424$\\
6 & $5813966948877075968$\\
7 & $1798572228470620672$\\
8 & $237957125502622976$\\
\bottomrule\end{tabular}
\end{center}

\paragraph{$\widehat h_4$}\label{cp:h4}

Interval \([11/40,7/20]\), degree \(8\).
\[
 \gamma_{\mathrm{cen}}=\frac{78479037}{64}.
\]

\begin{center}\small
\setlength{\tabcolsep}{6pt}\renewcommand{\arraystretch}{1.05}
\begin{tabular}{@{}rr@{}}\toprule $j$ & Coefficient\\\midrule
0 & $2311466417851073$\\
1 & $20656909253855896$\\
2 & $78742151281016384$\\
3 & $168368742138493792$\\
4 & $221638419257862560$\\
5 & $184237355033667712$\\
6 & $94510131275646464$\\
7 & $27362890914136576$\\
8 & $3423845549735168$\\
\bottomrule\end{tabular}
\end{center}

\paragraph{$-\delta$}\label{cp:minus_delta}

Interval \([11/40,7/20]\), degree
\(11\).
\[
 \gamma_{\mathrm{cen}}=\frac{8719893}{2048000}.
\]

\begin{center}\small
\setlength{\tabcolsep}{6pt}\renewcommand{\arraystretch}{1.05}
\begin{tabular}{@{}rr@{}}\toprule $j$ & Coefficient\\\midrule
0 & $153299955561314131332$\\
1 & $1455149431812444898593$\\
2 & $6264900482167196756635$\\
3 & $16152037140633960311445$\\
4 & $27714171407370370566360$\\
5 & $33237795765252214934016$\\
6 & $28438443351447483073344$\\
7 & $17363829325745517025440$\\
8 & $7416668930651411287680$\\
9 & $2111278541877942493440$\\
10 & $360613154784791146752$\\
11 & $28007076832120953088$\\
\bottomrule\end{tabular}
\end{center}

\paragraph{$n_7$}\label{cp:n7}

Interval \([11/40,7/20]\), degree \(13\).
\[
 \gamma_{\mathrm{cen}}=\frac{2906631}{409600000}.
\]

\begin{center}\small
\setlength{\tabcolsep}{6pt}\renewcommand{\arraystretch}{1.05}
\begin{tabular}{@{}rr@{}}\toprule $j$ & Coefficient\\\midrule
0 & $588719762057536380053223$\\
1 & $6619549187402455991637006$\\
2 & $34261244281968028100559784$\\
3 & $108051677086902836057953952$\\
4 & $231591077168218146501379720$\\
5 & $356094645396889479207871224$\\
6 & $403909078763418777740505408$\\
7 & $341986472981471739466359552$\\
8 & $215980950066300606856463616$\\
9 & $100392171305335687069697280$\\
10 & $33343460083310693678588928$\\
11 & $7481981040817450080526336$\\
12 & $1014589035926413510129664$\\
13 & $62622162394457219602432$\\
\bottomrule\end{tabular}
\end{center}

\paragraph{$n_8$}\label{cp:n8}

Interval \([11/40,7/20]\), degree \(13\).
\[
 \gamma_{\mathrm{cen}}=\frac{2906631}{409600000}.
\]

\begin{center}\small
\setlength{\tabcolsep}{6pt}\renewcommand{\arraystretch}{1.05}
\begin{tabular}{@{}rr@{}}\toprule $j$ & Coefficient\\\midrule
0 & $680699735394324858852423$\\
1 & $7676598793163499888391206$\\
2 & $39858343862780068511651584$\\
3 & $126133869608971315291938752$\\
4 & $271360964870701439268910720$\\
5 & $418985550829265629034779824$\\
6 & $477486002537012147845147008$\\
7 & $406473580057807358330045952$\\
8 & $258330612626454563903606016$\\
9 & $120977238742691456326241280$\\
10 & $40543363584825946131376128$\\
11 & $9198288197785237524545536$\\
12 & $1264565420995833341886464$\\
13 & $79426408493729791455232$\\
\bottomrule\end{tabular}
\end{center}

\Needspace{10\baselineskip}
\paragraph{$a_{8,1}$}\label{cp:f8_1}

Interval \([11/40,7/20]\), degree
\(11\).
\[
 \gamma_{\mathrm{cen}}=\frac{2906631}{42949672960000000000}.
\]

\begin{center}\small
\setlength{\tabcolsep}{6pt}\renewcommand{\arraystretch}{1.05}
\begin{tabular}{@{}rr@{}}\toprule $j$ & Coefficient\\\midrule
0 & $4541652319736690104494432117189236$\\
1 & $58542095412901614051089877943357789$\\
2 & $324998487674504633799122463040031855$\\
3 & $1043623562119300718804790163867817985$\\
4 & $2175033647753308943294060595599804280$\\
5 & $3107982416109132729322798566997959168$\\
6 & $3119679107661850146060497926393538112$\\
7 & $2205900829338741313520258811596697120$\\
8 & $1079008328433631874590060444316480640$\\
9 & $348263840905556962673737005053861120$\\
10 & $66834324623165104928995451347184896$\\
11 & $5782804578420519649495410798343424$\\
\bottomrule\end{tabular}
\end{center}

\paragraph{$-b_2$}\label{cp:minus_b2}

Interval \([11/40,7/20]\), degree \(8\).
\[
 \gamma_{\mathrm{cen}}=\frac{2906631}{3276800}.
\]

\begin{center}\small
\setlength{\tabcolsep}{6pt}\renewcommand{\arraystretch}{1.05}
\begin{tabular}{@{}rr@{}}\toprule $j$ & Coefficient\\\midrule
0 & $1672561662620501073337$\\
1 & $12741060325508447015384$\\
2 & $42183356921716102259776$\\
3 & $79226723994626197901408$\\
4 & $92287208826992158203040$\\
5 & $68262770705332697932928$\\
6 & $31312520131502047008256$\\
7 & $8145377259155204862464$\\
8 & $920268808770087884032$\\
\bottomrule\end{tabular}
\end{center}

\paragraph{$a_{8,3}$}\label{cp:f8_3}

Interval \([11/40,7/20]\), degree
\(11\).
\[
 \gamma_{\mathrm{cen}}=\frac{26159679}{8589934592000000000}.
\]

\begin{center}\small
\setlength{\tabcolsep}{6pt}\renewcommand{\arraystretch}{1.05}
\begin{tabular}{@{}rr@{}}\toprule $j$ & Coefficient\\\midrule
0 & $76453881681257186140702664227308$\\
1 & $878339625779199015222882674766267$\\
2 & $4527210289619026898151559679511065$\\
3 & $13822979910765858554649734388862455$\\
4 & $27792289748991866788219487739120840$\\
5 & $38654847102809769780501597084962304$\\
6 & $37969317486653074200695611216537536$\\
7 & $26352431078101834973120074216923360$\\
8 & $12670423842026279168220039558145920$\\
9 & $4021051166581602833722404328815360$\\
10 & $758372304370740341963755016940288$\\
11 & $64418669243282334586384780796672$\\
\bottomrule\end{tabular}
\end{center}

\paragraph{$-a_{8,0}$}\label{cp:minus_f8_0}

Interval \([8/25,7/20]\), degree
\(11\).
\[
 \gamma_{\mathrm{cen}}=\frac{37786203}{65536000000000000000000000}.
\]

\begin{center}\small
\setlength{\tabcolsep}{6pt}\renewcommand{\arraystretch}{1.05}
\begin{tabular}{@{}rr@{}}\toprule $j$ & Coefficient\\\midrule
0 & $65661822358938930896095960635415789568$\\
1 & $804571677739873496735631026227008471040$\\
2 & $4374217500147998858068220717934712576000$\\
3 & $14011317553388042886277021194404283840000$\\
4 & $29496981985137035102998613378927608650000$\\
5 & $42972117132689255249106845640458871300000$\\
6 & $44294751408356389348410651763710070875000$\\
7 & $32353895204656603150677826207333495312500$\\
8 & $16429953515562670750978243685075525390625$\\
9 & $5529481341992931289183122583906279296875$\\
10 & $1110774471538655479395566374546123046875$\\
11 & $100957792634995608733873593767236328125$\\
\bottomrule\end{tabular}
\end{center}

\paragraph{$-b_0$}\label{cp:minus_b0}

Interval \([1/3,7/20]\), degree \(8\).
\[
 \gamma_{\mathrm{cen}}=\frac{753571}{373248000}.
\]

\begin{center}\small
\setlength{\tabcolsep}{6pt}\renewcommand{\arraystretch}{1.05}
\begin{tabular}{@{}rr@{}}\toprule $j$ & Coefficient\\\midrule
0 & $4054440550798400000000$\\
1 & $47187317716748160000000$\\
2 & $211456650455622264000000$\\
3 & $505582464566597997600000$\\
4 & $723756790646033070330000$\\
5 & $643830198503490269334000$\\
6 & $350334020657589415410600$\\
7 & $107156489574699559416780$\\
8 & $14154568045228414391481$\\
\bottomrule\end{tabular}
\end{center}

\paragraph{First axis comparison}\label{cp:axis7}

Interval \([11/40,23/80]\),
degree \(9\).
\[
 \gamma_{\mathrm{cen}}=\frac{13}{6871947673600000000000000}.
\]

\begin{center}\small
\setlength{\tabcolsep}{6pt}\renewcommand{\arraystretch}{1.05}
\begin{tabular}{@{}rr@{}}\toprule $j$ & Coefficient\\\midrule
0 & $104171434906209848373427818483456$\\
1 & $811856514243580062725001404465280$\\
2 & $2786747747758328247774381973235008$\\
3 & $5521362263424261451866470929644640$\\
4 & $6945030820414825216165879825923056$\\
5 & $5736531350258484717659330931219784$\\
6 & $3100455627846187245681976763643740$\\
7 & $1051944994871638082712854727926242$\\
8 & $201742951198783403932904480104390$\\
9 & $16453459326357952648537486682529$\\
\bottomrule\end{tabular}
\end{center}

\paragraph{Second axis comparison}\label{cp:axis8}

Interval \([23/80,8/25]\),
degree \(9\).
\[
 \gamma_{\mathrm{cen}}=\frac{169}{34359738368000000000000000000000}.
\]

\begin{center}\small
\setlength{\tabcolsep}{6pt}\renewcommand{\arraystretch}{1.05}
\begin{tabular}{@{}rr@{}}\toprule $j$ & Coefficient\\\midrule
0 & $18520536026711668187478698649792890625$\\
1 & $99315741694162172893577840879888562500$\\
2 & $211312524258677688776441238216518731250$\\
3 & $229316976177824657208994150290734050000$\\
4 & $169531163481473794976015657951082080000$\\
5 & $185123746200224150777390609398029337600$\\
6 & $241548605439332786707163665190788259840$\\
7 & $192921692778137500816663615840198590464$\\
8 & $78592036297944736512310241305861029888$\\
9 & $12826708630981104691462950972479766528$\\
\bottomrule\end{tabular}
\end{center}

\paragraph{Fixed correction
polynomials}\label{cp:fixed-correction-polynomials}

For \(C_7(x)\), the coefficient vector below is multiplied by
\[
\frac{1}{2621440000000}.
\]

\begin{center}\small
\setlength{\tabcolsep}{6pt}\renewcommand{\arraystretch}{1.05}
\begin{tabular}{@{}rr@{}}\toprule $j$ & Coefficient\\\midrule
0 & $697376047527$\\
1 & $583123712388$\\
2 & $870571114812$\\
3 & $163819248600$\\
4 & $-1893181923150$\\
5 & $6525531936$\\
6 & $717122545524$\\
7 & $-123042955464$\\
\bottomrule\end{tabular}
\end{center}

\Needspace{6\baselineskip}
For \(C_8(x)\), multiply the following coefficient vector by
\(1/167772160000000\):

\begin{center}\small
\setlength{\tabcolsep}{6pt}\renewcommand{\arraystretch}{1.05}
\begin{tabular}{@{}rr@{}}\toprule $j$ & Coefficient\\\midrule
0 & $41689951763431$\\
1 & $33304031870472$\\
2 & $51618041344884$\\
3 & $42437195877720$\\
4 & $-114580648922250$\\
5 & $-91237359423672$\\
6 & $107232450299436$\\
7 & $-1960330425768$\\
8 & $-7876838365110$\\
\bottomrule\end{tabular}
\end{center}

The nine Bernstein coefficients of \(C_8\), in order, are

\begin{center}\small
\setlength{\tabcolsep}{6pt}\renewcommand{\arraystretch}{1.05}
\begin{tabular}{@{}rr@{}}\toprule $j$ & Coefficient\\\midrule
0 & $\frac{41689951763431}{167772160000000}$\\
1 & $\frac{1146323893681}{4194304000000}$\\
2 & $\frac{12964865301913}{41943040000000}$\\
3 & $\frac{7558409401739}{20971520000000}$\\
4 & $\frac{7079733842209}{16777216000000}$\\
5 & $\frac{9838060814243}{20971520000000}$\\
6 & $\frac{3948897100181}{8388608000000}$\\
7 & $\frac{8890811888621}{20971520000000}$\\
8 & $\frac{60626494019143}{167772160000000}$\\
\bottomrule\end{tabular}
\end{center}

\subsection{Thirteen explicit positive expansions}
\label{octic:positive-expansions}
The following definitions and integer vectors establish the coefficient
signs \eqref{octic:coefficient-signs}, the two additional comparisons
for $10/29\le\sigma\le3/8$ used in the Hermite-majorant argument, and
the moment signs \eqref{octic:moment-signs}. Write $\Xi({\sigma})={\sigma}(2{\sigma}+1)>0$ and define
\begin{equation}\label{octic:polynomials}
 \begin{aligned}
 P_1&=N_0,&P_2&=N_1/\Xi({\sigma}),\\
 P_3&=N_2/\Xi({\sigma}),&P_4&=N_4/\Xi({\sigma}),\\
 P_5&=(32N_2N_4-9N_3^2)/\Xi({\sigma})^2,\\
 P_6&=N_0-(1-4{\sigma}^2-c){D_{\mathrm{prod}}},&P_7&=N_1-4{D_{\mathrm{prod}}}/5,\\
 P_8&=-E_3/\Xi({\sigma}),&P_9&=-E_4/\Xi({\sigma}),\\
 P_{10}&=K/\Xi({\sigma}),\\
 P_{11}&=(4E_0K-E_1^2)/\Xi({\sigma}),\\
 P_{12}&=(E_1+E_2/2+3E_3/16+E_4/16)/\Xi({\sigma}),\\
 P_{13}&=(E_1+2E_2+3E_3+4E_4)/\Xi({\sigma}).
 \end{aligned}
\end{equation}
These are polynomials in ${\sigma}$, of respective degrees
\begin{equation}\label{octic:degrees}
 (7,5,5,5,10,9,7,5,5,5,12,5,5).
\end{equation}
For $i\ne6,7$, set $a_i=10/29$, $b_i=2/5$; for $i=6,7$, set
$a_i=10/29$, $b_i=3/8$. Write $v$ for the nonnegative
variable in the fractional-linear change of coordinates, distinct from
all radial parameters. The vectors below, read row by row, give the
coefficients in increasing order of powers of $v$ in
\begin{equation}\label{octic:positive-identities}
 (1+v)^{\deg P_i}
 P_i\left(\frac{a_i+b_i v}{1+v}\right)
 =\gamma_i^{\mathrm{oct}}\sum_{j=0}^{\deg P_i}(\mathbf b_i)_j v^j.
\end{equation}
The exact positive rational scale factors are as follows.
\begin{center}\small
\begin{tabular}{@{}cr@{}}\toprule $i$&$\gamma_i^{\mathrm{oct}}$\\\midrule
1&$\frac{192}{1347646586640625}$\\
2&$\frac{256}{102555745}$\\
3&$\frac{128}{102555745}$\\
4&$\frac{2654208}{102555745}$\\
5&$\frac{1358954496}{10517680832505025}$\\
6&$\frac{147}{1901480637349101568}$\\
7&$\frac{6125}{282621973446656}$\\
8&$\frac{324}{2563893625}$\\
9&$\frac{1323378}{2563893625}$\\
10&$\frac{3}{20511149000}$\\
11&$\frac{78}{17276112461204543017578125}$\\
12&$\frac{3}{4102229800}$\\
13&$\frac{6}{2563893625}$\\
\bottomrule\end{tabular}
\end{center}

All displayed integers and scale factors are positive. Since
$v\ge0$ parametrizes $[a_i,b_i)$, and the leading
coefficient gives the positive endpoint value at $b_i$ after rescaling,
all thirteen signs follow, including endpoints. For a general polynomial $P_{\mathrm{aux}}(\sigma)=\sum_{h=0}^n p_h\sigma^h$, the coefficient
of $v^j$ on the left is
\[
 \sum_{h=0}^n p_h\sum_{k=0}^h
 \binom{h}{k}\binom{n-h}{j-k}a_i^{h-k}b_i^k,
\]
where a binomial coefficient with lower index outside its defining
range is zero. Applying this formula to \eqref{octic:polynomials},
using the numerators in Section~\ref{cp:rational-identities-and-monotonicities}
and the moments \eqref{octic:moments}, yields the exact vectors below.
Only $M_8$ contributes a term of degree four in $x$ to
\eqref{octic:residual}; its coefficient is $8169/409600$. Thus
$-E_4=(8169/409600)N_4$ and $P_9=(8169/409600)P_4$.
The intervals and degrees for these two polynomials coincide, and
$\gamma_9^{\mathrm{oct}}/\gamma_4^{\mathrm{oct}}=8169/409600$.
Consequently $\mathbf b_9=\mathbf b_4$ in
\eqref{octic:positive-identities}.

\begingroup\small

\[\mathbf b_{1}=\left[\begin{array}{rr}
13719381833849555390625 & 70112474831905952140625\\[1pt]
149675409958034502228125 & 172042623456105750378125\\[1pt]
114078545291696556571875 & 43137067838568391609075\\[1pt]
8466851370826497776095 & 648372892381286576431
\end{array}\right].\]
\[\mathbf b_{2}=\left[\begin{array}{rrr}
17685124972790625 & 63469278101019375 & 89389354243262750\\[1pt]
61015080370973250 & 19768602845565825 & 2334476813146831
\end{array}\right].\]
\[\mathbf b_{3}=\left[\begin{array}{rrr}
2324489036559375 & 47696845662350625 & 147747112071637250\\[1pt]
188134718389998750 & 110173343745942255 & 24554261504818753
\end{array}\right].\]
\[\mathbf b_{4}=\left[\begin{array}{rrr}
1348318059375 & 8499727250625 & 19965948253250\\[1pt]
22453756019550 & 12227088459855 & 2591399075809
\end{array}\right].\]
\[\mathbf b_{5}=\left[\begin{array}{r}
17977471882918260128349609375\\[1pt]
735891407817243454256238281250\\[1pt]
5423200747163370796562738671875\\[1pt]
19077506469237867983337316875000\\[1pt]
39783106057876545172187199218750\\[1pt]
53555091222649700816097576487500\\[1pt]
48071167727930491332583774698750\\[1pt]
28723236366996325693937128043000\\[1pt]
11008405204634115457621634679075\\[1pt]
2454742244931889412452424959650\\[1pt]
242611279221992907889853353727
\end{array}\right].\]
\[\mathbf b_{6}=\left[\begin{array}{rr}
3503042114057486783741952 & 41424305023342070943973376\\[1pt]
199867414729731460493484032 & 534584922353007458867302144\\[1pt]
888559280226029520017516096 & 961068343116160934334523424\\[1pt]
680525037828595298514439012 & 305428699003387230594222669\\[1pt]
79063948412045739629645598 & 9012453811569641181763827
\end{array}\right].\]
\[\mathbf b_{7}=\left[\begin{array}{rr}
365490049372611674112 & 2243827499497504899072\\[1pt]
5827801429303683915776 & 8277379512708428235520\\[1pt]
6914438684065878905520 & 3375100881893051888184\\[1pt]
881929639006732500039 & 93351661267039545717
\end{array}\right].\]
\[\mathbf b_{8}=\left[\begin{array}{rrr}
56294301474403125 & 252079716061606875 & 452527894076725750\\[1pt]
404301579481219050 & 178987018361006805 & 31347528496727099
\end{array}\right].\]
\[\mathbf b_{9}=\left[\begin{array}{rrr}
1348318059375 & 8499727250625 & 19965948253250\\[1pt]
22453756019550 & 12227088459855 & 2591399075809
\end{array}\right].\]
\[\mathbf b_{10}=\left[\begin{array}{rr}
123998694931675396875 & 582126434888218033125\\[1pt]
1085691238190734714250 & 1004285105512285209750\\[1pt]
460515167041115838315 & 83747972581400737669
\end{array}\right].\]
\[\mathbf b_{11}=\left[\begin{array}{r}
1993239889341987122024387219150807373046875\\[1pt]
20202208869688399985879411594697408398437500\\[1pt]
93939419925912853091013145075589755488281250\\[1pt]
264531017382489791664641043077322730757812500\\[1pt]
501142114986862368227048204311739630812890625\\[1pt]
670405164876305988846519954605199826424375000\\[1pt]
645944666585901820534737785303384700388687500\\[1pt]
448134462332480298755789924523543539197225000\\[1pt]
219448550394799821643652527941211455597438125\\[1pt]
72392614593516394060303131097549597486651500\\[1pt]
14600898261983946149459291110220447817161650\\[1pt]
1424226761190095415071992833971912398994180\\[1pt]
20804873393841416731905199094805627532311
\end{array}\right].\]
\[\mathbf b_{12}=\left[\begin{array}{rr}
7784171019750103125 & 27811526612529586875\\[1pt]
38589171293855616950 & 25674979465523962410\\[1pt]
8008782095826403317 & 891877618355888923
\end{array}\right].\]
\[\mathbf b_{13}=\left[\begin{array}{rr}
5504604400573753125 & 22879487178880336875\\[1pt]
37351773624251785750 & 29922916261799368650\\[1pt]
11750883542894315205 & 1807006034798572619
\end{array}\right].\]
\endgroup

The displayed coefficients and scale factors specify each polynomial identity explicitly.

\subsection{Expanded arithmetic for the small-radius comparisons}\label{sec:explicit-small-arithmetic}
In the small-radius bound, let \(z=\sigma-4/17\). The four columns below give the complete shifted monomial expansions of the polynomials \(M_{\mathrm{bot}},N_E,N_2,N_4\) in Section~\ref{sec:bottom-strip}. Each column is divided by its indicated positive denominator; the entry in row \(j\) is the numerator of the coefficient of \(z^j\). A zero beyond the degree of a polynomial denotes an absent term. Thus all the termwise estimates there follow directly from the displayed coefficients.
\[
\begin{array}{c|rr}
 j & M_{\mathrm{bot}} & N_E\\\hline
\text{denominator}&2015993900449&6975757441\\\hline
0&963573418694034720&50510844444752640\\
1&4196255274696971088&1834428821805127872\\
2&-18437644568110572000&5696038303373261664\\
3&46248830208553490352&-20947422336378165048\\
4&-204582172091717908116&-72461110773355912640\\
5&419645623889087495694&66955928448868035275\\
6&-208652656203941830056&363060345990165954576\\
7&206154493211835067840&-37827851344877082788\\
8&-42486922617201913824&17267749900700658896\\
9&10733827002535002912&0\\
10&-612604098517238528&0\\
\end{array}\]
\[
\begin{array}{c|rr}
 j & N_2 & N_4\\\hline
\text{denominator}&6975757441&6975757441\\\hline
0&43011475269120&157719011328000\\
1&-706855058193024&1329669811814400\\
2&-408606581472768&-85483392307200\\
3&14815804925226816&-21821707453209600\\
4&-10537891568937120&-20152139312384000\\
5&-53820612431049000&127801646617968000\\
6&-4037751079972992&173671855820851200\\
7&-11464915046970144&-257364415705600\\
8&393767556029568&7500334400563200\\
\end{array}\]
\subsection{Endpoint arithmetic in the central product estimate}\label{sec:explicit-central-arithmetic}
The entries in the next table list the exact differences between the upper or lower expressions in the five scalar-comparison rows of Section~\ref{cp:scalar-majorization-by-a-cubic-taylor-estimate} and the stated rational bounds. More precisely, for row \(i\) with endpoints \(a,b\), exponent \(k\), and listed constants \(n,A,B,C,q\), the differences, in order, are
\begin{align*}
 \Delta_{i,1}&=n_k(a)-n,\\
 \Delta_{i,2}&=\widehat h_2(a)+\binom{k}{2}\delta(b)-K_{b,q^2}-A,\\
 \Delta_{i,3}&=\widehat h_3(b)-\binom{k}{3}\delta(a)-B,\\
 \Delta_{i,4}&=\min\{0,\widehat h_4(a)\}-C,\\
 \Delta_{i,5}&=c-f_a(-q).
\end{align*}
Each listed quotient is an exact fraction in lowest terms, with positive denominator. In particular, all five comparisons in every row are strict except where a zero is explicitly displayed.
\begingroup\small\renewcommand{\arraystretch}{2.2}\begin{longtable}{@{}ccr@{}}\toprule Row & Difference & Exact value \\ \midrule\endhead\bottomrule\endfoot
1&$\Delta_{1,1}$&$\dfrac{789512318548789}{858740309216608000}$\\
1&$\Delta_{1,2}$&$\dfrac{341090312022355467328463296712731460310239219}{28531129554203163394502831916928017192602580000}$\\
1&$\Delta_{1,3}$&$\dfrac{1552072972949123542280557}{933791692744076242808732000}$\\
1&$\Delta_{1,4}$&$\dfrac{25148069252427}{2683563466301900}$\\
1&$\Delta_{1,5}$&$\dfrac{1267487}{2490010000}$\\
2&$\Delta_{2,1}$&$\dfrac{739931892312527687}{900473510398962500000}$\\
2&$\Delta_{2,2}$&$\dfrac{14317428245128054411769174172405015133020162413}{1838027415125762146348346839567113020162818120000}$\\
2&$\Delta_{2,3}$&$\dfrac{66077597825162434681966483}{481294812707200161218512964000}$\\
2&$\Delta_{2,4}$&$\dfrac{4340418005991}{720378808319170}$\\
2&$\Delta_{2,5}$&$\dfrac{22007}{9610000}$\\
3&$\Delta_{3,1}$&$\dfrac{4605130229104933089}{5058420991687831040000}$\\
3&$\Delta_{3,2}$&$\dfrac{22399720906598568551044893245221}{6330590709092295658644521669780000}$\\
3&$\Delta_{3,3}$&$\dfrac{1228399827973252553}{268365686381397605000}$\\
3&$\Delta_{3,4}$&$\dfrac{42668574890411}{7903782799512236}$\\
3&$\Delta_{3,5}$&$\dfrac{105147}{16810000}$\\
4&$\Delta_{4,1}$&$\dfrac{12971483117}{5731445720000}$\\
4&$\Delta_{4,2}$&$\dfrac{5589955310389232959751}{1150230847683042165040000}$\\
4&$\Delta_{4,3}$&$\dfrac{252150572670940023}{49733260449802930000}$\\
4&$\Delta_{4,4}$&$\dfrac{67153301}{7164307150}$\\
4&$\Delta_{4,5}$&$\dfrac{383}{90000}$\\
5&$\Delta_{5,1}$&$\dfrac{9612564601796021}{5838333542660937500}$\\
5&$\Delta_{5,2}$&$\dfrac{12101711102589271}{2746352098467705000}$\\
5&$\Delta_{5,3}$&$\dfrac{180036721062181}{84072003014317500}$\\
5&$\Delta_{5,4}$&$0$\\
5&$\Delta_{5,5}$&$\dfrac{83567}{8410000}$\\
\end{longtable}\endgroup
The remaining endpoint comparisons used in the moment estimates are given below, with the difference oriented so that positivity is the required assertion. In the first interval, the ``upper bound'' and ``lower bound'' mean the exact right-hand expressions in Section~\ref{cp:the-moment-comparison-in-four-intervals}.
\begingroup\small\renewcommand{\arraystretch}{2.2}\setlength{\tabcolsep}{3pt}\begin{longtable}{@{}p{.43\linewidth}r@{}}\toprule Difference & Exact value\\\midrule\endhead\bottomrule\endfoot
$c_{7,1}-c_{8,1}$&$\dfrac{100397143059}{4194304000000}$\\
$\begin{aligned}&-1/20-a_{8,1}(23/80)\\ &-\delta(11/40)(c_{7,1}-c_{8,1})\end{aligned}$&$\dfrac{49431448435457964711095662444692166871}{59027256711139147431480708300800000000000}$\\
$-1/40$ minus the upper bound for the sum in the first interval&$\dfrac{4675091816178323459508487804509355049}{1383451329167323767925329100800000000000}$\\
$\text{lower bound for }a_{7,3}$&$\dfrac{327742125371582155397379}{3784388483585335061945600}$\\
$-\widehat h_4(23/80)$&$\dfrac{483319281246912}{1975945699878059}$\\
$-a_{8,4}(8/25)$&$\dfrac{74272853704024548410561601}{7164414751105299251200000000}$\\
$\begin{aligned}&-a_{8,2}(23/80)-a_{8,3}(8/25)\\ &-a_{8,6}(23/80)\end{aligned}$&$\dfrac{361515421523098697736586283692726023}{34042460400905586316547194880000000000}$\\
$a_{8,3}(23/80)$&$\dfrac{37953449326747646253832277142237}{331508678111253695037440000000000}$\\
$\frac{1099}{2000}-a_{8,0}(\frac{8}{25})$&$\dfrac{969263901110280578118834457}{2893321341792524697600000000000}$\\
$\frac{1}{30}-a_{8,1}(\frac{1}{3})$&$\dfrac{2435808897969270923887}{2844635273625600000000000}$\\
$- \frac{171}{500}-a_{8,2}(\frac{8}{25})$&$\dfrac{553667349325797382388466631}{9403294360825705267200000000000}$\\
$\frac{81}{500}-a_{8,3}(\frac{1}{3})$&$\dfrac{2700690286510323869}{27091764510720000000000}$\\
$\frac{1}{125}-a_{8,4}(\frac{1}{3})$&$\dfrac{527100399250773013}{4334682321715200000000}$\\
$\frac{13}{1000}-a_{8,6}(\frac{8}{25})$&$\dfrac{31757626915610723617013869}{447775921944081203200000000000}$\\
$\frac{109}{200}-b_0(\frac{1}{3})$&$\dfrac{16808625359}{39065120640000}$\\
$\frac{77}{1000}-b_1(\frac{10}{29})$&$\dfrac{5107752361651}{6216583533888000}$\\
$- \frac{43}{100}-b_2(\frac{1}{3})$&$\dfrac{153140593}{65108534400}$\\
$\frac{187}{1000}-b_3(\frac{10}{29})$&$\dfrac{741153179651}{1036097255648000}$\\
$\frac{19}{1000}-b_4(\frac{10}{29})$&$\dfrac{229714786807}{296027787328000}$\\
$b_3(1/3)$&$\dfrac{44397953}{275591680}$\\
$b_4(1/3)$&$\dfrac{5532747}{551183360}$\\
$n_8(1/3)-2/5$&$\dfrac{7003483349}{305196255000}$\\
$\delta(10/29)+1/80$&$\dfrac{178988120106737}{544598619999980000}$\\
$\widehat h_2(1/3)+4/25$&$\dfrac{115500689}{24415700400}$\\
$\widehat h_3(10/29)$&$\dfrac{8557325300}{32378039239}$\\
$\widehat h_4(1/3)$&$\dfrac{379280}{753571}$\\
\end{longtable}\endgroup

\section*{Declarations}

\subsection*{Conflict of interest}

The author declares no conflict of interest.

\subsection*{Data availability}

No empirical datasets were used. The exact algebraic data and rational
bounds used in the proof are included in the article.

\subsection*{Computational aspects}

Computer algebra and exact rational arithmetic were used in the
generation and verification of the explicit algebraic identities,
the positivity certificates and the rational evaluations, including
the finite-series bounds for transcendental constants.

The mathematical arguments, the explicit certificates, and the
coefficient formula~\eqref{eq:positive-expansion-rule}, from which
each certificate can be recomputed from the polynomial it certifies,
are included in the article. All finite algebraic checks are based on
exact identities and rational inequalities. The bounds for
transcendental quantities are obtained through elementary analytic
inequalities and finite series with rigorous remainder estimates.
No step of the proof relies on numerical sampling, numerical
quadrature, floating-point arithmetic, or interval arithmetic.

\subsection*{Declaration concerning AI-assisted work}

The author used AI tools for assistance with English-language editing, editorial organization, and LaTeX formatting. The mathematical results, the proof strategy, and the analytical arguments presented in the body of the document are the work of the author, who assumes full responsibility for the mathematical content of the article.

\end{document}